\documentclass[11pt,reqno]{amsart}

\usepackage{amsfonts}
\usepackage{eurosym}
\usepackage{amssymb}
\usepackage{amsthm}
\usepackage{amsmath}
\usepackage{bm}
\usepackage{cite}
\usepackage{mathrsfs}
\usepackage{xcolor}
\usepackage[OT1]{fontenc}
\usepackage[left=2cm, right=2cm, top=3cm]{geometry}
\usepackage{hyperref}
\hypersetup{colorlinks=true, linkcolor=blue, citecolor=red}

\usepackage{times}

\usepackage[shortlabels]{enumitem}
\usepackage{enumitem}
\newcommand{\subscript}[2]{$#1 _ #2$}
\numberwithin{equation}{section}
\usepackage[square,numbers,sectionbib]{natbib}
\usepackage{varioref}
\usepackage{listings}
\usepackage{color}
\definecolor{dkgreen}{rgb}{0,0.6,0}
\definecolor{gray}{rgb}{0.5,0.5,0.5}
\definecolor{mauve}{rgb}{0.58,0,0.82}

\allowdisplaybreaks
\newtheorem{theorem}{Theorem}[section]

\newtheorem{lemma}[theorem]{Lemma}
\newtheorem{proposition}[theorem]{Proposition}
\theoremstyle{remark}
\newtheorem{remark}[theorem]{Remark}

\def\d{{\rm d}}

\def\V{{\mathbf{V}}}
\def\H{\mathbf{H}}
\def\W{\mathbf{W}}
\def\A{\mathbf{A}}
\def\L{\mathbf{L}}

\def\g{\textbf{\textit{g}}}
\def\ddt{\frac{\d}{\d t}}
\def\C{\mathcal{C}}

\def\n{\mathbf{n}}
\def\E{\mathcal E}

\newcommand{\numberset}{\mathbb}

\newcommand{\R}{\numberset{R}}
\def\nablag{\nabla_\Gamma}
\def\divg{\mathrm{div}_{\Gamma}}
\def\J{\mathbf{J}}
\def\cd{\cdot}

\def\vn{v_\n}
\def\intg{\int_{\Gamma(t)}}
\def\ints#1{\int_{#1(t)}}
\def\Vn{v_\n}
\def\VVn{\V_\n}
\def\dtb{\partial^\bullet_t}
\def\dt{\partial^\circ_t}
\def\P{\mathbf{P}}
\def\v{\mathbf{v}}
\def\ig{\intg}
\def\u{\mathbf{u}}
\def\vphi{\varphi}

\def\D{\mathbf{D}}
\def\D{\mathbf{D}^-}
\def\tphi{\widetilde{\phi}}

\def\A{\mathbf{A}}
\def\phmt{\phi_{-t}}
\def\phit_{\phi_t}
\def\dtast{\partial^\ast_t}
\def\tphimt{\tphi_{-t}}
\def\nablagz{\nabla_{\Gamma_0}}
\def\D{\mathbf{D}}
\def\Dm{\D^-}
\def\tphit{\tphi_t}

\def\tphimt{\tphi_{-t}}
\def\tvphi{\widetilde{\vphi}}
\def\a{\mathbf{a}}
\def\b{\mathbf{b}}
\def\c{\mathbf{c}}
\def\di{\mathbf{d}}
\def\nablagphi{\nabla_{\Gamma(t)}^\phi}
\def\nablagphis{\nabla_{\Gamma(s)}^\phi}
\def\tT{\widetilde{T}}
\def\H{\mathbf{H}}
\def\Zut{Z^1_{\tT}}
\def\Zdt{Z^2_{\tT}}
\def\Ng{\mathcal{N}_\Gamma}
\def\w{\mathbf{w}}
\def\XT{X_{\widetilde{T}}}
\def\YT{Y_{\widetilde{T}}}
\def\norm#1{\left\Vert#1\right\Vert }
\def\norma#1{\left\vert#1\right\vert }
\def\trho{\rho}
\def\f{\mathbf f}
\def\dtn{\partial_t}
\def\ddt{\frac{d}{dt}}
\def \Pt{P(t)}
\def\LQs{L^2(0,\widetilde{T};\L^2_\sigma(\Gamma_0))}
\def\LQ{L^2(0,\widetilde{T};\L^2(\Gamma_0))}
\def\CT{C(T)}
\def\thetaq{\theta(1-\frac1q)}
\def\gz{\Gamma_0}
\def\LLQ(#1,#2){L^{#1}(0,\widetilde{T};\L^{#2}(\Gamma_0))}
\def\EE{\mathcal{E}^{\gz}}
\def\div#1{\text{div}_{#1}}
\def\divv{\text{div}}
\def\dts{\partial_t^*}
\def\Lqp{L^q(0,\tT;L^p(\gz))}
\def\Lqpb{L^q(0,\tT;\L^p(\gz))}
\def\gam{\Gamma}
\def\gt{(\gam(t))}

\def\non{\nonumber}
\def\Phitn{\Phi_t^n}
\def\z{\mathbf z}
\def\X{\mathbf X}
\def \ds{d \sigma}
\def\A{\mathbf{A}}

\def\rev#1{{\color{black}#1}}
\everymath{\displaystyle}

\def \no#1#2#3 {{\bf #1} (#3), #2.}
\def \eds#1#2#3 {#1, #2, #3.}

\author[H.~Abels]{Helmut Abels}\author[H.~Garcke]{Harald Garcke}
\address[H.~Abels, H.~Garcke]{Fakult\"{a}t f\"{u}r Mathematik\\
Universit\"{a}t Regensburg\\
93040 Regensburg, Germany}
\email[H.~Abels]{helmut.abels@ur.de}
\email[H.~Garcke]{harald.garcke@ur.de}
\author[A.~Poiatti]{Andrea Poiatti}
\address[A.~Poiatti]{Faculty of Mathematics, University of Vienna, 1090 Vienna, Austria}
\email[A.~Poiatti]{andrea.poiatti@univie.ac.at}

\begin{document}

\title
[Diffuse Interface Model for Two-Phase Flows on Evolving Surfaces]
{Diffuse Interface Model for Two-Phase Flows on Evolving Surfaces with Different Densities: Local Well-Posedness}

\date{\today }


\maketitle
\begin{abstract}
A Cahn--Hilliard--Navier--Stokes system for two-phase flow on an evolving surface with non-matched densities is derived using 
methods from rational thermodynamics. For a Cahn--Hilliard energy with {\color{black}a singular (logarithmic)} potential short time well-posedness of strong solutions together
with a separation property is shown, {\color{black}under the assumption of \textit{a priori }prescribed surface evolution.} The problem is reformulated with the help of a pullback to the initial surface. Then 
a suitable linearization and a contraction mapping argument for
the pullback system are used.  {\color{black}In order to deal with the linearized
system}, it is necessary to show maximal $L^2$-regularity for the {\color{black}surface }Stokes {\color{black}operator in the case of variable viscosity} and to obtain maximal $L^p$-regularity for the linearized Cahn-Hilliard system. 
\end{abstract}
\section{Introduction}
Diffuse interface models for two-phase flows in the Euclidean space are classical topics in modeling, analysis and numerics. Several models have been derived and analysed, starting with the famous model H introduced  by Hohenberg and Halperin \cite{HH77}.
The model H is valid for a situation in which the two fluid phases have the same density. 
Extensions to the case of non-matched densities have been introduced by Lowengrub and Truskinovsky \cite{LowengrubQuasiIncompressible}, by Abels, Garcke and Grün \cite{AGG}, {\color{black}and others, see ten Eikelder et al.~\cite{EikelderUnifiedNSCH} for an overview and further references}. The goal of this work is to generalize the model introduced by Abels, Garcke and Grün  \cite{AGG} to the case of a two-phase flow on an evolving fluidic surface and to show local existence of strong solutions {\color{black} in the case of \textit{a priori }prescribed surface evolution}. 
In the Euclidean case, the model from \cite{AGG} is given as
\begin{align*}
     \begin{cases}
     \rho \partial_t \u +((\rho\u+\J_\rho)\cdot{\nabla})\u-2\divv(\nu(\vphi)\E(\u))+\nabla \pi=-\divv(\nabla\vphi\otimes\nabla\vphi),\\
         \divv\u=0,
         \\\partial_t \vphi-\divv(m(\vphi)\nabla\mu)+\divv(\vphi\u)=0,\\
          \mu=-\Delta\vphi+\Psi'(\vphi),
     \end{cases}
 \end{align*}
where $\J_\rho=-\frac{\widetilde{\rho}_1-\widetilde{\rho}_2}{2}\nabla\mu$. The unknowns are the fluid velocity $\u$,
the pressure $\pi$, the phase field variable $\vphi$ and the chemical potential $\mu$. In addition, $\E(\u)$ is the
rate-of-strain tensor, \rev{which is} 
 the symmetrized velocity gradient, $\widetilde{\rho}_1,\widetilde{\rho}_2 $ are the  mass densities of the two fluids, $\rho=\rho(\vphi)$
affine linearly interpolates between the two densities, $\nu$ is the viscosity of the mixture, $m$ the mobility and $\Psi$ is a double well potential. For further details, we refer to \cite{AGG}.

  In order to generalize the model above to an evolving surface, it is first of all necessary to
  generalize the fluid model. Navier--Stokes equations defined on stationary manifolds is a classic topic in analysis,  see \cite{ArnoldK, EbinM, MitreaTaylor, Temam, Simonett}. However, only recently there was a growing interest in studying Navier--Stokes equations on evolving surfaces.
     Derivations of evolving surface Navier--Stokes equations have been given in \cite{HuZhangE, V1, KobaLiuGiga, Miura2018, NitschkeRV}.
     The first two of these papers derive the equations from mass and momentum balance laws on evolving surfaces and the two latter 
     use balance laws in bulk regions and derive the Navier--Stokes system as a thin film limit. In \cite{KobaLiuGiga} an 
     energetic variational
approach is used to derive the equations.
For a comparison of different approaches to derive the Navier-Stokes equations on an evolving surface we refer to
    \cite{BrandnerReuskenSchwering}. We  mention \cite{WangZZ}, who show  local well-posedness of a Navier-Stokes system for an evolving hypersurface taking an elastic regularisation into account. {\color{black}Also, in the recent contribution~\cite{V2}, the authors show existence and uniqueness of weak solutions to a Navier--Stokes system on \textit{a priori }prescribed evolving surfaces.}

  {\color{black}On the other hand, concerning only the advective Cahn-Hilliard system with singular potential in the case of \textit{a priori} prescribed evolving surface setting, the first analytical results about weak well-posedness of solutions can be traced back to Caetano and Elliott \cite{DE}, whereas well-posedness of global strong solutions, together with the first result on the strict separation property on evolving surfaces, has been shown in \cite{CEGP}.}

    Only very recently coupled Cahn--Hilliard--Navier--Stokes systems have  appeared  in the literature.
    Such systems have been first studied by \cite{NitschkeVoigtWensch} using a stream function formulation and by
    \cite{BGNfluidmembrane} in which a Cahn--Hilliard--Navier--Stokes system on an evolving surface is coupled to bulk Navier--Stokes equations.

    {\color{black}Combining some techniques developed in \cite{DE,V2}}, the first well-posedness result {\color{black} of weak solutions} for a diffuse interface model for two-phase flow on a given evolving surface has  recently been shown by Elliott and Sales \cite{ES}, which is also the first analytical result for such {\color{black}coupled systems} at all. In this paper we generalize the model studied in \cite{ES} to the case of non-matched densities and we will hence obtain a generalization 
    of the model introduced in \cite{AGG} to an evolving surface.
    Our derivation will be based on mass and momentum balance laws on evolving surfaces and we will restrict possible constitutive
    relations by exploiting a local formulation of the free energy inequality.
    The system we derive in the case where we do not prescribe the normal velocity is given by
    \begin{align}
     \begin{cases}
     \rho \dt \u +((\rho\u+\J_\rho)\cdot\widehat{\nabla}_\Gamma)\u-2\divg(\nu(\vphi)\E_S(\u))+\nablag \pi+H\pi\n=-\divg(\nablag\vphi\otimes\nablag\vphi),\\
         \divg\u=0,
         \\\dt\vphi-\divg(m(\vphi)\nablag\mu)+\divg(\vphi\u)=0,\\
          \mu=-\Delta_\Gamma\vphi+\Psi'(\vphi),
     \end{cases}
     \label{unproj}
 \end{align}
 where $H$ is the mean curvature, i.e., the sum of the principal curvatures,  $\dt$ is the normal material time derivative, $\E_S(\u)$ is the surface rate of deformation tensor
 and $\widehat{\nabla}_\Gamma$, $\divg$, $\nablag$ and $\Delta_\Gamma$ are variants of the flat differential operators on a surface.
 We will define these operators precisely in the following section.
 The system above still has the normal part of the velocity $\u$ as an unknown. Similar as in 
 \cite{BrandnerReuskenSchwering, ES} we will now prescribe the normal velocity of the surface, which fixes one degree of freedom 
 at each point, and project the first equation in \eqref{unproj}
 to the tangent space, which lead to one equation less in each point. We then obtain  the system 
  \begin{align}
     \begin{cases}
     \rho \P\dt \v +((\rho\v+\J_\rho)\cdot\nabla_\Gamma)\v+\rho \Vn\mathbf{H}\v+\Vn\H\J_\rho-2\P\divg(\nu(\vphi)\E_S(\v))+\nablag \pi\\
     =-\P\divg(\nablag\vphi\otimes\nablag\vphi)+2\P\divg(\nu(\vphi)\Vn\mathbf{H})+\frac\rho2\nablag(\Vn)^2,\\
         \divg\v=-H\Vn,
         \\\dt\vphi-\divg(m(\vphi)\nablag\mu)+\nablag\vphi\cdot \v=0,\\
         \mu=-\Delta_\Gamma\vphi+\Psi'(\vphi),
     \end{cases}
     \label{AGPproj}
 \end{align}
 where $\P$ is the projection to the tangent space, $\Vn$ is the normal velocity and  $\mathbf{H}$ is the  Weingarten map. \rev{Notice that additional tangential forcing terms can be considered in model \eqref{AGPproj}. Due to the \textit{a priori} prescription of the normal velocity of the surface evolution, here we admittedly cannot account for normal forcing terms, so that some effects, like the ones due to gravitational forces, cannot be completely described by this model.}

    We refer to 
    \cite{PZQuainiO} for a numerical study of a Cahn--Hilliard--Navier--Stokes system which generalizes the model of Abels, Garcke and
    Grün \cite{AGG} to a stationary surface. In addition, we mention that related are also the models discussed in \cite{BachiniKV,BachiniKNVoigt}
    who study numerical approaches for Cahn--Hilliard--Navier--Stokes systems with matched densities on evolving surfaces.

    The main result of this paper is Theorem \ref{strong1a} where the short time well-posedness of strong solutions together with a separation property is stated. Main ideas used in the paper are a suitable surface Piola transform which yields that the divergence-free property of tangential vector fields is preserved, together with the evolving Banach and Hilbert space framework introduced in \cite{Def, AE}. We {\color{black}obtain}  equivalent pullback equations which are formulated on the initial surface. The concrete functional setting for the pullback problem then  uses ideas introduced in \cite{AWe}, where the existence of strong solutions in the related case of a {\color{black}stationary} domain in Euclidean space is shown. 
    In order to prove existence of strong solutions locally in time we use  a suitable linearization and a contraction mapping argument for the pullback system. In the proof it is necessary to show maximal $L^2$-regularity for the Stokes part of the linearized system and to obtain maximal $L^p$-regularity for the linearized Cahn-Hilliard system. Compared to \cite{AWe}
    several new difficulties arise due to the additional geometric terms. {\color{black}In particular, a nontrivial issue is to establish a suitable regularity theory for the surface Stokes operator in the case of space dependent viscosities}.
    We {\color{black}} also refer to the forthcoming paper \cite{AGP3},
where  the existence and uniqueness of a global strong
solution to  \eqref{AGPproj} will be shown using heavily the results of the present paper, {\color{black}together with some ideas from \cite{CEGP}}.

    The outline of the paper is as follows. In the following section we derive the governing equations. In Section \ref{sec:Piola}
    we introduce the surface Piola transform and give the precise functional setting. We then formulate the main  local existence result
in Section \ref{mainresult}. To proceed, we formulate an equivalent system  using a divergence free velocity and state the pullback equations.    The main  technical part is Section \ref{sec:shortpullback} where short time existence of a regular solution to the pullback problem is shown. It is then possible to prove the main result by lifting the pullback solution to the evolving surface.
In two appendices regularity for Laplace equations on evolving surfaces and details on how to obtain the pullback equations are stated.
    \section{Derivation of the model}
Here we derive the Cahn--Hilliard--Navier--Stokes system on evolving surfaces using the framework of rational thermodynamics. The presentation adapts the approach in \cite{AGG}. In what follows, we consider smooth compact  orientable hypersurfaces $\Gamma(t)$ embedded in $\R^{n+1}$ evolving in time (in a sufficiently smooth way). For such an evolving hypersurface we introduce  a vector field  $\VVn$ on $\Gamma(t)$, such that $v_n=\VVn\cdot\n$ is the normal velocity of $\Gamma(t)$, where $\n$ is a normal vector, and $\VVn=\vn\n$. 
Clearly the geometric evolution of the surface is only in the normal direction as tangential directions of the velocity do not change the shape of the surface. 

We denote by $\dtb$ the total material derivative with respect to an advective velocity $\hat{\V}$, which has normal component $\VVn$. \rev{In particular, this is defined as 
\begin{align}
\partial^\bullet_t f =
\frac{\partial f_e}{\partial t} + \hat \V\cdot\nabla f_e\quad \text{ on }\bigcup_{t\in[0,T]}\{t\}\times \Gamma(t),
\end{align}
where $f_e$
is a smooth extension of $f : \bigcup_{t\in[0,T]}\{t\}\times \Gamma(t)\to \R$ into a spatial neighborhood of $\bigcup_{t\in[0,T]}\{t\}\times \Gamma(t)$.} Moreover, by $\dt$ we mean the normal material time derivative \rev{(i.e., with respect to $\hat\V=\VVn$)}. We are making a slight abuse of notation in that $\dtb$ can represent various total material time derivatives, according to the transport velocity $\hat{\V}$ we are considering, whereas $\dt$ is always the same, since the normal component of $\hat{\V}$ is the same. In any case, we have the decomposition
\begin{align}
\dtb f=\dt f+\nablag f\cdot (\hat{\V}-\V_\n), 
\label{scalar}
\end{align}
whereas, for a vector field $\mathbf{f}$, 
$$
\P\dtb \mathbf{f}=\P\dt \mathbf{f}+(\nablag \mathbf{f})(\hat{\V}-\V_\n) , 
$$
where $\P$ is the projector on the tangent space to $\Gamma$, namely $\P:=\mathbf{I}-\n\otimes \n$. This means that we will always make clear what $\dtb$ means with respect to $\hat{\V}$. Here notice that $\nabla_\Gamma$, a column tangential vector if applied to scalar quantities, \rev{when applied to a vector field is defined as the generalization of the classical covariant derivative to non-tangential vector fields. Indeed, it is defined as $\nabla_\Gamma \mathbf f:=\P \nabla \mathbf f_e \P$}, where $\nabla$ is intended to be the gradient of the function \rev{$\mathbf f_e$} defined from $\mathbf f$ by extension in an $n+1$ -neighborhood $U$ of $\Gamma$. In particular, we have, given a vector function $\mathbf f_e$ extended to $U$, $$
\nabla \mathbf f_e=\begin{bmatrix}
\nabla^T f_{e,1}\\
\ldots\\
\nabla^T f_{e,n+1}
\end{bmatrix}.
$$
It can be shown that the covariant gradient
only depends on the values of the function on $\Gamma\cap U$. We also introduce $$\widehat{\nabla}_\Gamma\f:=\begin{bmatrix}\nablag^T f_1\\\ldots\\\nablag^T f_{n+1}
\end{bmatrix},
$$
so that $\nablag \f=\P \widehat{\nabla}_\Gamma\f$.

For a vector $\v$, the tangential divergence is then defined as $\divg \v:=\mathrm{tr}(\nablag \v )$. On the other hand, in case of a 3-by-3 matrix $\A$, it is defined as $\divg \A:=(\divg(\mathbf e_i^T\A))_{i=1,2,3}$, with $(\mathbf e_i)_i$ as the canonical Euclidean basis in $\R^3$. For the basic properties of these operators we refer, for instance, to \cite{V1,V2}.
 
Assume now to have two fluids which are mixed up together and flowing on an evolving surface $\Gamma$. Their normal velocity coincides with $\VVn$. Thanks to the inextensibility constraint, we can study a mass balance equation for each species as follows. Consider an element of surface $\Sigma(t)\subset \Gamma(t)$ which follows the purely geometric evolution of $\Gamma$, i.e., it evolves with velocity $\VVn$. Then the mass balance for the two fluids $i=1,2$ can be assumed to be written as
$$
\frac{d}{dt}\ints{\Sigma} \rho_i\ds=\rev{-\ints{\partial\Sigma}\widehat{\mathbf J}_i\cdot \boldsymbol\tau\d\mathcal H^{n-1}}=-\ints{\Sigma} \divg(\hat{\J}_i)\ds,
$$
\rev{where $\boldsymbol\tau$ denotes the outer unit conormal to $\partial\Sigma(t)$, and in the last identity we used the classical divergence theorem on smooth manifolds (see, e.g., \cite{BGNhandbook}).}
Here $\hat{\J}_i$ are the tangential mass fluxes of the fluids and $\rho_i$ are the densities of the fluids. By the transport theorem for the evolving surface $\Sigma(t)$ \rev{(see, e.g., \cite[Theorem 32]{BGNhandbook})}, \rev{it holds 
$$
\frac{d}{dt}\ints{\Sigma} \rho_i \ds=\ints{\Sigma}\left(\dt\rho_i+\rho_i\divg\VVn  \right)\ds,
$$
and thus }we obtain
$$
\frac{d}{dt}\ints{\Sigma} \rho_i \ds=\ints{\Sigma}\left(\dt\rho_i+\rho_i\divg\VVn  \right)\ds=\ints{\Sigma}\left(\dt \rho_i+\divg (\rho_i\VVn)-\nablag\rho_i\cdot\VVn  \right)\ds
$$
and, together with the balance equation for the masses, exploiting the arbitrariness of $\Sigma(t)$, and recalling that $\nablag\rho_i\cdot\VVn= v_\n \nablag\rho_i\cdot \n=0$, we infer the local mass balance
$$
\dt\rho_i+\divg (\rho_i\VVn)+\divg(\hat{\J}_i)=0,\quad i=1,2.
$$
Introducing the tangential velocities $\v_i:=\frac{\hat{\mathbf{J}}_i}{\rho_i}$ we easily end up with 
\begin{align}
\dt\rho_i+\divg(\rho_i(\v_i+\VVn))=0,\quad i=1,2.
\label{mass2}
\end{align}
Furthermore, introducing the tangential mass flux relative to $\v$ (which will be defined later on as to be the tangential volume averaged velocity), i.e.,
$$
\J_i:=\hat{\J}_i-\rho_i\v,
$$
the same mass balance can be rewritten as
\begin{align}
\dt\rho_i+\divg(\J_i)+\divg(\rho_i(\v+\VVn))=0,\quad i=1,2.
\label{rho2}
\end{align}
We now define the total mass density as $\rho:=\rho_1+\rho_2$, to get, summing up the equations in \eqref{rho2},
\begin{align}
\label{ma}
\dt\rho+\divg(\J_\rho)+\divg(\rho\u)=0,
\end{align}
where $\J_\rho:=\J_1+\J_2$ is the total flux and $\u=\v+{\VVn}$ is the complete (volume averaged) velocity of the fluids (tangential component + normal component).

We now define the volume fractions $\varphi_i$ of the two fluids, i.e., $\varphi_i:=\dfrac{\rho_i}{\widetilde{\rho}_i}$, having defined the
specific (constant) density of the unmixed fluids by $\widetilde{\rho}_i$. Now, by dividing \eqref{rho2} by $\widetilde{\rho}_i$ we end up with 
\begin{align}
\dt\varphi_i+\divg(\J_{\varphi_i})+\divg(\varphi_i\u)=0
,\label{conc}
\end{align}
where we defined $\J_{\varphi_i}:=\frac{\J_i}{\widetilde{\rho}_i}$. We also need  to fix the tangential velocity $\v$. In particular, we choose it as to be the volume averaged tangential velocity, i.e.
$$
\v:=\v_1\varphi_1+\v_2\varphi_2.
$$
\rev{We now assume that the excess volume due to mixing
is zero (see \cite{AGG}), which results in the assumption}
\begin{align}
\varphi_1+\varphi_2=1.
\label{sum}
\end{align}
\rev{As a consequence, multiplying \eqref{mass2} by $\frac1 {\widetilde{\rho}_i}$ and adding the two equations then gives
$$
\divg \u=\divg\left(\v_1\varphi_1+\v_2\varphi_2+\VVn\right)=-\dt (\varphi_1+\varphi_2)=-\dt1=0.
$$
From now on we hence require
$$
\divg \u=0,
$$
recalling $\u=\v+\VVn$.
Notice that this pointwise constraint is much stronger than the inextensibility (integral)
constraint on $\Gamma(t)$, since it also requires information on the tangential component v. In fact,  it is an incompressibility
constraint for the fluid. }
\rev{Using this constraint,} summing over $i=1,2$ in the equation \eqref{conc}, we infer
\begin{align}
\divg(\J_{\vphi_1}+\J_{\vphi_2})=0,\label{sumJ}
\end{align}
and we will thus assume $\J_{\vphi_1}+\J_{\vphi_2}=\mathbf 0$ to satisfy this constraint. \rev{This assumption is quite standard for classical Cahn Hilliard models (see, e.g., \cite{CahnH1971,HH77, DeGennes, Otto, AGG}). In fact, in \cite{CahnH1971} it is assumed that the net flux of component $1$ is opposite to the net flux  of component $2$.  An alternative, which we do not consider in the present contribution, is to simply assume \eqref{sumJ} as a constraint. The relaxation to divergence free total flux was initially
proposed in \cite{WP} and then studied in \cite{Otto, Cances}.}

We now introduce the order parameter $\varphi:=\varphi_1-\varphi_2$ and, by taking the difference in \eqref{conc}, infer
\begin{align}
\dt\varphi+\divg(\J_\varphi)+\divg(\varphi\u)=0,
    \label{diff}
\end{align}
where $\J_\vphi:=\J_{\vphi_1}-\J_{\vphi_2}$. Here notice that $\rho$ can be written as $\rho=\widetilde{\rho}_1\frac{1+\vphi}2+\widetilde{\rho}_2\frac{1-\vphi}2$ and $\J_\rho=\frac{\widetilde{\rho}_1-\widetilde{\rho}_2}{2}\J_\vphi$.

Let us now consider a conservation of linear momentum law. To this aim, we introduce the baricentric tangential velocity $\widetilde{\v}$ as
$$
\rho\widetilde{\v}:=\rho_1\v_1+\rho_2\v_2=\J_\rho+\rho\v.
$$
Then we consider a surface element $\widehat{\Sigma}(t)\subset \Gamma(t)$ evolving with velocity $\widetilde{\v}+\VVn$ and write the following linear momentum balance equation:
$$
\frac{d}{dt}\ints{\rev{\widehat{\Sigma}}} \rho \u\ds=\rev{\int_{\partial\widehat\Sigma(t)}\mathbf T\boldsymbol\tau\ \d\mathcal H^{n-1}}=\ints{\rev{\widehat{\Sigma}}} \divg\mathbf{T}\ds,
$$
where $\mathbf{T}$ is a symmetric tangential (i.e., $\mathbf T^T=\mathbf T$ and $\mathbf{T}=\P\mathbf{T}\P$, see \cite{gurtin}) stress tensor to be defined later on. \rev{Note that in the last identity we have again used the classical divergence theorem for manifolds.} Observe that, \rev{again by the transport of \cite[Theorem 32]{BGNhandbook}, we have}
\begin{align}
&\frac{d}{dt}\ints{\rev{\widehat{\Sigma}}} \rho \u\ds=\ints{\rev{\widehat{\Sigma}}} (\dtb (\rho \u)+\rho\u\divg(\widetilde{\v}+\VVn))\ds,
\label{al}
\end{align}
where here $\dtb$ stands for the material time derivative with respect to the advective velocity $\hat{\V}=\widetilde{\v}+\VVn$.
Recalling \eqref{scalar}, we have 
$$
\dtb (\rho u_i)= \dt (\rho u_i)+\nablag(\rho u_i)\cdot \widetilde{\v},
$$
so that component-wise we can write
\begin{align*}
    &\dtb (\rho u_i)+\rho u_i\divg(\widetilde{\v}+\VVn)\\&=\dt (\rho u_i)+\nablag(\rho u_i)\cdot \widetilde{\v}+\rho u_i\divg(\widetilde{\v}+\VVn)\\&
    =\dt (\rho u_i)-\nablag(\rho u_i)\cdot\VVn+\divg(\rho u_i (\widetilde{\v}+\VVn))\\&
    =\dt (\rho u_i)-\nablag(\rho u_i)\cdot\VVn+\divg( u_i (\rho({\v}+\VVn)+\J_\rho))\\&
    =u_i\dt\rho +\rho \dt u_i-\nablag(\rho u_i)\cdot\VVn+\divg( u_i (\rho\u+\J_\rho)).
\end{align*}
Note that 
$$
\{\divg (u_i\rho\u)\}_{i=1,\ldots,n}=\divg(\rho \u\otimes \u), \qquad \{\divg (u_i\rho\J_\rho)\}_{i=1,\ldots,n}=\divg(\rho \u\otimes \J_\rho),
$$
and, since $-\nablag(\rho u_i)\cdot\VVn=0$, exploiting equation \eqref{ma} we infer 
\begin{align*}
   & \dtb (\rho \u)+\rho\u\divg(\widetilde{\v}+\VVn)\\&
   =\u\dt \rho +\rho\dt \u+\divg( \u\otimes (\rho\u+\J_\rho))\\&
   =-\u\divg(\rho \u) -\u\divg \J_\rho +\rho\dt \u+\divg( \u\otimes (\rho\u+\J_\rho))\\&=
   \rho \dt \u +((\rho\u+\J_\rho)\cdot\widehat{\nabla}_\Gamma)\u.
\end{align*}
Coming back to \eqref{al}, by the arbitrariness of $\widehat{\Sigma}(t)$, we have then obtained the desired local equation for $\u$:
\begin{align}
\rho \dt \u+((\rho\u+\J_\rho)\cdot\widehat{\nabla}_\Gamma)\u-\divg\mathbf{T}=\mathbf0.
    \label{balancefinal}
\end{align}


Let us now introduce the energy density $e$ as 
\begin{align}
e(\u,\vphi,\nabla\vphi)=\frac{\rho}2\vert \u\vert^2+\Psi(\vphi)+\frac12\vert \nablag\vphi\vert^2,
\label{energy density}
\end{align}
with $\Psi:\R\to \R$ a suitable double well potential. 
Consider a surface element $\widetilde{\Sigma}(t)\subset \Gamma(t)$ evolving with velocity $\u$. To describe the change of the free energy due to diffusion, we introduce
a chemical potential $\mu$ such that
$$
-\ints{\partial\widetilde{\Sigma}}\mu\J_\vphi\cdot\boldsymbol{\xi}\ds-\ints{\partial\widetilde{\Sigma}}\frac{\vert \v\vert^2}{2}\J_\rho\cdot\boldsymbol{\xi}\ds,
$$
with $\boldsymbol{\xi}$ as the outer unit conormal  to $\partial\widetilde{\Sigma}(t)$, is the energy transported into $\widetilde{\Sigma}(t)$ by diffusion. Furthermore, we assume the existence of a generalized vectorial tangential surface force $\boldsymbol\zeta$ such that
$$
\ints{\partial\widetilde{\Sigma}}\dtb \vphi\boldsymbol{\zeta}\cdot\boldsymbol{\xi}\ds
$$
is the working due to microscopic stresses. Here and in the sequel $\dtb$ stands for the material time derivative with respect to the advective velocity $\hat{\V}=\u$.  In conclusion,
$\ints{\partial\widetilde{\Sigma}} \mathbf{T}\boldsymbol{\xi}\cdot \u\ds$
 represents the working on $\widetilde{\Sigma}(t)$  due to the macroscopic stresses in the fluid. By the second law of thermodynamics in an isothermal situation (see for instance \cite{poli}), we can then infer the following dissipation inequality:
 \begin{align}
     \frac{d}{dt}\int_{\widetilde{\Sigma}(t)}e(\u,\vphi,\nablag\vphi)\ds\leq \ints{\partial\widetilde{\Sigma}} \mathbf{T}\boldsymbol{\xi}\cdot \u\ds+\ints{\partial\widetilde{\Sigma}}\dtb \vphi\boldsymbol{\zeta}\cdot\boldsymbol{\xi}\ds-\ints{\partial\widetilde{\Sigma}}\mu\J_\vphi\cdot\boldsymbol{\xi}\ds-\ints{\partial\widetilde{\Sigma}}\frac{\vert \v\vert^2}{2}\J_\rho\cdot\boldsymbol{\xi}\ds
     \label{complete}
 \end{align}
for every surface element $\widetilde{\Sigma}(t)$ transported by the flow.
Let us develop the left-hand side. Recalling \eqref{energy density} and, e.g., \cite[Proposition 2.8]{DE}, we have
\begin{align}
   \nonumber &\frac 12\frac{d}{dt}\ints{\widetilde{\Sigma}}\vert\nablag\vphi\vert^2\ds\\&=\non\ints{\widetilde{\Sigma}}\nabla_\Gamma\vphi\cdot \nablag\dtb\vphi\ds-\frac12\ints{\widetilde{\Sigma}}(\widehat{\nabla}_\Gamma\u+\widehat{\nabla}^T_\Gamma\u)\nablag\vphi\cdot \nablag\vphi\ds\\&
    \non=\ints{\widetilde{\Sigma}}\nabla_\Gamma\vphi\cdot \nablag\dtb\vphi\ds-\frac12\ints{\widetilde{\Sigma}}\P (\widehat{\nabla}_\Gamma\u+\widehat{\nabla}_\Gamma^T\u)\P\nablag\vphi\cdot \nablag\vphi\ds\non\\
    &=\ints{\widetilde{\Sigma}}\nabla_\Gamma\vphi\cdot \nablag\dtb\vphi\ds-\ints{\widetilde{\Sigma}}\E_S(\u)\nablag\vphi\cdot \nablag\vphi\ds,
    \label{o1}
\end{align}
where $\E_S(\u):=\frac12\P(\widehat{\nabla}_\Gamma\u+\widehat{\nabla}_\Gamma^T\u)\P=\frac12(\nablag\u+\nablag^T\u)$. Here we exploited the facts that $\divg\u=0$, $\P\nablag\vphi=\nablag\vphi$ and $\P$ is selfadjoint with respect to the Euclidean inner product.
Then, similarly, by the standard transport theorem and recalling again that $\divg\u=0$, we have 
\begin{align}
\label{o2}
\frac{d}{dt}\ints{\widetilde{\Sigma}}\Psi(\vphi)\ds=\ints{\widetilde{\Sigma}}\Psi'(\vphi)\dtb\vphi\ds=\ints{\widetilde{\Sigma}}\Psi'(\vphi)\dt\vphi+\ints{\widetilde{\Sigma}}\Psi'(\vphi)\nablag\vphi\cdot (\u-\VVn)\ds.
\end{align}
In conclusion, we have 
\begin{align}
&\non\frac 1 2\frac d{dt}\ints{\widetilde{\Sigma}} \rho\vert \u\vert^2\ds=\frac12\ints{\widetilde{\Sigma}} \dtb(\rho\vert\u\vert^2)\ds\\&\non=\ints{\widetilde{\Sigma}} (\frac12\dt\rho\vert\u\vert^2+\rho\u\cdot\dt \u)\ds+\frac12\ints{\widetilde{\Sigma}} \nablag(\rho\vert\u\vert^2)\cdot(\u-\VVn)\ds\\&\non
=\ints{\widetilde{\Sigma}}\left(-\frac12\vert\u\vert^2\divg\J_\rho-\frac12\vert\u\vert^2\divg(\rho\u)\right)\ds\\&\non\quad 
-\ints{\widetilde{\Sigma}}\left(\u\cdot\left((\rho\u+\J_\rho)\cdot\widehat{\nabla}_\Gamma)\u\right)\non+\u\cdot\divg\mathbf{T}\right)\ds\\&\quad\non+\ints{\widetilde{\Sigma}}\left(\frac12\vert\u\vert^2\nablag\rho\cdot \u+\rho\u\cdot ((\u\cdot\widehat{\nabla}_\Gamma)\u)\right)\ds\\&\non
=\ints{\widetilde{\Sigma}}\left(\divg\left(\mathbf T\u-\frac 12\vert \u\vert^2\J_\rho\right)-\mathbf T: \widehat{\nabla}_\Gamma \u\right)\ds\\&\label{o3}
=\ints{\widetilde{\Sigma}}\left(\divg\left(\mathbf T\u-\frac 12\vert \u\vert^2\J_\rho\right)-\mathbf T: \E_S(\u)\right)\ds,
\end{align}
where we have used the facts that $\mathbf{T}=\mathbf T^T$ and $\P\mathbf T\P=\mathbf{T}$, as well as $\nablag(\rho\vert\u\vert^2)\cdot(\u-\VVn)=\nablag(\rho\vert\u\vert^2)\cdot\u$. We can then conclude the argument putting everything together in \eqref{complete}. First, we notice that the symmetry of $\mathbf T $ and the divergence theorem ensures:
\begin{align}
   \non&\ints{\partial\widetilde{\Sigma}} \mathbf{T}\boldsymbol{\xi}\cdot \u\ds+\ints{\partial\widetilde{\Sigma}}\dtb \vphi\boldsymbol{\zeta}\cdot\boldsymbol{\xi}\ds-\ints{\partial\widetilde{\Sigma}}\mu\J_\vphi\cdot\boldsymbol{\xi}\ds-\ints{\partial\widetilde{\Sigma}}\frac{\vert \v\vert^2}{2}\J_\rho\cdot\boldsymbol{\xi}\ds\\&=\ints{\widetilde{\Sigma}} \divg(\mathbf{T}\u)\ds+\ints{\widetilde{\Sigma}}\divg(\dtb \vphi\boldsymbol{\zeta})\ds-\ints{\widetilde{\Sigma}}\divg(\mu\J_\vphi)\ds-\ints{\widetilde{\Sigma}}\divg(\frac{\vert \v\vert^2}{2}\J_\rho)\ds.\label{o4}
\end{align}
\rev{Then, recalling the definition of the energy density given in \eqref{energy density}, plugging \eqref{o1}-\eqref{o4} in \eqref{complete}, we end up with}
\begin{align*}
    &\ints{\widetilde{\Sigma}}\left(\divg\left(\mathbf T\u-\frac 12\vert \u\vert^2\J_\rho\right)-\mathbf T: \E_S(\u)\right)\ds\\&+\ints{\widetilde{\Sigma}}\Psi'(\vphi)\dt\vphi\ds+\ints{\widetilde{\Sigma}}\Psi'(\vphi)\nablag\vphi\cdot \u\ds\\&+\ints{\widetilde{\Sigma}}\nabla_\Gamma\vphi\cdot \nablag\dtb\vphi\ds-\ints{\widetilde{\Sigma}}\E_S(\u)\nablag^T\vphi\cdot \nablag\vphi\ds\\&\leq\ints{\widetilde{\Sigma}} \divg(\mathbf{T}\u)\ds+\ints{\widetilde{\Sigma}}\divg(\dtb \vphi\boldsymbol{\zeta})\ds-\ints{\widetilde{\Sigma}}\divg(\mu\J_\vphi)\ds-\ints{\widetilde{\Sigma}}\divg(\frac{\vert \v\vert^2}{2}\J_\rho)\ds.
\end{align*}
Now notice that, since $\boldsymbol\zeta$ is tangential, 
$$
\divg(\dtb \vphi\boldsymbol{\zeta})=\nablag\dtb\vphi\cdot \boldsymbol \zeta+\dtb\vphi\divg\boldsymbol\zeta=\nablag\dtb\vphi\cdot \boldsymbol \zeta+\dt\vphi\divg\boldsymbol\zeta+(\nablag\vphi\cdot \u)\divg\boldsymbol \zeta,
$$
and, by \eqref{diff}, 
\begin{align*}
    &\divg(\mu\J_\vphi)=\nablag \mu \cdot \J_\vphi+\mu\divg \J_\vphi=\nablag \mu \cdot \J_\vphi-\dt\vphi\mu-\divg(\vphi\u)\mu
    =\nablag \mu \cdot \J_\vphi-\dt\vphi\mu-\mu\nablag\vphi\cdot \u,
\end{align*}
where in the last identity we used $\divg\u=0$.
By the arbitrariness of $\widetilde{\Sigma}(t)$ we can thus deduce the following pointwise inequality:
\begin{align*}
    &(\dt\vphi+\nablag\vphi\cdot \u)\left(\Psi'(\vphi)-\divg\boldsymbol\zeta-\mu\right)+\nablag\dtb\vphi\cdot \left(\nablag\vphi-\boldsymbol\zeta\right)\\&-(\boldsymbol T+\nablag\vphi\otimes\nablag\vphi):\E_S(\u)
    +\nablag\mu\cdot \J_\vphi \leq 0.
\end{align*}
We only need now to make constitutive assumptions to satisfy the inequality. First, since the quantity $\nablag\dtb\vphi$ is arbitrarily independent from the rest, we get 
$$
\boldsymbol\zeta=\nablag\vphi.
$$
Then, we consider the Fick's law for $\mu$, imposing
$$
\J_\vphi=-m(\vphi)\nablag\mu,
$$
were $m\geq0$ is a suitable mobility constant. We then define
$$
\mu:=\Psi'(\vphi)-\divg\boldsymbol\zeta=\Psi'(\vphi)-\Delta_\Gamma\vphi.
$$
In conclusion, we are left to satisfy
\begin{align}
-(\widetilde{\mathbf{S}}+\nablag\vphi\otimes\nablag\vphi):\E_S(\u)\leq 0,
\label{finish}
\end{align}
where we substituted $\mathbf T$ with $ \widetilde{\mathbf{S}}=\mathbf T+\pi\P$, where $p$ is the pressure term. This is necessary since $\mathbf T$ is actually defined only up to a constant $p\mathbf P$, due to the constraint $\divg \u=0$.
Now, the term ${\textbf{S}}:=\widetilde{\mathbf{S}}+\nablag\vphi\otimes\nablag\vphi$ is the viscous stress tensor,
 since it corresponds to irreversible changes of energy due to friction. Notice that it is symmetric and satisfies $\P\mathbf{S}\P=\mathbf{S}$. Indeed, since $\P^T=\P$, we have $\P\nablag\vphi\otimes\nablag\vphi\P=(\P\nablag\vphi)\otimes (\P\nablag\vphi)=\nablag\vphi\otimes\nablag\vphi$. We can then set
 $$
{\textbf{S}}=\widetilde{\mathbf{S}}+\nablag\vphi\otimes\nablag\vphi:=2\nu(\vphi)\E_S(\u),
 $$
 for some $\vphi$-dependent viscosity $\nu(\vphi)\geq0$, to satisfy inequality \eqref{finish}. Therefore, we have now all the ingredients to come up with the final system: we have 
 \begin{align}
     \begin{cases}
     \rho \dt \u +((\rho\u+\J_\rho)\cdot\widehat{\nabla}_\Gamma)\u-2\divg(\nu(\vphi)\E_S(\u))+\nablag \pi+H\pi\n=-\divg(\nablag\vphi\otimes\nablag\vphi),\\
         \divg\u=0,
         \\\dt\vphi-\divg(m(\vphi)\nablag\mu)+\divg(\vphi\u)=0,\\
          \mu=-\Delta_\Gamma\vphi+\Psi'(\vphi),
     \end{cases}
     \label{gensystem}
 \end{align}
where $\J_\rho=\frac{\widetilde{\rho}_1-\widetilde{\rho}_2}{2}\J_\vphi=-\frac{\widetilde{\rho}_1-\widetilde{\rho}_2}{2}\nablag\mu$. Notice that, since $\divg\u=0$, we  compute, using the transport theorem on evolving hypersurfaces, see
Theorem 32 in \cite{BGNhandbook},
$$
\frac d{dt}\vert \Gamma(t)\vert=\intg \divg \u\ds=0
$$
where $\vert \cdot\vert$ is the $n$-dimensional Hausdorff measure. This means that the total surface area is conserved, which is basically a consequence of mass conservation.
Using the divergence theorem on hypersurfaces, we obtain, see Theorem 21 in \cite{BGNhandbook},
\begin{equation}0= \intg \divg \u\ds = \intg H \u\cdot \n\ds= \intg H\Vn\ds, \label{constraint} \end{equation}
where $\Vn$ is the normal velocity. Observe that the global inextensibility constraint is actually only related to the geometric evolution of $\Gamma$, i.e., to $\VVn=\Vn \n$.

The system \eqref{gensystem} is in particular an equation both for the normal and the tangential part of the velocity $\u$. We will now, similar as discussed in \cite{BrandnerReuskenSchwering} and \cite{ES}, consider a given family of
hypersurfaces $(\Gamma(t))_{t\geq 0}$. This mean that the normal part of the velocity is prescribed and only the tangential part of the velocity is an unknown. As we have less unknowns, we also have to reduce the number of equations and hence
project the first equation in \eqref{gensystem} onto the tangent space, so that we find an equation only for $\v$, i.e., for the tangential volume averaged velocity (we have  $\u=\v+\VVn$). 
Due to \eqref{constraint}, we need to require $\intg H\Vn\ds= 0$ for the normal part of the velocity.

\noindent It is easy to observe that 
$$
\nablag\VVn=\mathbf{H}\Vn,
$$
where $\mathbf{H}:=\nablag\n$, as well as (see \cite{V1})
$$
\P\dt\VVn=-\frac12\nablag(\Vn)^2.
$$
Moreover, observe that 
\begin{align*}
&\P((\rho\u+\J_\rho)\cdot\widehat{\nabla}_\Gamma)\u=\P(\widehat{\nabla}_\Gamma\u)(\rho\u+\J_\rho)=(\nablag\u)(\rho\u+\J_\rho)\\&=(\nablag\v)(\rho\v+\J_\rho)+(\nablag\VVn)(\rho\v+\J_\rho)
=(\nablag\v)(\rho\v+\J_\rho)+\rho \Vn\mathbf{H}\v+\Vn\H\J_\rho.
\end{align*}
We also notice that, since $\mathbf{H}$ is symmetric,
$$
2\P\divg(\nu(\vphi)\E_S(\VVn))=2\P\divg(\nu(\vphi)\Vn\mathbf{H}).
$$
Therefore, we can thus conclude that the system for the sole tangential velocity $\v$ reads: \\
Given $\VVn=\Vn\n$ satisfying the compatibility assumption $$ \intg H\Vn\ds=0  \qquad \text{ (globally inextensible material surface)}, $$
find $(\v,p,\vphi)$ such that
 \begin{align}
     \begin{cases}
     \rho \P\dt \v +((\rho\v+\J_\rho)\cdot\nabla_\Gamma)\v+\rho \Vn\mathbf{H}\v+\Vn\H\J_\rho-2\P\divg(\nu(\vphi)\E_S(\v))+\nablag \pi\\
     =-\P\divg(\nablag\vphi\otimes\nablag\vphi)+2\P\divg(\nu(\vphi)\Vn\mathbf{H})+\frac\rho2\nablag(\Vn)^2,\\
         \divg\v=-H\Vn,
         \\\dt\vphi-\divg(m(\vphi)\nablag\mu)+\nablag\vphi\cdot \v=0,\\
         \mu=-\Delta_\Gamma\vphi+\Psi'(\vphi),
     \end{cases}
     \label{mainp}
 \end{align}
with $\J_\rho=-\frac{\widetilde{\rho}_1-\widetilde{\rho}_2}{2}\nablag\mu$.
From now on we set for simplicity $m(\vphi)\equiv 1$, since we are interested in the constant mobility case.

\section{Surface Piola transform and functional setting}\label{sec:Piola}
\subsection{Flow map and Surface Piola Transform}
\label{regflowmap1}
In a similar way as in \cite{V2}, for some $T>0$, let us assume to consider an initial closed orientable surface $\Gamma_0$ of class $C^6$ (with unit normal denoted by $\n_0$) embedded in $\R^3$, together with a vector field $\w\in C^6([0,T]\times \R^3;\R^3)$ so that $\max_{[0,T]\times \R^3}\norma{\w}\leq C$. 
Then the ODE system 
\begin{align*}
&\ddt \mathbf x(t,\z)=\w(t,\mathbf x(t,\z)),\\&
\mathbf x(0,\z)=\z,
\end{align*}
has a unique solution for any $z\in \Gamma_0\subset\R^3$, defining a one-to-one mapping $\Gamma_0\to \gam(t)$ for all $t\in[0,T]$ (see, for instance Theorems II.1.1, V.3.1 and remark to
Theorem V.2.1 in \cite{Hartman}), where $\Gamma(t)$ is still closed, orientable, and embedded in $\R^3$ for any $t\in[0,T]$. Moreover, this mapping is $C^6(\mathcal{S}_0;\R^4)$ (see \cite[Corollary V.4]{Hartman}), where $\mathcal{S}_0:=[0,T]\times\gam_0$. Therefore, we can define $\mathcal{S}$ as the $C^6$ manifold  image of $\mathcal{S}_0\in \C^6$ under a $C^6$-mapping. Clearly both $\mathcal{S}$ and $\mathcal{S}_0$ are closed manifolds. 
Now, with the same arguments as in \cite{V2}, we can construct, through the signed distance function to $\Gamma_0$, a globally $C^5$-smooth extension of the spatial normal $\n(t,x)$, $(t,x)\in\mathcal{S}$. One can then obtain an extended (and not relabeled) normal vector ${\n}\in C^5([0,T]\times \R^3;\R^3)$. Then we introduce the normal component of the flow $\VVn:=(\w\cdot\n)\n\in C^5([0,T]\times\R^3;\R^3)$. Now we can introduce the normal flow map, which is what we need in this context, since we are only interested in the hypersurface $\Gamma(t)$ evolution. In particular, we set $\Phi^n_{(\cdot)}$ as the solution to
 \begin{align*}
&\ddt \Phi_t^n(\z)=\VVn(t,\Phi_t^n(\z)),\\&
\Phi_0^n(\z)=\z,
\end{align*}
which by standard regularity theorems is such that $\Phi^n_{(\cd)}\in C^5(\mathcal{S}_0,\mathcal{S})$. We then need to introduce the Surface Piola transform (see also \cite{V2}). First, for any $t\in[0,T]$ the map $\Phitn: \gz\to \gam(t)$ is a $C^5$-diffeomorphism, so that we can introduce the inverse flow map as to be $\Phi_{-t}^n: \gam(t)\to\gz$ for any $t\in[0,T]$ with the same regularity properties, and define the differential of $\Phi_t^n$ as to be $D\Phi_t^n(\z): T_{\z}\gz\to T_{\mathbf x}\gam(t)$, where $\mathbf x:=\Phitn(\z)$ and $T_{\z}\gz$ is the tangent space to $\gz$ at $\z$ (same, \textit{mutatis mutandis}, for $T_{\mathbf x}\gam(t)$). Clearly $D\Phi_t^n$ is invertible over these spaces. We also set $J=J(t,\z):=\norma{\text{det} D\Phitn(\z)}$ and $J^{-1}:=\norma{\text{det} D\Phi_{-t}^n(\z)}$. Note that $J J^{-1}\equiv 1$. We can then define the matrices $\mathbf{D}=\mathbf{D}(t,\z):=D\Phitn(\z)\P_0:\R^3\to\R^3$, where $\P_0=\mathbf I-(\n_0\otimes\n_0)$ is the projector over $T_{\z}\gz$, and $\mathbf{D}^{-}=\mathbf{D}^{-}(t,\mathbf x):= D\Phi_{-t}^n(\mathbf x)\P:\R^3\to \R^3$, where $\P=\mathbf I-(\n\otimes\n)$ is the projector over $T_{\mathbf x}\gam(t)$. We also denote by $\D^{-T}$ the transpose matrix $(\Dm)^{T}$. Clearly, by construction, $\mathbf D^{-}\mathbf D=\P_0$ as well as $\mathbf D\mathbf D^{-}=\P$. We then introduce the operator $\A:\R^3\to \R^3$ as
\begin{align}
\A(t,\z):=J^{-1}\mathbf D(t,\z)+\n\otimes \n_0,
\end{align}
which is invertible with inverse $\A^{-1}:=J\mathbf{D}^{-}+\n_0\otimes \n$. By construction it holds $\A_{\vert T_\z\gz}=J^{-1}D\Phitn$, for any $\z\in \gz$. Note also that $\text{det}\A\equiv 1$. Furthermore, thanks to the regularity of the flow map, we can easily prove, in a similar way to \cite[Lemma 3.1]{V2}, that 
\begin{align}
\norm{J}_{C^4(\mathcal{S}_0)}+\norm{J^{-1}}_{C^4(\mathcal{S})}+\norm{\mathbf D}_{C^4(\mathcal{S}_0)}+\norm{\mathbf{D}^{-}}_{C^4(\mathcal{S})}+\norm{\A}_{C^4(\mathcal{S}_0)}+\norm{\A^{-1}}_{C^4(\mathcal{S})}\leq C(T).
    \label{regularity_basic}
\end{align}
Notice that, since $\mathcal S$ and $\mathcal{S}_0$ are compact in $\R\times\R^3$, we also have $C^{4}(\mathcal S)\subset C^{2,1}(\mathcal S)$ as well as $C^{4}(\mathcal S_0)\subset C^{2,1}(\mathcal S_0)$, which is the minimal regularity to perform the local existence argument in the proof of our main result.  
Recalling that $\H:=\nablag\n$ and $H=\text{tr}\H$, we also deduce the following regularity
\begin{align}
\norm{H}_{C^4(\mathcal{S})}+\norm{\vn}_{C^4(\mathcal{S})}\leq C(T).
    \label{regularity_basic2}
\end{align}
Now we can introduce some pullback and pushforward maps. First, we set $\tphi_t$ such that for a given $\f: \mathcal{S}_0\to \R^q$, for some $q\geq1$,
$$
(\tphi_{t}\f)(t,\mathbf x)=\f(t,\Phi_{-t}^n(\mathbf x)),
$$
together with its inverse $\tphimt$, defined as, for any $\mathbf g: \mathcal{S}\to \R^q$,
$$
(\tphimt\mathbf g)(t,\z):=\mathbf g(t,\Phitn(\z)).
$$
Then we also introduce another map, which is essential to preserve the divergence-free property of tangential vector fields. In particular, we define the surface Piola transform as
for a given $\f: \mathcal{S}_0\to \R^3$,
$$
(\phi_{t}\f)(t,\mathbf x)=\A(t,\Phi_{-t}^n(\mathbf x))\f(t,\Phi_{-t}^n(\mathbf x))=\A(t,\Phi_{-t}^n(\mathbf x))(\tphi_t\f)(t,\mathbf x),
$$
as well as its inverse as, for any $\mathbf g: \mathcal{S}\to \R^3$,
$$
(\phi_{-t}\mathbf g)(t,\z):=\A^{-1}(t,\Phitn(\z))\mathbf g(t,\Phitn(\z))=\A^{-1}(t,\Phitn(\z))(\tphimt\g)(t,\z).
$$
This map has the property of mapping tangent vectors into tangent vectors, but also that $\text{div}_{\gam(t)}\f=0$ for almost any $\mathbf x\in \gam(t)$ if and only if $\text{div}_{\gz}\phi_{-t}\f=0$  for almost any $\z\in \gz$ (see, for instance, \cite{41V2}). Note that, setting $\mathbf x=\Phitn(\z)$, in general it holds, for $\f$ tangent vector (i.e. $\f(t,\mathbf x)\cdot\n(t,\mathbf x)=0$ for any $x\in \gam(t)$),

$$(\text{div}_{\gam(t)}\f)(t,\mathbf x)=J^{-1}(t,\z)(\text{div}_{\gz}\phi_{-t}\f)(t,\z).$$
\subsection{Evolving Banach and Hilbert spaces}
Here we introduce the functional setting necessary to perform the analysis of problems on evolving surfaces.  We do it in a summarised way and refer the reader to \cite{Def,AE} for a more rigorous and detailed explanation of how these are constructed and how they can be abstracted to a more generalised setting.
First, for any $t\in[0,T]$, we introduce the Hilbert and Banach spaces we will use. In particular,  the Sobolev spaces are denoted as usual by $W^{k,p}\gt$%
 , where $k\in \mathbb{N}$ and $1\leq p\leq \infty $, with norm $\Vert \cdot
 \Vert _{W^{k,p}\gt}$. The Hilbert space $W^{k,2}\gt$ is denoted
 by $H^{k}\gt$ with norm $\Vert \cdot \Vert _{H^{k}\gt}$. We then set $H^{-1}\gt=(H^1\gt)'$. Moreover, $B^s_{
pq}(\Gamma(t))$ denotes, for any $t\in[0,T]$, the standard Besov space, where $s \in \R$, $1 \leq  p, q \leq \infty$. Given a vector space $X(t)$ of functions defined on the surface $\Gamma(t)$, we denote by $\mathbf{X}(t)$ the generic space of tangential vectors or matrices, with each component belonging to $X(t)$. Recall that a tangential vector $\v$ is such that $\v\cdot\n=0$ on $\gam(t)$, i.e., $\P(t)\v=\v$. Analogously, a tangential matrix $\mathbf M$ is a matrix such that $\P(t)\mathbf M\P(t)=\mathbf M$. In this multi-component case $\vert \mathbf{v} \vert$ is the Euclidean norm
 of $\mathbf{v}\in  \mathbf{X}(t)$, i.e., $\vert\mathbf{v}\vert^2=\sum_{j}\|v_j\|_{X(t)}^2$.
 We then denote by $(\cdot,\cdot )$ the inner product in $\mathbf{L}^{2}\gt$ and by $\Vert \cdot \Vert $
 the induced norm. 
 \begin{align*}
     &\L^2_\sigma\gt:=\{\v\in \L^2\gt: \divg \v=0\text{ a.e. on }\gam(t)\},
     \\
     &\H^1_\sigma\gt:=\{\v\in \H^1\gt: \divg \v=0\text{ a.e. on }\gam(t)\}.
 \end{align*}
Notice that, for any $t\in[0,T]$, it holds the Hilbert triplet $\H^1_\sigma\gt\hookrightarrow \L^2_\sigma\gt \hookrightarrow \H^1_\sigma\gt'$ \rev{(see, e.g., \cite{Hebey})}. Now suppose, for $t\in [0,T]$, $X(t)$ to be a generic scalar Banach space of functions over $\Gamma(t)$ (and $\mathbf X(t)$ a space of tangential vector fields over $\Gamma(t)$), and let $X_0=X(0)$ ($\mathbf X_0=\mathbf X(0)$, respectively). Let us consider $t=0$. We denote by $L^q (0,T; X_0 )$ the Bochner space of $X_0$-valued $q$-integrable (or essentially bounded) functions. We then denote by $BC([0,T];X_0)$ the Banach space of bounded continuous functions on $[0,T]$, equipped with the supremum norm. The space $BUC([0,T];X_0)$ is then its subspace of bounded and uniformly continuous functions. The function space 
 $C^{0,\alpha} (0,T ; X_0)$, for $\alpha\in[0,1]$ denotes the vector space of all $\alpha$-H\"{o}lder continuous functions $f : [0,T]\to X_0$. Finally, $W^{1,p} (0, T ; X_0)$, $1 \leq p < \infty$, is the space of functions $f$ such that $\partial_t f\in L^p(0,T;X_0)$ and $f\in L^p(0,T;X_0)$, where $\partial_t$ denotes the vector-valued distributional derivative of $f$. We then set $H^1 (0, T ; X_0) =W^{1,2}(0,T;X_0)$.
  
 We now pass to consider the  setting of evolving surfaces: We study functions of the form 
\begin{align*}
u\colon [0, T]&\to \bigcup_{t\in[0,T]} X(t)\times \{t\}, \quad t\mapsto (\bar{u}(t), t)
\end{align*}
and identify $u(t)\equiv \overline u(t)$ (same definition for vector fields and matrices); roughly speaking, we want to see $u(t)$ as an element of the space $X(t)$. The function spaces are defined as follows. For $p\in [1,\infty]$, $u\in L^p_{X}$ if  $t\mapsto \tphi_{-(\cdot)}\bar{u}(\cdot)\in L^p(0, T; X_0)$ with the norm
 \begin{align*}
     \|u\|_{L^p_X} := \begin{cases}
         \left( \int_0^T \|u(t)\|_{X(t)}^p dt\right)^{1/p} &\text{ if } p<\infty, \\
         \text{ess sup}_{t\in [0,T]} \|u(t)\|_{X(t)} &\text{ if } p=\infty.
     \end{cases}
 \end{align*}
 If $X(t)=H(t)$ are Hilbert spaces, then so is $L^2_H$ with the corresponding inner product
\begin{align*}
(u,v)_{L^2_{H}} &:= \int_0^T (u(t), v(t))_{H(t)}dt. 
\end{align*} 
In the case of $\mathbf X(t)=\L^2_\sigma\gt, \H^1_\sigma\gt,\H^2\gt$, which are vector spaces of tangential vector fields, we adopt as pullback the map $\phi_t$ and define, for $p\in [1,\infty]$, $\v\in L^p_{\mathbf X}$ if  $t\mapsto \phi_{-(\cdot)}\bar{\v}(\cdot)\in L^p(0, T; \mathbf X_0)$ with the norm
 \begin{align*}
     \|\v\|_{L^p_\mathbf X} := \begin{cases}
         \left( \int_0^T \|\v(t)\|_{\mathbf X(t)}^p dt\right)^{1/p} &\text{ if } p<\infty, \\
         \text{ess sup}_{t\in [0,T]} \|\v(t)\|_{\mathbf X(t)} &\text{ if } p=\infty.
     \end{cases}
 \end{align*}
 We now observe that, thanks to \eqref{regularity_basic}, the pairs $(\tphi_t, X(t))_t$ are compatible according to the definition in \cite{Def} for $X(t)=H^{-1}\gt,\ W^{k,p}\gt$, $k=0,1,\ldots,4$, $p\in[1,\infty]$. Analogously, the pairs $(\phi_t, \mathbf X(t))_t$ are compatible for $\mathbf X(t)=\L^2_\sigma\gt, \H^1_\sigma\gt,\H^2\gt$.
 
Then, $u\in C^1_{X}$ if $t\mapsto \tphi_{-t}u(t)\in C^1([0,T]; X_0)$ and we define its time derivative as 
\begin{align}
\partial^{\bullet}_t u(t)=\dt u(t)= \tphi_t \dfrac{d}{dt} \tphi_{-t} u(t).\label{dtclass}
\end{align}
Note that, since $\Phitn$ is the normal flow map, the time derivative coincides with the normal material derivative.
 These spaces are endowed with the norm
\begin{align*}
\|u\|_{C^1_{X}} := \sup_{t\in [0,T]} \|u(t)\|_{X(t)} +  \sup_{t\in [0,T]} \|\dt u(t)\|_{X(t)}.
\end{align*} 
In the case of spaces of divergence free vector fields $\mathbf X(t)$ like $\L^2_\sigma\gt$ and $\H^1_\sigma\gt$, we define, as in \cite{V2}, two types of time derivatives. We set the standard time derivative as in \eqref{dtclass}, i.e.,
$\v\in C^1_{\mathbf X}$ if $t\mapsto \tphi_{-t}u(t)\in C^1([0,T]; \X_0)$ and we define its time derivative as 
\begin{align}
\partial^{\bullet}_t \v(t)=\dt \v(t)= \tphi_t \dfrac{d}{dt} \tphi_{-t} \v(t).\label{dtclass1}
\end{align}
 This derivative is not divergence free, though. Therefore, we also introduce a new kind of time derivative, which is also consistent with the abstract framework developed in \cite{Def} when applied to the spaces $\L^2_\sigma\gt$ and $\H^1_\sigma\gt$. Namely, we set
 \begin{align*}
     \dts \v:=\phi_t\ddt\phi_{-t}\v(t),
 \end{align*}
 which satisfies, for a divergence free tangential vector field $\v$, $\P\dts\v=\dts\v$ and $\divg\dts\v=0$ almost everywhere in $\gam(t)$ \rev{(see, e.g., \cite{41V2} for the properties of surface Piola transform, or \cite[Eq. (3.20)]{V2})}. Thanks to \cite[Lemma 3.6]{V2}, we have a precise relation between the two notions, since we know that, for smooth vector fields,
 \begin{align}
\dt\v=\dts\v-\A(\dt\A^{-1})\v,\quad \P\dt\v=\dts\v-\P\A(\dt\A^{-1})\v.
    \label{relation_base1}
\end{align}
Now we can also introduce the notions of weak time derivatives. In particular, let $u\in L^2_{H^1}$. A function $v\in L^2_{H^{-1}}$ is the weak time derivative of $u$, and we can write $v=\dt u$, if, for any $\eta\in \mathcal{D}_{H^1}$ (where $\eta\in \mathcal{D}_{X}$ if $t\mapsto \phi_{-t}u(t)\in C_0^\infty((0,T); X_0)$) it holds
\begin{align*}
&\int_0^T \langle v(t), \eta(t)\rangle_{H^{-1}(\Gamma(t)), H^1(\Gamma(t))} dt = -\int_0^T (u(t),\dt\eta(t))dt - \int_0^T \int_{\Gamma(t)} u(t)\eta(t) v_\n H\ds dt.
\end{align*}
Here, as discussed above, $\dt \eta$ is the strong material derivative of $\eta$. We will use $\dt u$ for both the strong and weak material derivatives, since the two notions coincide whenever $u$ is sufficiently regular. 
The same definition holds for the time derivative $\dt\v$, where $\v$ is tangential vector field. We can write $\w=\dt\u$, if, for any $\boldsymbol\eta\in \mathcal{D}_{{\H^1}}$ it holds
\begin{align*}
&\int_0^T \langle \w(t), \boldsymbol\eta(t)\rangle_{\H^{1}(\Gamma(t))',\H^1\gt)}dt = -\int_0^T (\v(t),\dt\boldsymbol\eta(t))dt - \int_0^T \int_{\Gamma(t)} \v(t)\cd\boldsymbol\eta(t) v_\n H\ds dt.
\end{align*}
Concerning the $\dts$ time derivative, by exploiting \eqref{relation_base1} for smooth solutions, we can also define a suitable weak notion as:
$\tilde{\w}=\dts\u$, if, for any $\boldsymbol\eta\in \mathcal{D}_{{\H^1}}$ it holds
\begin{align*}
&\int_0^T \langle \tilde{\w}(t), \boldsymbol\eta(t)\rangle_{\H^{1}(\Gamma(t))',\H^1\gt)}dt\\& = -\int_0^T (\v(t),\dt\boldsymbol\eta(t))dt - \int_0^T \int_{\Gamma(t)} \v(t)\boldsymbol\eta(t) v_\n H\ds dt+\int_0^T\ig \A\dt(\A^{-1})\v\cd \boldsymbol\eta dt,
\end{align*}
entailing by comparison in the formulations that the relations \eqref{relation_base1} also hold between the weak notions of time derivatives. 

To conclude this section, we finally settle down the following notation: if a function is not defined on the entire interval of time $[0,T]$, but only on $[0,T_1]$, with $0<T_1\leq T$, we denote all the corresponding spaces as, for instance, $L^p_{X(T_1)}$, to highlight the dependence on the smaller time $T_1$.

\section{Assumptions and main result}\label{mainresult}
We consider the following regularity assumptions:
\begin{enumerate}[label=(\subscript{L}{{\arabic*}})] 
\item \label{l1}$\Psi\in C^{4}(\R)$,
\item \label{l2}$\nu\in C^{2}(\R)$, such that $\nu(s) \geq \nu_* > 0$  for every $s \in \R$
and some $\nu_* > 0$. 
\item \label{ro}
$  \rho(s):=\frac{\widetilde{\rho}_1+\widetilde{\rho}_2}{2}+\frac{\widetilde{\rho}_2-\widetilde{\rho}_1}{2}s,\ \forall s\in \R,\quad \widetilde{\rho}_1,\widetilde{\rho}_2>0$.
\end{enumerate}
We now state our main result, concerning the short time existence and uniqueness of a strong solution in the following sense: 
 
\begin{theorem}[Short time well-posedness of strong solutions]
	\label{strong1a}
	Assume the regularity assumptions on $\Phitn$ stated in Section \ref{regflowmap1}, and assumptions \ref{l1}-\ref{ro} for $\Psi$, $\nu$, and $\rho$. Let $\mathbf{v}_0 \in \H^1(\gz)$
	and $\vphi_0\in B_{p,q}^{4-\frac4q}(\gz)$, for some $2<q\leq 4$ and $p>4$. Then there exists a $\tT\in(0,T]$ and a unique  strong solution $(\mathbf{v},\pi, \vphi)$ to \eqref{mainp} defined on $[0,\tT]$,
 such that 
 \begin{align*}
	&\mathbf{v} \in L ^2_{\H^2(\tT)},\quad \dt \v\in L^2_{\L^2(\tT)} ,\quad \pi\in L^2_{H^1(\tT)},\quad 
	\vphi\in L^q_{W^{4,p}(\tT)},\quad  \dt\vphi\in L^q_{L^p(\tT)},
	\end{align*}
	and it satisfies \eqref{mainp} in the almost everywhere sense.
 Furthermore, if $\norm{\vphi_0}_{L^\infty(\gz)}\leq1-2\delta_0$ for some $\delta_0>0$, then there exists $T_*\in(0,\tT]$ also depending on $\delta_0$, such that
 \begin{align}
    \sup_{t\in[0,T_*]} \norm{\vphi}_{C\gt}\leq 1-\delta_0. \label{separaz}
 \end{align}
\end{theorem}
\begin{remark}
Recalling the notion of evolving space equivalence in the sense of \cite[Definition 2.31]{Def} for the corresponding spaces, it is easy to infer that $\v \in C_{\H^1(\tT)}$, and $\vphi\in C_{(L^p,W^{4,p})_{1-\frac1q,q}(\tT)}$ (see also \eqref{uniform} below), so that the initial conditions $\v(0)=\v_0$ and $\vphi(0)=\vphi_0$ make sense.
\end{remark}
\begin{remark}
    Clearly, from the solution $(\v,p,\vphi)$ given in the theorem above, one can also retrieve the chemical potential $\mu\in L^q_{W^{2,p}(\tT)}$ exploiting its definition in system \eqref{mainp}.
\end{remark}

\begin{remark}
\label{Psiessential}
    \rev{Let us assume that $\Psi$ is a singular potential, i.e., $\Psi'$ is singular at the pure phases $\pm1$} (as, for instance, in the case of logarithmic potential), and that there exists $\delta_0$ such that $\norm{\vphi_0}_{L^\infty(\gz)}\leq1-2\delta_0$. We can consider a potential $\widetilde{\Psi}$, satisfying assumption \ref{l2}, such that $\widetilde{\Psi}_{\vert[-1+\delta_0,1+\delta_0]}=\Psi_{\vert [-1+\delta_0,1+\delta_0]}$, and then extended in a smooth way outside $[-1+\delta_0,1+\delta_0]$. With this potential $\widetilde{\Psi}$ in place of $\Psi$, Theorem \ref{strong1a} above applies. Since, on an interval $[0,T_*]$, the final solution $\vphi$ is contained in $[-1+\delta_0,1-\delta_0]$, the solution $(\v,{p},\vphi)$ is also a solution to the original problem \eqref{systf}, with $\Psi$ in place of $\tilde{\Psi}$, on the same interval $[0,T_*]$. This allows to show the short time existence of a strong solution also for the cases when $\Psi$ is singular, as long as the initial datum is strictly separated from the pure phases. Clearly, uniqueness cannot be proven in the same way, since the application of Theorem \ref{strong1a} depends on the choice of the regularization of $\Psi$.
\end{remark}
In the next sections we will provide the main ingredients leading to the final proof of Theorem \ref{strong1a} in Section \ref{proof1}.
\section{Equivalent formulation for a divergence-free velocity}
First, we observe that we can write system \eqref{mainp} for a divergence-free velocity. In particular, we solve the problem
\begin{align*}
    -\Delta_\Gamma \Pi=Hv_\n
\end{align*}
for $\Pi$ such that $\overline{\Pi}=0$. This admits a unique solution by standard elliptic theory.  Then we define $\widehat{\u}=\nablag\Pi$ and set $\v=\V+\widehat{\u}$. Note that, as shown in Lemma \ref{appendixb}, the assumptions on the flow map in Section \ref{regflowmap1} (see in particular \eqref{regularity_basic2}) entail that we have $\Pi\in C^1_{W^{4,p}}$, and thus  $\widehat{\u}\in C^0_{\W^{3,p}}$ with $\P\dt\widehat{\u}\in C^0_{\W^{3,p}}$, for any $p\in[2,+\infty)$, entailing that $\widehat{\u}\in C^0_{\C^2}$ and $\P\dt\widehat{\u}\in C^{0}_{\C^2}$, as it is sufficient for all our purposes. In this way we find the system for $(\V,\pi,\vphi)$
 \begin{align}
     \begin{cases}
     \rho \P\dt \V +((\rho\V+\J_\rho)\cdot\nabla_\Gamma)\V+\rho(\widehat{\u}\cdot\nabla_\Gamma)\V+\rho(\V\cdot\nabla_\Gamma)\widehat{\u}+\rho \Vn\mathbf{H}\V+\Vn\H\J_\rho\\-2\P\divg(\nu(\vphi)\E_S(\V))+\nablag \pi
     =-\P\divg(\nablag\vphi\otimes\nablag\vphi)+2\P\divg(\nu(\vphi)\Vn\mathbf{H})+\frac\rho2\nablag(\Vn)^2\\
     \quad-\rho \P\dt \widehat{\u} -((\rho\widehat{\u}+\J_\rho)\cdot\nabla_\Gamma)\widehat{\u}-\rho \Vn\mathbf{H}\widehat{\u}+2\P\divg(\nu(\vphi)\E_S(\widehat{\u})),\\
         \divg\V=0,
         \\\dt\vphi-\Delta_\Gamma\mu+\nablag\vphi\cdot \V
         +\nablag\vphi\cdot (\widehat{\u}+\V_\n)=0,\\
         \mu=-\Delta_\Gamma\vphi+\Psi'(\vphi).
     \end{cases}
     \label{detract}
 \end{align}
Here we clearly have $\V(0)=\V_0:=\v_0-\widehat{\u}(0)$ and $\vphi(0)=\vphi_0$. 
\section{Pullback equations}
 We now need to apply a pullback operator to all the equations of \eqref{detract} on the initial surface $\Gamma_0$.  We define $\u=\phmt{(\V)}$, $\widetilde{\vphi}=\tphi_{-t}(\vphi)$, $\widetilde{\mu}=\tphi_{-t}(\mu)$, and $\widetilde{\rho}=\rho(\tvphi)$. First we observe that
 \begin{align}
\dt\vphi=\tphi_t\left(\frac{d}{dt}\tphi_{-t}(\vphi)\right),
 \label{dtv2}
 \end{align}
 so that $\tphi_{-t}(\dt\vphi)=\partial_t\widetilde{\vphi}$. Note that it is easy to see that this last relation also holds in case of weak time derivatives. Moreover, we recall that
 $$
\partial^\ast_t\V=\phi_t\left(\frac{d}{dt}\phmt(\V)\right),
 $$
 so that $\phmt(\dtast\V)=\partial_t\u$. Now, recalling \eqref{relation_base1}, we have
 $$
 \P\dt\V=\dtast\V+\A\P_0\dt(\A^{-1})\V,
 $$
 where $\P_0=\P(0):\mathbb{R}^3\to T\Gamma_0$. This entails that
 \begin{align}
  \phmt(\P\dt\V)=\partial_t\u+\P_0\tphi_{-t}(\dt(\A^{-1}))\A\u,
  \label{dtv1}
 \end{align}
 which clearly holds also in case of weak time derivatives.
 We also introduce the symbol $\nabla_{\Gamma(t)}^\phi$, which stands for the following 
\begin{align}
\nabla_{\Gamma(t)}^\phi f:=\phmt(\nablag(\tphi_t(f)))=\tphimt(\A^{-1}\D^{-T})\nablagz f ,\quad f:\Gamma_0\to\R,
\label{pressure}
\end{align}
so that, by introducing the pressure term $\widetilde{\pi}:=\tphimt(\pi)$, we have
$$
\phmt(\nablag \pi)=\nabla_{\Gamma(t)}^\phi \widetilde{\pi}.
$$
To avoid tedious computations here, we postpone the complete derivation of the pullback system to Appendix \ref{pullback1}.

In the end, the pullback of the complete system reads: find $(\u,\widetilde{\pi},\tvphi)$ such that
\begin{align}
\label{syst1}
    \begin{cases}
\widetilde{\rho}\partial_t\u+\A_\ast(t;\tvphi,\nablagz\tvphi)(\u,\nablagz\u)
+\b_3(t)\u\cdot(\di_1(t;((\nablagz(\nablagz(\nablagz \tvphi)_i)_j)_k))+\di_2(t; ((\nablagz(\nablagz \tvphi)_i)_j)))\\+(\b_6(t)(\di_1(t;((\nablagz(\nablagz(\nablagz \tvphi)_i)_j)_k))+\di_2(t; ((\nablagz(\nablagz \tvphi)_i)_j)))\cdot \nablagz)\u
\\-\nu(\tvphi)\P_0\a_0(t;\nablagz((\nablagz\u_i)_j)_k)-\nu(\tvphi)\a_1(t)\u-\nu(\tvphi)a_2(t)\u
\\+\b_7(t)(\di_1(t;((\nablagz(\nablagz(\nablagz \tvphi)_i)_j)_k))+\di_2(t; ((\nablagz(\nablagz \tvphi)_i)_j)))+\nabla_{\Gamma(t)}^\phi\widetilde{\pi}=\mathbf{0},\\\, \\
\mathrm{div}_{\Gamma_0}\u=0,\\ \, \\
\partial_t\tvphi+B(t;\tvphi,\nablagz\tvphi)(\u)+f_0(t;(\nablagz(\nablagz(\nablagz(\nablagz \tvphi)_i)_j)_k)_l)+f_1(t;((\nablagz(\nablagz(\nablagz \tvphi)_i)_j)_k))\\+f_2(t; ((\nablagz(\nablagz \tvphi)_i)_j)
-\Psi''(\tvphi)f_3(t;(\nablagz(\nablagz\tvphi)_i)_j)=0,
    \end{cases}
\end{align}
where all the functions $\A_\ast,B, \mathbf a_i,\mathbf b_i,\mathbf d_i,f_i$ are at most 
(bi)linear in their arguments (apart from time $t$ and $\tvphi$). See Appendix \ref{pullback1} for the details.
In order to treat these equations in a simpler way, we also introduce an operator acting like a modified \rev{Helmholtz} projector. In particular, let us define the standard orthogonal (in $\L^2(\Gamma(t))$) \rev{Helmholtz} projector $P^{\Gamma(t)}:\L^2(\Gamma(t))\to \L^2_\sigma(\Gamma(t))$ such that, for any $\v\in \L^2(\Gamma(t))$  
$$
\v=P^{\Gamma(t)}\v+\nablag p,
$$
for some $p\in H^1(\Gamma(t))$, where clearly $\divg P^{\gam(t)}\v=0$. If we now apply $\phi_{-t}$ to this identity, and set $\v=\phi_t\u$, we obtain 
$$
\u=\phi_{-t}(P^{\Gamma(t)}\phi_{t}\u)+\nablag^{\phi}\widetilde{\pi},
$$
with $\widetilde{\pi}=\tphimt(p)$. We thus introduce the modified (not orthogonal in general) operator $P(t)$ as 
\begin{align*}
P(t)(\cdot):=\phi_{-t}(P^{\Gamma(t)}\phi_{t}(\cdot)): \L^2(\Gamma_0)\to \L^2_\sigma(\Gamma_0).
\end{align*}
Notice that it holds $P(0)=P^{\Gamma_0}$ and the spaces in boldface contain only tangential vector fields, e.g. $\L^2(\Gamma_0):=L^2(\Gamma_0;T\Gamma_0)$. Notice that, as expected, $P(t)\nablag^{\phi}\widetilde{\pi}=0$. Therefore, we can apply $P(t)$ to the equation for $\u$, obtaining  
\begin{align}
    \begin{cases}
P(t)\widetilde{\rho}\partial_t\u+P(t)\A_\ast(t;\tvphi,\nablagz\tvphi)(\u,\nablagz\u)
\\+P(t)\b_3(t)\u\cdot(\di_1(t;((\nablagz(\nablagz(\nablagz \tvphi)_i)_j)_k))+\di_2(t; ((\nablagz(\nablagz \tvphi)_i)_j)))\\+P(t)(\b_6(t)(\di_1(t;((\nablagz(\nablagz(\nablagz \tvphi)_i)_j)_k))+\di_2(t; ((\nablagz(\nablagz \tvphi)_i)_j)))\cdot \nablagz)\u
\\-P(t)\nu(\tvphi)\P_0\a_0(t;\nablagz((\nablagz\u_i)_j)_k)-P(t)\nu(\tvphi)\a_1(t)\u-P(t)\nu(\tvphi)a_2(t)\u
\\+P(t)\b_7(t)(\di_1(t;((\nablagz(\nablagz(\nablagz \tvphi)_i)_j)_k))+\di_2(t; ((\nablagz(\nablagz \tvphi)_i)_j)))=\mathbf{0},\\\, \\
\mathrm{div}_{\Gamma_0}\u=0,\\ \, \\
\partial_t\tvphi+B(t;\tvphi,\nablagz\tvphi)(\u)+f_0(t;(\nablagz(\nablagz(\nablagz(\nablagz \tvphi)_i)_j)_k)_l)+f_1(t;((\nablagz(\nablagz(\nablagz \tvphi)_i)_j)_k))\\+f_2(t; ((\nablagz(\nablagz \tvphi)_i)_j)
-\Psi''(\tvphi)f_3(t;(\nablagz(\nablagz\tvphi)_i)_j)=0,
    \end{cases}
    \label{systf}
\end{align}
with the initial data $\u_0=\phi_{-t}(\V_0)_{\vert t=0}=\V_0$ and $\tvphi_0=\tphimt(\vphi_0)_{\vert t=0}=\vphi_0$.
\section{Short time existence of a regular solution to the pullback problem}\label{sec:shortpullback}
We now aim at showing that, if the initial datum is sufficiently regular and $\tvphi_0$ is strictly separated from pure phases, there exists a unique local-in-time regular solution to \eqref{systf}. This will be used to find a solution to \eqref{syst1}. In order to prove this result, 
We then need to introduce a suitable functional setting in the spirit of \cite{AWe}. 
Let us define the spaces ($\tT\leq T$)
\begin{align*}
    &Z^1_{\tT}:=L^2(0,\tT; \H^2(\Gamma_0))\cap W^{1,2}(0,\tT;\L^2_\sigma(\Gamma_0)),\\&
    Z_{\tT}^2:=L^q(0,\tT;W^{4,p}(\Gamma_0))\cap W^{1,q}(0,\tT;L^p(\Gamma_0)),\quad\text{ for some } 2<q\leq  4,\ p>4,
\end{align*}
equipped with the norms
\begin{align}
   & \Vert \f \Vert_{Z^1_{\tT}}:=\Vert \f \Vert_{L^2(0,\tT; \H^2(\Gamma_0))}+\Vert \f\Vert_{ W^{1,2}(0,\tT;\L^2_\sigma(\Gamma_0))}+\Vert \f(0) \Vert_{\H^1(\Gamma_0)},\\&
   \Vert f\Vert_{\Zdt}:=\Vert f \Vert_{L^q(0,\tT;W^{4,p}(\Gamma_0))}+\Vert f \Vert_{W^{1,q}(0,\tT;L^p(\Gamma_0))}+\Vert f(0) \Vert_{(L^p(\Gamma_0),W^{4,p}(\Gamma_0))_{ 1-\frac1q,q}}.
\end{align}
Note that, by \cite[Lemma 2]{Saal1}, we can immediately deduce that there exists $C(T)>0$ such that
\begin{align}
\Vert \f\Vert_{BUC([0,\tT];\H^1(\Gamma_0))}\leq C(T)\Vert \f\Vert_{\Zut},\quad \Vert f\Vert_{BUC([0,\tT];(L^p(\Gamma_0),W^{4,p}(\Gamma_0))_{1-\frac1q,q})}\leq C(T)\Vert f\Vert_{\Zdt},\quad \forall \tT\leq T.
    \label{uniform}
\end{align}
Moreover, we point out that, since $q>2$ and $p>4$, 
\begin{align}
\label{embedding}
(L^p(\Gamma_0),W^{4,p}(\Gamma_0))_{1-\frac 1q,q}=B_{p,q}^{4-\frac4q}(\Gamma_0)\hookrightarrow W^{2,p}(\Gamma_0)\hookrightarrow C^1(\Gamma_0), 
\end{align}
and, recalling the following Sobolev-Gagliardo-Nirenberg's inequality (which is valid also on $\Gamma_0$ by a standard partition of unity argument)
\begin{align*}
    \Vert f\Vert_{W^{3,p}(\Gamma_0)}\leq C\Vert f\Vert_{W^{4,p}(\Gamma_0)}^{\frac12}\Vert f\Vert_{W^{2,p}(\Gamma_0)}^\frac12,
\end{align*}
we infer from \eqref{embedding}, since $q>2$, 
\begin{align}
  \Vert f\Vert_{L^4(0,\tT;W^{3,p}(\Gamma_0))}\leq C\Vert f\Vert_{L^2(0,T;W^{4,p}(\Gamma_0))}^\frac12\Vert f\Vert_{L^\infty(0,\tT;W^{2,p}(\Gamma_0))}^\frac12\leq C\tT^\frac{q-2}{4q}\Vert f\Vert_{\Zdt}.  \label{fundamental}
\end{align}
The inequality above is the justification to the choice of the space $\Zdt$.
In conclusion, we also have the following interpolation and embedding results (see, for instance, \cite{Triebel}):
\begin{align*}
    (B^{4-\frac4q}_{p,q}(\Gamma_0),L^p(\Gamma_0))_{\theta,q}\hookrightarrow B^{s}_{p,q}(\Gamma_0)\hookrightarrow W^{2,p}(\Gamma_0),
\end{align*}
where $\theta\in(0,1)$ is such that $s:=(1-\theta)(4-\frac 4q)>2$.
Now, since $\Zdt\hookrightarrow W^{1,q}(0,\tT;L^p(\Gamma_0)\hookrightarrow C^{0,1-\frac 1q}([0,\tT];L^p(\Omega))$, by \cite[Lemma 1]{AWe} we also infer, by the interpolation result above and \eqref{uniform}-\eqref{embedding},
\begin{align}
    \label{holder}
    \Zdt\hookrightarrow C^{0,\theta(1-\frac1q)}([0,\tT];B_{p,q}^s(\Gamma_0))\hookrightarrow C^{0,\theta(1-\frac1q)}([0,\tT];W^{2,p}(\Gamma_0)).
\end{align}
We then introduce, for further use, the following spaces: 
$$
X^1_{\widetilde{T}}:=\{\v\in \Zut:\ \v_{|t=0}=\u_0\},\quad X^2_{\widetilde{T}}:=\{\vphi\in \Zdt:\ \vphi_{|t=0}=\tvphi_0\}, \quad \XT:= \XT^1\times \XT^2,
$$
as well as 
$$
Y^1_{\widetilde{T}}:=L^2(0,\widetilde{T};\L^2_\sigma(\Gamma_0)), \quad Y_{\widetilde{T}}^2=L^q(0,\widetilde{T};L^p(\Gamma_0)),\quad \YT:= Y_{\widetilde{T}}^1\times Y_{\widetilde{T}}^2.
$$
We can then state our main result of this section: 
\begin{theorem}
    Let $\Psi,\nu,\rho$  be as in Assumptions \ref{l1}-\ref{ro}. Moreover, let $\u_0\in \H^1_\sigma(\Gamma_0)$, and $\tvphi_0\in B_{pq}^{4-\frac4q}(\Gamma_0)$ for some $2<q\leq 4$ and $p>4$, be given
with $\norm{\tvphi_0}_{L^\infty(\gam_0)}\leq 1-2\delta_0$, for some $\delta_0>0$. Then, there exists $\tT>0$ such that
\eqref{systf} has a unique solution $(\u,\tvphi)$ such that 
$$
(\u,\tvphi) \in \XT^1\times \XT^2, 
$$
and 
\begin{align}
\Vert\tvphi(t)\Vert_{C^0(\Gamma_0)}\leq 1-\delta_0,\quad \forall t\in[0,\tT].
\label{strspr}
\end{align}
\label{thm1}
\end{theorem}

In order to prove this theorem, by means of a fixed point argument, we need some preliminary lemmas. First, the following lemma holds:
\begin{lemma}
\label{lemmaproj}
    It holds
    \begin{align}
    \label{est1}
    \Vert P(t)-P(s)\Vert_{\mathcal{L}(\L^2(\Gamma_0),\L^2_\sigma(\Gamma_0))}\leq C(T)\vert t-s\vert,\quad \forall t,s\in[0,T],
    \end{align}
    entailing, in particular, 
    \begin{align}
    \Vert P(t)\Vert_{\mathcal{L}(\L^2(\Gamma_0),\L^2_\sigma(\Gamma_0))}\leq C(T),\quad \forall t\in[0,T].
    \label{est0}
    \end{align}
\end{lemma}
\begin{remark}
    Estimate \eqref{est0} is a direct consequence of \eqref{est1}, recalling that $P(0)=P^{\Gamma_0}\in \mathcal{L}(\L^2(\Gamma_0),\L^2_\sigma(\Gamma_0))$.
\end{remark}
\begin{proof}
    The proof can be carried out with similar arguments to \cite[Lemma 3.1]{Saal1}. Let us first observe that, by the regularity of the flow map $\Phi^n_t$, we can easily deduce that the norm $\Vert \nabla_{\Gamma(t)}^\phi(t) f\Vert_{L^q(\Gamma_0)}$ (for any $t\in[0,T]$) and $\Vert \nablagz f\Vert_{L^q(\Gamma_0)}$ are equivalent for any $q>1$. In particular, there exists $C(T)>0$ such that 
    \begin{align}
        C(T)\Vert \nabla_{\Gamma(t)}^\phi(t) f\Vert_{L^q(\Gamma_0)}\leq \Vert \nablagz f\Vert_{L^q(\Gamma_0)}\leq C(T)\Vert \nabla_{\Gamma(t)}^\phi(t) f\Vert_{L^q(\Gamma_0)},\quad \forall t\in[0,T],\quad \forall f\in W^{1,q}(\gz).
        \label{equivalence}
    \end{align}
  Now, by the definition of $P(t)$, it holds that, for any $\v\in \L^2(\Gamma_0)$,
  $$
  P(t)\v=\v-\nabla_{\Gamma(t)}^\phi(t) p(t),
  $$
  where  $p(t)\in H^1(\Gamma_0)$ with $\overline{p}(t)=0$ (unique, thanks to the properties of $P^{\Gamma(t)}$) satisfies
    \begin{align}
\label{essential}
\int_{\Gamma_0}\nabla_{\Gamma(t)}^\phi(t) p(t) \cdot \nablagz \eta \ds=\int_{\Gamma_0}\v\cdot \nablagz\eta\ds,\quad \forall \eta\in H^1(\Gamma_0).
  \end{align}
Indeed, it holds, for any $t\in[0,T]$,
\begin{align*}
    &\int_{\Gamma_0}P(t)\v \cdot \nablagz \eta \ds= \int_{\Gamma_0}\tphimt(J\A^{-1}P^{\Gamma(t)}\tphit(\A\v))\cdot \nablagz \eta \ds\\&=
  \int_{\Gamma_0}\tphimt(J P^{\Gamma(t)}\tphit(\A\v))\cdot \nablagz \eta\tphimt(\Dm) \ds= \int_{\Gamma_0}\tphimt(J P^{\Gamma(t)}\tphit(\A\v))\cdot \tphimt(\nablag\tphit(\eta)) \ds\\&
  =\int_{\Gamma(t)}P^{\Gamma(t)}\tphit(\A\v)\cdot \nablag(\tphit(\eta))\ds=0,\quad \forall \eta \in H^1(\gz).
\end{align*}
Moreover, we can write
  \begin{align}
  \label{Pt}
  P(t)\v-P(s)\v=\nablagphi(t)p(t)-\nablagphis(s)p(s)=\nablagphi(t)(p(t)-p(s))+(\nablagphi(t)-\nablagphis(s))p(s).
  \end{align}
  Observe now that also, again by the regularity of $\Phi_t^n$ and the definition of $\phi_t$,
  \begin{align*}
      &\Vert \nablagphi(t)w\Vert\leq C(T)\Vert \nablag \tphi_{t}(w)\Vert\leq C(T)\sup_{\eta\in H^1(\Gamma(t)):\ \nablag\eta\not=0}\dfrac{\vert\ints{\Gamma}\nablag \tphi_t(w)\cdot\nablag\eta\ds\vert}{\Vert \nablag\eta\Vert},
  \end{align*}
  where in the last step we exploited the properties of the (orthogonal) \rev{Helmholtz} projector $P^{\Gamma(t)}$, so that, for any tangential $\mathbf{f}\in \L^{2}(\Gamma(t))$,
$$\ints{\Gamma}\nablag \tphi_t(w)\cdot \mathbf{f} \ds= \ints{\Gamma}\nablag \tphi_t(w)\cdot\nablag \eta\ds,\quad \text{with }\nablag\eta=\mathbf{f}-P^{\Gamma(t)}\mathbf{f}.
$$ 
We also have, by the definition of $\phi_t$ and $\nablagphi(t)$,
\begin{align*}
&{\ints{\Gamma}\nablag \tphi_t(w)\cdot\nablag\eta}\ds=\int_{\Gamma_0}\nablagz w\cdot \tphimt(J\Dm\D^{-T})\nablagz\tphimt(\eta)\ds\\&=\int_{\Gamma_0}\tphimt(J\Dm\D^{-T})\nablagz w\cdot \nablagz\tphimt(\eta)\ds   \\&=\int_{\Gamma_0}(\tphimt(\A^{-1}\D^{-T})\nablagz w)\cdot \nablagz\tphimt(\eta)\ds=\int_{\Gamma_0}\nablagphi(t) w\cdot \nablagz\tphimt(\eta)\ds.
\end{align*}
We can then conclude that
  \begin{align*}
      &\Vert \nablagphi(t)w\Vert\leq C(T)\sup_{\eta\in H^{1}(\Gamma_0):\ \nablagz\eta\not=0}\dfrac{\vert\int_{\Gamma_0}\nablagphi(t) w\cdot\nablagz\eta\ds\vert}{\Vert \nablagz\eta\Vert},
  \end{align*}
by exploiting the fact that $\Vert \nablag\eta\Vert\geq C(T)\Vert \nablagz\tphimt(\eta)\Vert$. As a consequence, we have 
\begin{align*}
   \Vert \nablagphi(t)(p(t)-p(s))\Vert&\leq C(T)\sup_{\eta\in H^1(\Gamma_0):\ \nablagz\eta\not=0}\dfrac{\vert\int_{\Gamma_0}\nablagphi(t) (p(t)-p(s))\cdot\nablagz\eta\ds\vert}{\Vert \nablagz\eta\Vert}\\&=C(T)\sup_{\eta\in H^1(\Gamma_0):\ \nablagz\eta\not=0}\dfrac{\vert\int_{\Gamma_0}(\nablagphi(t)-\nablagphis(s))p(s) \cdot\nablagz\eta\ds\vert}{\Vert \nablagz\eta\Vert}\\&\leq C\Vert (\nablagphi(t)-\nablagphis(s))p(s)\Vert, 
\end{align*}
where in the last step we exploited the fact that, due to \eqref{essential}, $$
\int_{\Gamma_0}(\nablagphi(t)p(t)-\nablagphis(s)p(s))\cdot\nablagz\eta\ds=0.
$$
Coming back to \eqref{Pt}, we thus see that we only need to control $\Vert (\nablagphi(t)-\nablagphis(s))p(s)\Vert$. In particular, thanks to the regularity of the flow map, we know that the maps $\A^{-1}$ and $\Dm$ belong to $C^{0,1}(\mathcal{S})$, entailing, by definition of $\nablagphi$, 
$$\Vert (\nablagphi(t)-\nablagphis(s))p(s)\Vert\leq C(T)\vert t-s\vert \Vert \nablagz p(s)\Vert.
$$
Since by \eqref{equivalence} and \eqref{essential}  we have
$$
\Vert p(t)\Vert_{H^1(\Gamma_0)}\leq C(T)\Vert \v\Vert, \quad \forall t\in[0,T],
$$
estimate \eqref{est1} is proven.
\end{proof} 
Secondly, we need to study a surface Stokes operator with variable viscosity. In particular, we have the following lemma, which is an adaptation of \cite[Lemma 4]{AARMA} to the case of compact surfaces. For the sake of generality, we state the Lemma in the case of  $\Gamma(t)$, $t\in[0,T]$ showing that the constants are uniform in time.
\begin{lemma}
\label{regularity}    Let  $\vphi  \in W^{1,\infty}(\Gamma(t) )$, and $\f(t) \in  \L^2(\Gamma(t))$, for almost any $t\in[0,T]$. Let $\omega>0$ be arbitrary. Then, there exists $C(T,\omega)>0$, depending on  $T$ and $\omega$, such that the unique $\u=\u(t)\in \H^2(\Gamma(t))$ satisfying, for almost any $t\in[0,T]$,
    \begin{align}
       &\nonumber2\ints{\Gamma}\nu(\varphi(t))\E_S(\u(t)):\E_S(\boldsymbol\eta)\ds+\omega\int_{\Gamma(t)}\u(t)\cdot\boldsymbol\eta\ds\\&\label{weakform}=\ints{\Gamma}\f(t)\cdot \boldsymbol\eta\ds,\quad \forall \boldsymbol\eta\in\H^1(\Gamma(t))\cap \L^2_\sigma(\Gamma(t)), 
    \end{align}
    is such that
\begin{align}
\Vert \u(t)\Vert_{\H^2(\Gamma(t))}\leq C(T,\omega)\left(1+\Vert\vphi(t)  \Vert_{W^{1,\infty}(\Gamma(t))}\right)\Vert \f(t)\Vert_{\L^2(\Gamma(t))},\quad \text{for a.a. } t\in[0,T].
    \label{estStokes}
\end{align}
\end{lemma}
\begin{remark}
    In the case $t=0$, it can be easily seen from the proof that, if $\vphi\in W^{1,\infty}(\Gamma_0)$ and $\f\in \L^2(\Gamma_0)$, the constant $C>0$ in \eqref{estStokes} does not depend on $T$ but only on $\omega>0$. Indeed, the operator $\phi_t$ at $t=0$ is the identity on $\Gamma_0$.
\end{remark}
\begin{proof}
 First we observe that, by the Lax-Milgram Lemma (recalling the norm equivalence between the quadratic form associated to this problem and the $\H^1(\Gamma(t))$ norm in \cite[(3.11)]{V2}), it is immediate to show that there exists a unique solution $\u(t)\in \H^1(\Gamma(t))\cap \L_\sigma^2(\Gamma(t))$ to \eqref{weakform}, for almost any $t\in[0,T]$. Moreover, it holds, for some $C(T)>0$, 
 \begin{align}
  C(T)\Vert \u(t)\Vert_{\H^1(\Gamma(t))}\leq 2\nu_*\Vert \E_S(\u(t))\Vert^2_{\L^2(\Gamma(t))}+\omega\Vert \u(t)\Vert_{\L^2(\Gamma(t))}^2\leq \Vert \f(t)\Vert^2_{\L^2(\Gamma(t))},\quad \text{for a.a. }t\in[0,T].\label{contest}
\end{align}
    Note that $C(T)$ is independent of $T$ if we only consider the case $t=0$. Now we first consider the case when $\nu=\nu_*>0$ is a constant. In this case we know from \cite[Corollary 3.4]{Simonett} that the Stokes operator $\A_\sigma^{\Gamma(t)}:=-2\nu_* P^{\Gamma(t)}\P\text{div}_{\Gamma(t)}(\E_S(\cdot))$ when $t=t_0\geq0$ is fixed has $L^p$-maximal regularity, and thus the operator  $\A_\sigma^{\Gamma(t), \omega}:=\A_\sigma^{\Gamma(t)}+\omega: \mathfrak{D}(\A_\sigma^{\Gamma(t)})=\H^2(\Gamma(t))\cap \L^2_\sigma(\Gamma(t))\to \L^2_\sigma(\Gamma(t))$ is an isomorphism, for for any $t\in[0,T]$ and any $\omega>0$.  Notice that $\omega>0$ can be fixed as to be the same for any time since the spectral bound of the operator $\A_\sigma^{\Gamma(t)}$ is nonnegative for any $t\in[0,T]$ (see \cite[Proposition 4.1]{Simonett}, see also \cite[Remark 32]{Simonett2}), and thus, for a fixed $\omega>0$, $\A_\sigma^{\omega,\Gamma(t)}$ is an isomorphism for any $t\in[0,T]$ (clearly the operator norm depends on $\omega>0$ and potentially diverges as $\omega\to 0$). Now we introduce the composed operator 
    $$
\A^{\omega}_\sigma(t)=\phi_{-t}\A_\sigma^{\Gamma(t),\omega}\phi_t.
$$ 
Thanks to the regularity of the map $\phi_t$, the operators  $
\A^{\omega}_\sigma(t)
$ are still isomorphisms from $\mathfrak{D}(\A_\sigma^{\Gamma_0})$ to $\L^2_\sigma(\Gamma_0)$ for any $t\in[0,T]$. We now need to show that $\A^{\omega}_\sigma(t)$ are continuous in the space $\mathcal{L}(\mathfrak{D}(\A_\sigma^{\Gamma_0}),\L^2_\sigma(\Gamma_0))$. We can use Lemma \ref{lemmaproj}. In particular, observe that, in the notation of Appendix \ref{pullback1}, it holds 
\begin{align}
\A_\sigma^\omega(t)\u=P(t)\P_0\a_0(t;\nablagz((\nablagz\u_i)_j)_k)+P(t)\a_{1}(t;(\nablagz\u_i)_j)+P(t)\a_{2}(t)\u+P(t)a_3(t)\u.
\end{align}
Therefore by Lemma \ref{lemmaproj} and the regularity of the flow map we have 
\begin{align*}
   & \Vert \A_\sigma^\omega(t)-\A_\sigma^\omega(s)\Vert_{\mathcal{L}(\mathfrak{D}(\A_\sigma^{\Gamma_0}),\L^2_\sigma(\Gamma_0))}\\&\leq \Vert P(t)-P(s)
 \Vert_{\mathcal{L}(\L^2(\Gamma_0),\L^2_\sigma(\Gamma_0))}\Vert \P_0\a_0(t;\nablagz((\nablagz\u_i)_j)_k)+\a_1(t;(\nablagz\u_i)_j)+\a_2(t)\u+a_3(t)\u\Vert_{\L^2(\Gamma_0)}\\&+\Vert P(s)
 \Vert_{\mathcal{L}(\L^2(\Gamma_0),\L^2_\sigma(\Gamma_0))}(\Vert \P_0\a_0(t;\nablagz((\nablagz\u_i)_j)_k)-\P_0\a_0(s;\nablagz((\nablagz\u_i)_j)_k)\Vert_{\L^2(\Gamma_0)}\\&+\Vert \a_1(t;(\nablagz\u_i)_j)-\a_1(s;(\nablagz\u_i)_j)\Vert_{\L^2(\Gamma_0)}\\&+\Vert (\a_2(t)-\a_2(s))\u\Vert_{\L^2(\Gamma_0)}+\Vert (a_3(t)-a_3(s))\u\Vert_{\L^2(\Gamma_0)})\leq C(T)\vert t-s\vert\Vert \u\Vert_{\H^2(\Gamma_0)},
\end{align*}
 which entails the desired continuity on $[0,T]$. This means that it holds
 \begin{align*}
\A_\sigma^\omega(t)\in C^{0,1}{([0,T];\mathcal{L}(\mathfrak{D}(\A_\sigma^{\Gamma_0}),\L^2_\sigma(\Gamma_0)))}.
 \end{align*}
Now it is immediate to verify that, for any fixed $t_0\in[0,T]$, it holds 
\begin{align*}
    \Vert (\A_\sigma^\omega(t_0))^{-1}\A_\sigma^\omega(t)-Id\Vert_{\mathcal{L}(\L^2_\sigma(\Gamma_0), \L^2_\sigma(\Gamma_0)))}&\leq \Vert (\A_\sigma^\omega(t_0))^{-1}\Vert_{\mathcal{L}(\L^2_\sigma(\Gamma_0), \mathfrak{D}(\A_\sigma^{\Gamma_0}))}\Vert \A_\sigma^\omega(t)-\A_\sigma^\omega(t_0)\Vert_{\mathcal{L}( \mathfrak{D}(\A_\sigma^{\Gamma_0}),\L^2_\sigma(\Gamma_0))}\\&\leq C(t_0,T,\omega)\vert t-t_0\vert,\quad \forall t\in[0,T],
\end{align*}
where $Id$ is the identity operator in $\L^2_\sigma(\Gamma_0)$, so that there exists $\delta(t_0)>0$ such that 
 \begin{align}
   &\Vert (\A_\sigma^\omega(t_0))^{-1}\A_\sigma^\omega(t)-Id\Vert_{\mathcal{L}(\L^2_\sigma(\Gamma_0), \L^2_\sigma(\Gamma_0)))}\leq \frac{1}{2},\quad \forall t\in[t_0-\delta(t_0),t_0+\delta(t_0)].
    \label{deltai}
\end{align}
This allows to use the Neumann series (see also \cite[Lemma 3.4]{Saal1}) to deduce that 
\begin{align}
     &\nonumber\Vert (\A_\sigma^\omega(t))^{-1}\Vert_{\mathcal{L}(\L^2_\sigma(\Gamma_0), \mathfrak{D}(\A_\sigma^{\Gamma_0}))}\\&\leq\dfrac{\Vert \A_\sigma^\omega(t_0))^{-1}\Vert_{\mathcal{L}(\L^2_\sigma(\Gamma_0), \mathfrak{D}(\A_\sigma^{\Gamma_0})))}}{1-\Vert (\A_\sigma^\omega(t_0))^{-1}\A_\sigma^\omega(t)-Id\Vert_{\mathcal{L}(\L^2_\sigma(\Gamma_0), \L^2_\sigma(\Gamma_0)))}}\leq C(t_0,T,\omega),\quad \forall t\in[t_0-\delta(t_0),t_0+\delta(t_0)]. \label{deltai2}
\end{align}
Let us consider a cover of $[0,T]$ made of balls of center in some $t_0^i$ and radius the corresponding $\delta_0^i$ such that \eqref{deltai2} holds. Since $[0,T]$ is compact, these balls can be chosen in finite number, say $m(T)$. By defining then $C_0(T,\omega):=\max_{i=1,\dots,m(T)}\{C(t_0^i,T,\omega)\}$, we infer that 
\begin{align}
     \Vert (\A_\sigma^\omega(t))^{-1}\Vert_{\mathcal{L}(\L^2_\sigma(\Gamma_0), \mathfrak{D}(\A_\sigma^{\Gamma_0}))}\leq C_0(T,\omega),\quad \forall t\in[0,T]. \label{deltai3}
\end{align}

    By the regularity of the map $\phi_t$, since it holds $\A_\sigma^{\Gamma(t),\omega}=\phi_t(\A_\sigma^\omega(t))^{-1}\phi_{-t}$, we immediately obtain
\begin{align*}
     \Vert (\A_\sigma^{\Gamma(t),\omega(t)})^{-1}\Vert_{\mathcal{L}(\L^2_\sigma(\Gamma(t)), \mathfrak{D}(\A_\sigma^{\Gamma(t)}))}\leq C(T,\omega),\quad \forall t\in[0,T],
\end{align*}
    from which we obtain, coming back to the solution $\u(t)$ satisfying \eqref{contest} in the case of constant $\nu=\nu_*$,  that it also verifies
    \begin{align}
        \Vert \u(t)\Vert_{\H^2(\Gamma(t))}\leq C(T,\omega)\Vert \f\Vert_{L^2(\Gamma(t))},\quad \text{for a.a. }t\in [0,T].\label{nuconst}
    \end{align}
    Let us now pass to consider the case of variable viscosity $\nu$.
First we need to study the problem for the pressure: given $\v(t)\in\L^2(\Gamma(t))$ for a fixed $t\in[0,T]$, find $p$ such that (in weak formulation)
$$
-\Delta_\Gamma p(t)=-\divg(\v(t)).
$$
By standard Lax-Milgram lemma, there exists a unique weak solution $p(t)\in H^1(\Gamma(t))$ for almost any $t\in[0,T]$, such that $\overline{p}(t)=0$ and
$$
\Vert p(t)\Vert_{H^1(\Gamma(t))} \leq C(T)\Vert \v(t)\Vert_{\L^2(\Gamma(t))},
$$
where $C(T)$ does not depend on $t\in[0,T]$.
Moreover, if $\divg(\v(t))\in L^2(\Gamma(t))$, the solution $p(t)$ is a strong solution and it holds (see, for instance, \cite[Theorem 3.3]{DzEll})
\begin{align}
    \Vert p(t)\Vert_{H^2(\Gamma(t))}\leq C(T)\Vert \divg(\v(t))\Vert_{L^2(\Gamma(t))}.
    \label{press}
\end{align}
Notice that the fact that the constant $C(T)$ does not depend on the specific $t\in[0,T]$ can be obtained arguing in a very similar way as to infer \eqref{nuconst}. We denote the surface gradient of the solution $p(t)$ by $$\Ng(\v(t)):=\nablag p(t).$$ We now consider the weak solution $\u(t)$ to the problem \eqref{weakform}, with nonconstant $\nu$, which exists uniquely in $\H^1(\Gamma(t))$ by the arguments above, and satisfies, for the fixed $t\in[0,T]$, equation \eqref{contest}. Now, let us take $\boldsymbol\eta=\frac{\mathbf{w}}{2\nu(\vphi(t))}-\Ng\left(\frac{\mathbf{w}}{2\nu(\vphi(t))}\right)$, for some $\w\in C^\infty(\Gamma(t))$ such that $\divg\w=0$. By construction it holds $\divg(\boldsymbol\eta)=0$, and thus $\boldsymbol\eta$ can be used as a test function in \eqref{weakform}. This leads to the following 
\begin{align}
&\nonumber\ints{\Gamma}\E_S(\u(t)):\E_S(\w)\ds+\omega\ints{\Gamma}\u(t)\cdot\w\ds\\&\nonumber=\ints{\Gamma}\dfrac{\nu'(\vphi(t))}{\nu(\vphi(t))}\E_S(\u(t)):(\w\otimes \nablag\vphi(t))\nonumber\ds\\&\nonumber+\ints{\Gamma}\nu(\vphi(t))\E_S(\u(t)):\E_S\left(\Ng\left(\frac{\mathbf{w}}{\nu(\vphi(t))}\right)\right)\ds\\&+\omega\ints{\Gamma}\u(t)\cdot\w\left(1-\frac {1}{2\nu(\vphi(t))}\right)\ds\nonumber\\&+\omega\ints{\Gamma}\u(t) \cdot \Ng\left(\frac{\mathbf{w}}{2\nu(\vphi(t))}\right)\ds\nonumber\\&+\ints{\Gamma}\left(\frac{\mathbf{w}}{2\nu(\vphi(t))}-\Ng\left(\frac{\mathbf{w}}{2\nu(\vphi(t))}\right)\right)\cdot \f(t)\ds.
    \label{comp}
\end{align}
Observe that, since $\nu\in W^{1,\infty}(\R)$ and by \eqref{contest}, 
\begin{align*}
&\left\vert\ints{\Gamma}\dfrac{\nu'(\vphi(t))}{\nu(\vphi(t))}\E_S(\u(t)):(\w\otimes \nablag\vphi(t))\ds\right\vert\\&\leq C\Vert \u\Vert_{\H^1(\Gamma(t))}\Vert \nablag\vphi(t)\Vert_{\L^\infty(\Gamma(t))}\Vert \w\Vert_{\L^2(\Gamma(t))}\\&\leq C(T)\Vert \f(t)\Vert_{\L^2(\Gamma(t))}\Vert \nablag\vphi(t)\Vert_{\L^\infty(\Gamma(t))}\Vert \w\Vert_{\L^2(\Gamma(t))}. 
\end{align*}
Furthermore, it holds, recalling \eqref{contest}, \eqref{press}, and $\divg\w=0$,
\begin{align*}
&\ints{\Gamma}\nu(\vphi(t))\E_S(\u(t)):\E_S\left(\Ng\left(\frac{\mathbf{w}}{\nu(\vphi(t))}\right)\right)\ds\\&\leq C\Vert \u(t)\Vert_{\H^1(\Gamma(t))}\left\Vert \Ng\left(\frac{\mathbf{w}}{\nu(\vphi(t))}\right)\right\Vert_{\H^1(\Gamma(t))}\\&
\leq C(T)\Vert \f\Vert_{\L^2(\Gamma(t))}\left\Vert \divg\left(\frac{\mathbf{w}}{\nu(\vphi(t))}\right)\right\Vert_{\L^2(\Gamma(t))}\\&
\leq C(T)\Vert \f\Vert_{\L^2(\Gamma(t))}\left\Vert\nablag\vphi(t)\right\Vert_{\L^\infty(\Gamma(t))}\norm{\w}_{\L^2(\Gamma(t))}.
\end{align*}
Similarly, we have,
\begin{align*}
&\omega\ints{\Gamma}\u(t)\cdot\w(1-\frac {1}{2\nu(\vphi(t))})+\omega\ints{\Gamma}\u(t) \cdot \Ng\left(\frac{\mathbf{w}}{2\nu(\vphi(t))}\right)\ds\\&+\ints{\Gamma}\left(\frac{\mathbf{w}}{2\nu(\vphi(t))}-\Ng\left(\frac{\mathbf{w}}{2\nu(\vphi(t))}\right)\right)\cdot \f(t)\ds\\&
    \leq C(T)\Vert \f\Vert_{\L^2(\Gamma(t))}\Vert \w\Vert_{\L^2(\Gamma(t))}.
\end{align*}
In conclusion, by introducing the vector field $\mathbf{F}(t)\in\L^2_\sigma(\Gamma(t))$ defined, for almost any $t\in[0,T]$, by 
\begin{align*}
\int_{\Gamma(t)}\mathbf{F}(t)\cdot \w\ds&=\ints{\Gamma}\dfrac{\nu'(\vphi(t))}{\nu(\vphi(t))}\E_S(\u(t))(\w\otimes \nablag\vphi(t)+\nablag\vphi(t)\otimes\w)\ds\\&+\ints{\Gamma}\nu(\vphi(t))\E_S(\u(t)):\E_S\left(\Ng\left(\frac{\mathbf{w}}{\nu(\vphi(t))}\right)\right)\ds+\omega\ints{\Gamma}\u(t)\cdot\w(1-\frac {1}{2\nu(\vphi(t))})\ds\nonumber\\&+\omega\ints{\Gamma}\u(t) \cdot \Ng\left(\frac{\mathbf{w}}{2\nu(\vphi(t))}\right)\ds\nonumber+\ints{\Gamma}\left(\frac{\mathbf{w}}{2\nu(\vphi(t))}-\Ng\left(\frac{\mathbf{w}}{2\nu(\vphi(t))}\right)\right)\cdot \f(t)\ds,
\end{align*}
for any $\w\in \L^2_\sigma(\Gamma(t))$, we obtain, recalling $W^{1,\infty}(\Gamma(t))\hookrightarrow H^1(\Gamma(t))$ (where the embedding constant can be chosen depending only on $T$), from the estimates above, by a density argument, that 
$$
\Vert \mathbf{F}(t)\Vert_{\L^2(\Gamma(t))}\leq C(T)(1+\Vert \vphi(t)\Vert_{W^{1,\infty}(\Gamma(t))})\Vert \f(t)\Vert_{\L^2(\Gamma(t))}.
$$
Therefore, $\u(t)$ satisfies 
$$
2\ints{\Gamma}\E_S(\u(t)):\E_S(\w)+\omega\ints{\Gamma}\u(t)\cdot\w\ds=\ints{\Gamma}\mathbf{F}(t)\cdot \w\ds,\quad \forall \w\in \H^1(\Gamma(t))\cap\L^2_\sigma(\Gamma(t)).
$$
We can then apply the regularity result \eqref{nuconst} (with $\nu_*=1$), to infer 
$$
\Vert \u(t)\Vert_{\H^2(\Gamma(t))}\leq C(T,\omega)\Vert \mathbf{F}(t)\Vert_{\L^2(\Gamma(t))}. 
$$
Since the same result, with the same constant $C(T,\omega)$, holds for almost any $t\in[0,T]$, this allows to deduce \eqref{estStokes}, thus concluding the proof. 
\end{proof}
We are now ready to start the proof of Theorem \ref{thm1}. To this aim we rewrite system \eqref{systf} as follows 
\begin{align*}
        \mathcal{L}(\u,\tvphi)=\mathcal{F}(\u,\tvphi),
\end{align*}
where we defined the  linear operator $\mathcal{L}:\XT\to\YT$ as 
\begin{align}
    \label{L}
    \mathcal{L}(\v,\vphi):=\begin{bmatrix}
        P^{\Gamma_0}\tilde{\rho}_0\partial_t\v-2P^{\Gamma_0}\P_0\mathrm{div}_{\Gamma_0}(\nu(\tvphi_0)\EE_S(\v))+\omega \v\\
        \partial_t\vphi+\Delta_{\Gamma_0}^2\vphi 
    \end{bmatrix},
\end{align}
for some arbitrary $\omega>0$, where we set $\tilde{\rho}_0=\rho(\widetilde{\vphi}_0)$. Note that $\Delta^2_{\gz}$ has $L^q-L^p$ maximal regularity on the finite interval $[0,T]$ (see, e.g., \cite[Theorem 6.4.3 (ii)]{SP}). We also define the nonlinear operator $\mathcal{F}:\XT\to\YT$ as
\begin{align} 
    \label{F}
    \mathcal{F}(\v,\vphi):=\begin{bmatrix}
       -P(t)\rho\partial_t\v
    +P^{\Gamma_0}\tilde{\rho}_0\partial_t\v+\omega\v-P(t)\A_\ast(t;\vphi,\nablagz\vphi)(\v,\nablagz\v)
\\-P(t)\b_3(t)\v\cdot(\di_1(t;((\nablagz(\nablagz(\nablagz \vphi)_i)_j)_k))+\di_2(t; ((\nablagz(\nablagz \vphi)_i)_j)))\\-P(t)(\b_6(t)(\di_1(t;((\nablagz(\nablagz(\nablagz \vphi)_i)_j)_k))+\di_2(t; ((\nablagz(\nablagz \vphi)_i)_j)))\cdot \nablagz)\v
\\+P(t)\nu(\vphi)\P_0\a_0(t;\nablagz((\nablagz\v_i)_j)_k)+P(t)\nu(\vphi)\a_1(t)\v\\+P(t)\nu(\vphi)a_2(t)\u-2P^{\Gamma_0}\P_0\mathrm{div}_{\Gamma_0}(\nu(\tvphi_0)\EE_S(\v))
\\-P(t)\b_7(t)(\di_1(t;((\nablagz(\nablagz(\nablagz \vphi)_i)_j)_k))+\di_2(t; ((\nablagz(\nablagz \vphi)_i)_j)))
\\ \\
-B(t;\vphi,\nablagz\vphi)(\v)+\Delta_{\gz}^2\vphi-f_0(t;(\nablagz(\nablagz(\nablagz(\nablagz \vphi)_i)_j)_k)_l)\\-f_1(t;((\nablagz(\nablagz(\nablagz \vphi)_i)_j)_k))-f_2(t; ((\nablagz(\nablagz \vphi)_i)_j)\\
+\Psi''(\vphi)f_3(t;(\nablagz(\nablagz\vphi)_i)_j)
    \end{bmatrix}.
\end{align}

We now state the following propositions, but we postpone their proofs to Sections \ref{l1p}-\ref{l2p} below.  
\begin{proposition}
    Let the assumptions of Theorem \ref{thm1} hold true. Then
there is a constant $C(T,\widetilde{T}, R)>0$ such that
$$\Vert\mathcal{F}(\v_1, \vphi_1) -\mathcal{F}(\v_2, \vphi_2)\Vert_{\YT}
\leq C(T,\widetilde{T}, R)\Vert(\v_1 - \v_2, \vphi_1 -\vphi_2)\Vert_{\XT},$$
for all $(\v_i, \vphi_i ) \in \XT$ with $\Vert(\v_i, \vphi_i )\Vert_{\XT}
\leq R$, $R>0$, and $i = 1, 2$. Furthermore it holds
$C(T,\widetilde{T}, R) \to0$ as $\widetilde{T} \to0$.
\label{prop1}
\end{proposition}
\begin{proposition}
   \label{prop2}
   Let $\widetilde{T} > 0$ and $\mathcal{L}$ as defined above. Then, $\mathcal{L}: \XT \to \YT$
is an isomoprhism. Moreover, there is a constant $C(T) > 0$ such that
$\Vert\mathcal{L}^{-1}\Vert_{\mathcal{L}(\YT ,\XT )} \leq C(T)$ for any $\widetilde{T} \in(0, T]$.
\end{proposition}
\subsection{Proof of Theorem \ref{thm1}}
With  Propositions \ref{prop1}-\ref{prop2} at hand, which correspond to \cite[Propostion 1, Theorem 4]{AWe}, Theorem \ref{thm1} is proved by following \textit{verbatim} the same fixed point argument (contraction mapping theorem) as in the proof of \cite[Theorem 2]{AWe}, where $T$ is here substituted by $\widetilde{T}$. Here the fixed point argument is applied to the problem satisfied by a solution $(\u,\tvphi)\in \XT$, for some $\tT\in(0,T]$, i.e.,
$$
(\u,\tvphi)=\mathcal{L}^{-1}\mathcal{F}(\u,\tvphi) \quad\text{in }\XT.
$$
First, one fixes $R$ so that, for some $(\overline{\v},\overline{\vphi})\in \XT$ it holds
$(\overline{\v},\overline{\vphi})\in\overline{B}_{ R}^{\XT}(0)$, where  $\overline{B}_R^{\XT}(0)$ is the closed ball of $\XT$ of radius $R$ and centered at zero. Moreover, we choose $R>2\norm{\mathcal{L}^{-1}\mathcal{F}(\overline{\v},\overline{\vphi})}_{\XT}$. Then one fixes $0<\tT\leq T$ (possibly depending also on $R$) such that the operator $\mathcal{L}^{-1}\mathcal{F}$ is, for instance, a $\tfrac12$-contraction mapping from $\XT$ to $\XT$. This is possible thanks to Propositions \ref{prop1} and \ref{prop2}. Then, thanks to the same propositions, one shows that $\mathcal{L}^{-1}\mathcal{F}$ is well defined from $\overline{B}_R^{\XT}(0)$ to itself, as long as we choose a (possibly) smaller $0<\tT\leq T$ (also depending on $R$). In conclusion, to apply Banach fixed point theorem one has to show that $\mathcal{L}^{-1}\mathcal{F}$ is a contraction, and this is again possible thanks to  Propositions \ref{prop1} and \ref{prop2}. Therefore, there exists a solution $(\u,\tvphi)\in\XT$ to the problem under study. By a standard argument it is also easy to show that the solution  $(\u,\tvphi)\in\overline{B}_R^{\XT}(0)\subset\XT$ we just found is unique in $\XT$. Indeed, assume that there exists another solution $(\v,\vphi)\in  \XT$. Then consider $\tilde{R}>0$ larger than $R$ used in the previous argument (otherwise it would be obvious from the uniqueness given in the contraction mapping theorem), so that $(\v,\vphi)\in\overline{B}_{\tilde{R}^{\XT}(0)}$. Then repeating the same reasoning we deduce by the contraction mapping theorem that there exists $\tT_1(\tilde{R})\in(0,\tT]$ such that the solution in $\overline{B}_{\tilde{R}^{X_{\tT_1}}(0)
}$ is unique on $[0,\tT_1(\tilde{R})]$. Therefore, it is immediate to infer $(\v,\vphi)_{\vert [0,\tT_1]}=(\u,\tvphi)_{\vert [0,\tT_1]}$. A continuation argument then allows to show that the identity holds on the entire interval $[0,\tT]$.  

The additional strict separation property \eqref{strspr} can be proven as follows. Observe that, since the solution $\tvphi\in \XT^2$, from \eqref{holder} and the embedding $W^{2,p}(\gz)\hookrightarrow C^{0,1}(\gz)$ we infer that $\tvphi\in C^{0,\thetaq}(\gz\times[0,\tT])$. Therefore, if the initial datum is strictly separated, i.e., there exists $\delta_0$ such that $\norm{\tvphi_0}_{C_0(\gz)}\leq 1-2\delta_0$, then there exists a $T_*\in(0,\tT]$, possibly smaller than $\tT$, such that \eqref{strspr} holds by continuity. The proof of Theorem \ref{thm1} is complete.
\subsection{Proof of Proposition \ref{prop1}}
\label{l1p}
We need to compute the quantity $\Vert\mathcal F(\v_1,\vphi_1)-\mathcal F(\v_2,\vphi_2)\Vert_{\YT}$. We now analyze term by term this quantity. Let us start with the first row. First, we have, indicating $\trho(\vphi_i)$ with $\trho_i$, $i=1,2$ and recalling \eqref{holder}, \eqref{est1}-\eqref{est0}, the embedding $W^{2,p}(\Gamma_0)\hookrightarrow C^1(\Gamma_0)$ ($p>4$ by assumption), $\vphi_1(0)=\vphi_2(0)=\tvphi_0$, and $P(0)=P^{\Gamma_0}$,  
\begin{align*}
&\norm{\Pt(\trho_1\dtn\v_1-\trho_2\dtn\v_2)-P^{\Gamma_0}\trho_0\dtn(\v_1-\v_2)}_{\LQs}\\& \leq \norm{(P(t)-P(0))\trho_1(\dtn\v_1-\dtn\v_2)}_{\LQs}\\&
+\norm{P(t)(\trho_1-\trho_0)(\dtn\v_1-\dtn\v_2)}_{\LQs}+C(T)\norm{(\trho_1-\trho_2)\dtn\v_2}_{\LQs}\\&
\leq \CT \tT\norm{\v_1-\v_2}_{\XT^1}+
\CT\norm{\vphi_1-\vphi_1(0)}_{L^\infty(0,T;L^\infty(\Gamma))}\norm{\dtn\v_1-\dtn\v_2}_{\LQ}\\&
+C(T)\norm{(\vphi_1-\vphi_2)-(\vphi_1(0)-\vphi_2(0))}_{L^\infty(0,T;L^\infty(\Gamma_0))}\norm{\dtn\v_2}_{\LQ}
\\&\leq \CT \tT\norm{\v_1-\v_2}_{\XT^1}+\CT\tT^{\thetaq}\norm{\v_1-\v_2}_{\XT^1}+\CT R\tT^{\thetaq}\norm{\vphi_1-\vphi_2}_{\XT^2}\\&
\leq C(T,\widetilde{T}, R)\Vert(\v_1 - \v_2, \vphi_1 -\vphi_2)\Vert_{\XT},
\end{align*}
where here and in the sequel, $ C(T,\widetilde{T}, R)$ is increasing with $R$ and goes to zero as $\widetilde{T}\to 0^+$. Proceeding in the estimates, we recall the decomposition \eqref{decomposition} of the operator $\A_\ast$ in $\A_i$, $i=0,\ldots,4$, provided in Appendix \ref{pullback1}. First, we have, denoting by $D\A_{1}$ the Jacobian of $\A_{1}$ with respect to the variables $\vphi,\nablagz\vphi$, which is clearly bounded in $L^\infty(0,\tT;L^\infty(\Gamma_0))$ by the assumptions on the flow map, the continuity of $\Psi''',\Psi'',\nu',\nu'',\trho$ and since, by \eqref{holder}, $\norm{\vphi_i}_{W^{1,\infty}(\Gamma_0)}\leq CR$, $i=1,2$, 
\begin{align*}
    &\norm{P(t)(\A_{1}(t;\vphi_1,\nablagz\vphi_1)(\v_1)-\A_{1}(t;\vphi_2,\nablagz\vphi_2)(\v_2))+\omega(\v_1-\v_2)}_{\LQs}\\&\leq 
    C(T)(1+\norm{D\A_{1}}_{L^\infty((0,\tT)\times [-CR,CR]\times [-CR,CR]^3)}\norm{(\vphi_1-\vphi_2,\nablagz(\vphi_1-\vphi_2))}_{L^\infty(0,\tT;\L^\infty(\Gamma_0))})\norm{\v_1}_{\LQ}\\&
    +\CT(1+\norm{\A_{1}(t;\vphi_2,\nablagz\vphi_2)}_{L^\infty(0,\tT;\L^\infty(\Gamma_0))})\norm{\v_1-\v_2}_{\LQ}\\&
    \leq C(T)\tT^{\thetaq+\frac12}
\norm{\v_1}_{L^\infty(0,\tT;\L^2(\Gamma_0))}\norm{\vphi_1-\vphi_2}_{\XT^2}+\CT R\tT^\frac12\norm{\v_1-\v_2}_{L^\infty(0,\tT;\L^2(\Gamma_0))}\\&
\leq C(T,\widetilde{T}, R)\Vert(\v_1 - \v_2, \vphi_1 -\vphi_2)\Vert_{\XT}.
\end{align*}
Here we also used $(\vphi_1-\vphi_2,\nablagz(\vphi_1-\vphi_2))=(\vphi_1-\vphi_2-(\vphi_1-\vphi_2)(0),\nablagz(\vphi_1-\vphi_2)-\nablagz(\vphi_1-\vphi_2)(0))$.
Note that we have set $$\norm{\A_{1}(t;\vphi_2,\nablagz\vphi_2)}_{L^\infty(0,\tT;\L^\infty(\Gamma_0))}:=\sup_{(x,t)\in\Gamma_0 \times [0,T]}\sup_{\u\in \R^3,\ \u\not=\mathbf 0}\frac{\norma{\A_{1}(t;\vphi_2(x,t),\nablagz\vphi_2(x,t))(\u)}}{\norma{\u}},$$ 
as well as 
$$
\norm{D\A_{1}}_{L^\infty((0,\tT)\times [-CR,CR]\times [-CR,CR]^3)}:=\sup_{(t,y,\mathbf z)\in(0,\tT)\times [-CR,CR]\times [-CR,CR]^3}\sup_{\u\in \R^3,\ \u\not=\mathbf 0}\frac{\norma{D\A_{1}(t;y,\mathbf z)(\u)}}{\norma{\u}}.
$$
We can define analogously, \textit{mutatis mutandis}, all the following quantities related to the operators $\A_i$.
In a completely analogous way, recalling $\v_i\in BUC([0,\tT];\H^1(\gz))$ by \eqref{uniform}, we obtain 
\begin{align*}
    &\norm{P(t)(\A_{2}(t;\vphi_1,\nablagz\vphi_1)(\nablagz\v_1)-\A_{2}(t;\vphi_2,\nablagz\vphi_2)(\nablagz\v_2))}_{\LQs}\\&\
\leq C(T,\widetilde{T}, R)\Vert(\v_1 - \v_2, \vphi_1 -\vphi_2)\Vert_{\XT}.
\end{align*}
Then, concerning $\A_0$, we have again the same result, namely
\begin{align*}
    &\norm{P(t)(\A_{0}(t;\vphi_1,\nablagz\vphi_1)-\A_{0}(t;\vphi_2,\nablagz\vphi_2))}_{\LQs}\\&\leq 
    C(T)\tT^\frac 12\norm{D\A_{0}}_{L^\infty((0,\tT)\times [-CR,CR]\times [-CR,CR]^2)}\norm{(\vphi_1-\vphi_2,\nablagz(\vphi_1-\vphi_2))}_{L^\infty(0,\tT;\L^\infty(\Gamma_0))}
  \\& 
\leq C(T,\widetilde{T}, R)\Vert(\v_1 - \v_2, \vphi_1 -\vphi_2)\Vert_{\XT}.
\end{align*}
Then it is straightforward to show that, about $\A_3$, recalling the regularity \eqref{uniform} for $\v_i$ and the fact that $\XT^1\hookrightarrow L^4(0,\tT;\mathbf{W}^{1,4}(\gz))$, 
\begin{align*}
    &\norm{P(t)(\A_3(t;\vphi_1,\nablagz\vphi_1)(\v_1,\v_1)-\A_3(t;\vphi_2,\nablagz\vphi_2)(\v_2,\v_2))}_{\LQs}\\&
    \leq C(T) \norm{(\vphi_1-\vphi_2,\nablagz(\vphi_1-\vphi_2))}_{L^\infty(0,\tT;\L^\infty(\Gamma_0))}\norm{\v_1}_{L^4(0,\tT;\L^4(\gz))}\norm{\v_1}_{L^4(0,\tT;\L^4(\gz))}\\&+
    C(T)\norm{\v_1-\v_2}_{L^4(0,\tT;\L^4(\gz))}(\norm{\v_1}_{L^4(0,\tT;\L^4(\gz))}+\norm{\v_2}_{L^4(0,\tT;\L^4(\gz))})\\&\leq C(T,\widetilde{T}, R)\Vert(\v_1 - \v_2, \vphi_1 -\vphi_2)\Vert_{\XT}.
\end{align*}
In conclusion, for $\A_4$, we have, recalling the embeddings $\H^2(\gz)\hookrightarrow \L^\infty(\gz)$ and $\H^1(\gz)\hookrightarrow\L^4(\gz)$, 
\begin{align*}
    &\norm{P(t)(\A_4(t;\vphi_1,\nablagz\vphi_1)(\v_1,\nablagz\v_1)-\A_4(t;\vphi_2,\nablagz\vphi_2)(\v_2,\nablagz\v_2))}_{\LQs}\\&\leq C(T)\norm{(\vphi_1-\vphi_2,\nablagz(\vphi_1-\vphi_2))}_{L^\infty(0,T;L^\infty(\Gamma_0))}\norm{\v_1}_{L^2(0,\tT;\L^\infty(\gz))}\norm{\nablagz \v_1}_{L^\infty(0,\tT;\L^2(\gz))}\\&+
    C(T)\norm{\v_1-\v_2}_{L^2(0,\tT;\mathbf{W}^{1,4}(\gz))}\norm{\v_1}_{L^\infty(0,\tT;\L^4(\gz))}\\&
    + C(T)\norm{\v_1-\v_2}_{L^\infty(0,\tT;\L^4(\gz))}\norm{\v_2}_{L^2(0,\tT;\mathbf{W}^{1,4}(\gz))}
    \\&\leq 
    C(T,R)\tT^{\thetaq}\norm{\vphi_1-\vphi_2}_{\XT^2}+C(T,R)\tT^\frac14\norm{\v_1-\v_2}_{\XT^2}
    \\&\leq C(T,\widetilde{T}, R)\Vert(\v_1 - \v_2, \vphi_1 -\vphi_2)\Vert_{\XT}.
\end{align*}
Then we have, by the embedding $\H^1(\gz)\hookrightarrow\L^4(\gz)$ and \eqref{fundamental},
\begin{align*}
 &  \Vert P(t)(\b_3(t)\v_1\cdot (\di_1(t;((\nablagz(\nablagz(\nablagz \vphi_1)_i)_j)_k))+\di_2(t; ((\nablagz(\nablagz \vphi_1)_i)_j)))\\&-\b_3(t)\v_2\cdot (\di_1(t;((\nablagz(\nablagz(\nablagz \vphi_2)_i)_j)_k))+\di_2(t; ((\nablagz(\nablagz \vphi_2)_i)_j)))\Vert_{\LQs}\\&\leq
 C(T)\norm{\v_1}_{\LLQ(4,4)}\norm{\vphi_1-\vphi_2}_{L^4(0,\tT;W^{3,4}(\gz))}+C(T)\norm{\v_1-\v_2}_{\LLQ(4,4)}\norm{\vphi_2}_{L^4(0,\tT;W^{3,4}(\gz))}\\&
 \leq C(T)\tT^\frac14\norm{\v_1}_{\LLQ(\infty,4)}\norm{\vphi_1-\vphi_2}_{L^4(0,\tT;W^{3,4}(\gz))}\\&+C(T)\tT^\frac14\norm{\v_1-\v_2}_{\LLQ(\infty,4)}\norm{\vphi_2}_{L^4(0,\tT;W^{3,4}(\gz))}
 \\&\leq C(T,\widetilde{T}, R)\Vert(\v_1 - \v_2, \vphi_1 -\vphi_2)\Vert_{\XT}.
\end{align*}
Furthermore, in a similar fashion, recalling \eqref{fundamental} and that $\XT^1\hookrightarrow L^4(0,\tT;\mathbf W^{1,4}(\gz))$, 
\begin{align*}
    &\Vert P(t)((\b_6(t)(\di_1(t;((\nablagz(\nablagz(\nablagz \vphi_1)_i)_j)_k))+\di_2(t; ((\nablagz(\nablagz \vphi_1)_i)_j)))\cdot \nablagz)\v_1\\&-(\b_6(t)(\di_1(t;((\nablagz(\nablagz(\nablagz \vphi_2)_i)_j)_k))+\di_2(t; ((\nablagz(\nablagz \vphi_2)_i)_j)))\cdot \nablagz)\v_2)\Vert_{\LQs}
    \\&\leq C(T)\norm{\nablagz\v_1}_{\LLQ(4,4)}\norm{\vphi_1-\vphi_2}_{L^4(0,\tT;W^{3,4}(\gz))}\\&+C(T)\norm{\nablagz\v_1-\nablagz\v_2}_{\LLQ(4,4)}\norm{\vphi_2}_{L^4(0,\tT;W^{3,4}(\gz))}\\&
 \leq C(T)\tT^\frac{q-2}{4q}\norm{\nabla\v_1}_{\LLQ(4,4)}\norm{\vphi_1-\vphi_2}_{\XT^2}\\&+C(T)\tT^\frac{q-2}{4q}\norm{\nabla\v_1-\nabla\v_2}_{\LLQ(4,4)}\norm{\vphi_2}_{\XT^2}
 \\&\leq C(T,\widetilde{T}, R)\Vert(\v_1 - \v_2, \vphi_1 -\vphi_2)\Vert_{\XT},
\end{align*}
where the exponent $q$ is the same as in \eqref{fundamental}. Now, we observe that
\begin{align*}
&2P^{\Gamma_0}\P_0\mathrm{div}_{\Gamma_0}(\nu(\tvphi_0)\EE_S(\v))\\&=2P^{\Gamma_0}\nu(\tvphi_0)\P_0\mathrm{div}_{\Gamma_0}(\EE_S(\v))+2P^{\Gamma_0}\nu'(\tvphi_0)\EE_S(\v)\nablagz \tvphi_0\\&
=P^{\Gamma_0}\nu(\tvphi_0)\P_0\div{\gz}(\nablagz\v)+P^{\Gamma_0}\nu(\tvphi_0)K_0\v+2P^{\Gamma_0}\nu'(\tvphi_0)\EE_S(\v)\nablagz \tvphi_0,
\end{align*}
where $K_0$ is the Gaussian curvature of $\Gamma_0$.
Above we used (see, e.g., \cite[Lemma 2.1]{V1}) that $$\P\text{div}_{\Gamma_0}(\nablagz^T\v)=\nablagz(\text{div}_{\Gamma_0}\v)+K_0\v=K_0\v,$$ recalling $\text{div}_{\Gamma_0}\v=0$.

Then, we can split the terms and compare them with the function $\a_0(t),\a_1(t;\cdot,\cdot),\a_2(t)$. Indeed, it is easy to see that 
\begin{align}
&\label{a1a}\P_0\a_0(0;\nablagz((\nablagz\v_{i})_j)_k)\v=\P_0\div{\gz}(\nablagz\v),\\&\label{a2a}
    \a_1(0;(\nablagz\v_i)_j,\nablagz\vphi)=2\EE_S(\v)\nablagz \vphi,\\&
    a_2(0)\v=K_0\v.\label{a3a}
\end{align}
Note also that, by the regularity of the flow map and of the surface, for any $t\in[0,T]$,
\begin{align*}
&\vert\partial_t\a_0(t;\nablagz((\nablagz\v_{i})_j)_k)\vert\leq C(T)\vert \div{\gz}(\nablagz \v)\vert,\\& \vert\partial_t\a_1(t;(\nablagz\v_i)_j,\nablagz\vphi)\vert\leq C(T)\vert {\nabla}_{\gz} \v\vert\vert \nablagz\vphi\vert,\\&
\vert\partial_ta_2(t)\vert\leq C(T).
\end{align*}
Therefore, concerning $\a_0$, recalling \eqref{holder}, \eqref{est1}, the continuity of $\nu,\nu'$ (recall that $\vphi_i$ are bounded in $L^\infty(0,\tT;L^\infty(\gz))$ by a constant depending on $R$), the estimate on $\partial_t\a_0$ above, and $\vphi_i(0)=\tvphi_0$ (with $\vert \tvphi_0\vert<1$ on $\gz$),
\begin{align*}
\Vert &P(t)\nu(\vphi_1)\P_0\a_0(t;\nablagz((\nablagz\v_{1,i})_j)_k)-P^{\Gamma_0}\nu(\tvphi_0)\P_0\div{\gz}(\nablagz\v_1)\\&-P(t)\nu(\vphi_2)\P_0\a_0(t;\nablagz((\nablagz\v_{2,i})_j)_k)+P^{\Gamma_0}\nu(\tvphi_0)\P_0\div{\gz}(\nablagz\v_2)\Vert_{\LQs}\\&\leq 
\norm{(P(t)-P^{\Gamma_0})(\nu(\vphi_1)-\nu(\vphi_2))\P_0\a_0(t;\nablagz((\nablagz\v_{1,i})_j)_k)}_{\LQs}
\\&+\norm{(P(t)-P^{\gz})\nu(\vphi_2)\P_0(\a_0(t;\nablagz((\nablagz\v_{1,i})_j)_k)-\a_0(t;\nablagz((\nablagz\v_{2,i})_j)_k))}_{\LQs}\\&
+\norm{P^{\gz}(\nu(\vphi_1)-\nu(\tvphi_0))\P_0(\a_0(t;\nablagz((\nablagz(\v_{1,i}-\v_{2,i}))_j)_k)}_{\LQs}\\&+
\norm{P^{\gz}\nu(\tvphi_0)\P_0(\a_0(t;\nablagz((\nablagz(\v_{1,i}-\v_{2,i}))_j)_k)-\a_0(0;\nablagz((\nablagz(\v_{1,i}-\v_{2,i}))_j)_k))}_{\LQs}
\\&\leq
C(T,R)\tT\norm{\v_1}_{L^2(0,\tT;\H^2(\gz))}\norm{\vphi_1-\vphi_2}_{L^\infty(0,\tT;L^\infty(\gz)))}\\&+C(T)\tT\norm{\nu(\vphi_2)}_{L^\infty(0,\tT;L^\infty(\gz))}\norm{\v_1-\v_2}_{L^2(0,\tT;\H^2(\gz))}\\&+C(T,R)\norm{\vphi_1-\tvphi_0}_{L^\infty(0,\tT;L^\infty(\gz)))}\norm{\v_1-\v_2}_{L^2(0,\tT;\H^2(\gz))}\\&
+C(T)\norm{\nu(\tvphi_0)}_{L^\infty(0,\tT;L^\infty(\gz)))}\norm{t\int_0^1 \partial_t\a_0(st;\nablagz((\nablagz(\v_{1,i}-\v_{2,i}))_j)_k))ds}_{\LQ}\\&\leq C(T,\widetilde{T}, R)\Vert(\v_1 - \v_2, \vphi_1 -\vphi_2)\Vert_{\XT},
\end{align*}
where again we used that, from \eqref{holder},
\begin{align}
\norm{\vphi_1-\tvphi_0}_{L^\infty(0,\tT;L^\infty(\gz)))}\leq \tT^{\thetaq}\norm{\vphi_1}_{C^{\thetaq}([0,\tT]; W^{2,p}(\gz))} \leq C(R,T)\tT^{\thetaq},
\label{est11}
\end{align}
and, from  \eqref{a1a},
\begin{align*}
    \norm{t\int_0^1 \partial_t\a_0(st;\nablagz((\nablagz(\v_{1,i}-\v_{2,i}))_j)_k))ds}_{\LQ}\leq C(T)\tT\norm{\v_1-\v_2}_{L^2(0,\tT;\H^2(\gz))}.
\end{align*}
Now, concerning $\a_1$, we have, exploiting \eqref{est1}, \eqref{est11}, $\vphi_1(0)-\vphi_2(0)=0$, the control over $\partial_t\a_1$, the regularity of the flow map, and the continuity $\nu',\nu''$, recalling that $\norm{\tvphi_i}_{L^\infty(0,\tT;L^\infty(\gz))}\leq CR$,
\begin{align*}
 & \Vert   P(t)\nu'(\vphi_1)\a_1(t;(\nablagz\v_{1,i})_j,\nablagz\vphi_1)-2P^{\gz}\nu'(\tvphi_0)\EE_S(\v_1)\nablagz^T \tvphi_0\\&
-P(t)\nu'(\vphi_2)\a_1(t;(\nablagz\v_{2,i})_j,\nablagz\vphi_2)+2P^{\gz}\nu'(\tvphi_0)\EE_S(\v_2)\nablagz^T \tvphi_0
 \Vert_{\LQs}\\&\leq
\norm{(P(t)-P^{\gz})(\nu'(\vphi_1)-\nu'(\vphi_2))\a_1(t;(\nablagz\v_{1,i})_j,\nablagz\vphi_1)}_{\LQs}\\&
\quad+\norm{(P(t)-P^{\gz})\nu'(\vphi_2)\a_1(t;(\nablagz\v_{1,i})_j,\nablagz(\vphi_1-\vphi_2))}_{\LQs}\\&\quad
+\norm{(P(t)-P^{\gz})\nu'(\vphi_2)\a_1(t;(\nablagz(\v_{1,i}-\v_{2,i}))_j,\nablagz(\vphi_2))}_{\LQs}\\&\quad
+\norm{P^{\gz}(\nu'(\vphi_1)-\nu'(\tvphi_0))\a_1(t;(\nablagz\v_{1,i})_j,\nablagz(\vphi_1-\vphi_2))}_{\LQs}\\&\quad
+\norm{P^{\gz}(\nu'(\vphi_1)-\nu'(\tvphi_0))\a_1(t;(\nablagz(\v_{1,i}-\v_{2,i}))_j,\nablagz(\vphi_2))}_{\LQs}\\&\quad
+\norm{P^{\Gamma_0}(\nu'(\vphi_1)-\nu'(\vphi_2))\a_1(t;(\nablagz\v_{2,i})_j,\nablagz(\vphi_2))}_{\LQs}
\\&\quad
+\norm{P^{\gz}\nu'(\tvphi_0)(\a_1(t;(\nablagz\v_{1,i})_j,\nablagz(\vphi_1-\tvphi_0))}_{\LQs}\\&\quad
+\norm{P^{\gz}\nu'(\tvphi_0)(\a_1(t;(\nablagz\v_{2,i})_j,\nablagz(\vphi_1-\vphi_2))}_{\LQs}\\&\quad
+\norm{P^{\gz}\nu'(\tvphi_0)(\a_1(t;(\nablagz(\v_{1,i}-\v_{2,i}))_j,\nablagz(\tvphi_0))-\a_1(0;(\nablagz(\v_{1,i}-\v_{2,i}))_j,\nablagz(\tvphi_0)))}_{\LQs}\\&
\leq C(T,R)\tT\norm{\vphi_1-\vphi_2}_{L^\infty(0,\tT;L^\infty(\gz))}\norm{\v_1}_{L^2(0,\tT;\H^1(\gz))}\norm{\nablagz\vphi_1}_{L^\infty(0,\tT;L^\infty(\gz))}\\&\quad
+C(T)\tT\norm{\nu'(\vphi_1)}_{L^\infty(0,\tT;L^\infty(\gz))}\norm{\v_1}_{L^2(0,\tT;\H^1(\gz))}\norm{\nablagz\vphi_1-\nablagz\vphi_2}_{L^\infty(0,\tT;L^\infty(\gz))}\\&\quad
+C(T)\tT\norm{\nu'(\vphi_1)}_{L^\infty(0,\tT;L^\infty(\gz))}\norm{\v_1-\v_2}_{L^2(0,\tT;\H^1(\gz))}\norm{\nablagz\vphi_2}_{L^\infty(0,\tT;L^\infty(\gz))}\\&\quad
+C(T,R)\tT^{\thetaq}\norm{\vphi_1}_{\XT^2}\norm{\v_1}_{L^2(0,\tT;\H^1(\gz))}\norm{\nablagz\vphi_1-\nablagz\vphi_2}_{L^\infty(0,\tT;L^\infty(\gz))}\\&\quad
+C(T,R)\tT^{\thetaq}\norm{\vphi_1}_{\XT^2}\norm{\v_1-\v_2}_{L^2(0,\tT;\H^1(\gz))}\norm{\nablagz\vphi_1}_{L^\infty(0,\tT;L^\infty(\gz))}\\&\quad
+C(T)\norm{\vphi_1-\vphi_2}_{L^\infty(0,\tT;L^\infty(\gz))}\norm{\v_2}_{L^2(0,\tT;\H^1(\gz))}\norm{\nablagz\vphi_2}_{L^\infty(0,\tT;L^\infty(\gz))}\\&\quad
+C(T)\norm{\nu'(\vphi_0)}_{L^\infty(0,\tT;L^\infty(\gz))}\norm{\v_1}_{L^2(0,\tT;\H^1(\gz))}\norm{\nablagz\vphi_1-\nablagz\tvphi_0}_{L^\infty(0,\tT;L^\infty(\gz))}\\&\quad
+C(T)\norm{\nu'(\vphi_0)}_{L^\infty(0,\tT;L^\infty(\gz))}\norm{\v_2}_{L^2(0,\tT;\H^1(\gz))}\norm{\nablagz\vphi_1-\nablagz\vphi_2}_{L^\infty(0,\tT;L^\infty(\gz))}
\\&\quad
+C(T)\norm{\nu'(\tvphi_0)}_{L^\infty(0,\tT;L^\infty(\gz))}\norm{t\int_0^1\partial_t\a_1(st;(\nablagz(\v_{1,i}-\v_{2,i}))_j,\nablagz\tvphi_0)ds}_{\LQ}
\\&\quad
\leq C(T,\widetilde{T}, R)\Vert(\v_1 - \v_2, \vphi_1 -\vphi_2)\Vert_{\XT},
\end{align*}
where we used again $(\vphi_1-\vphi_2,\nablagz(\vphi_1-\vphi_2))=(\vphi_1-\vphi_2-(\vphi_1-\vphi_2)(0),\nablagz(\vphi_1-\vphi_2)-\nablagz(\vphi_1-\vphi_2)(0))$ and, from \eqref{a2a},
\begin{align*}
    \norm{t\int_0^1\partial_t\a_1(st;(\nablagz\v_{i})_j,\nablagz \tvphi_0)ds}_{\LQ}\leq C(T)\tT\norm{\nablagz \tvphi_0}_{L^\infty(0,\tT;L^\infty(\gz))}\norm{\v}_{L^2(0,\tT;\H^1(\gz))}.
\end{align*}
Then, about $a_2$, the arguments are analogous, and, recalling also \eqref{a3a}, in the end we easily get
\begin{align*}
    &\norm{P(t)\nu(\vphi_1)a_2(t)\v_1-P^{\gz}\nu(\tvphi_0)K_0\v_1-P(t)\nu(\vphi_2)a_3(t)\v_2+P^{\gz}\nu(\tvphi_0)K_0\v_2}\\&\leq C(T,\widetilde{T}, R)\Vert(\v_1 - \v_2, \vphi_1 -\vphi_2)\Vert_{\XT}.
\end{align*}
We can then proceed with the analysis. We have, exploiting \eqref{fundamental} and \eqref{est0},
\begin{align*}
    &\Vert-P(t)\b_7(t)(\di_1(t;((\nablagz(\nablagz(\nablagz \vphi_1)_i)_j)_k))+\di_2(t; ((\nablagz(\nablagz \vphi_1)_i)_j)))\\&+P(t)\b_7(t)(\di_1(t;((\nablagz(\nablagz(\nablagz \vphi_2)_i)_j)_k))+\di_2(t; ((\nablagz(\nablagz \vphi_2)_i)_j)))\Vert_{\LQs}\\&\leq
    C(T)\norm{\vphi_1-\vphi_2}_{L^2(0,\tT;H^3(\gz))}\leq C(T)\tT^\frac14\norm{\vphi_1-\vphi_2}_{L^4(0,\tT;W^{3,p}(\gz))} \leq C(T,\widetilde{T}, R)\Vert(\v_1 - \v_2, \vphi_1 -\vphi_2)\Vert_{\XT}.
\end{align*}
In a similar way we have, recalling again \eqref{fundamental},
\begin{align*}
    &\Vert -P(t)\b_8(t;\nablagz\vphi_1,(\nablagz(\nablagz\vphi_1)_i)_j)+P(t)\b_8(t;\nablagz\vphi_2,(\nablagz(\nablagz\vphi_2)_i)_j)\Vert_{\LQ}\\&\leq
    C(T)(\norm{\nablagz\vphi_1-\nablagz\vphi_2}_{L^\infty(0,\tT;\L^\infty(\gz))}\norm{\vphi_1}_{L^2(0,\tT;H^2(\gz))}+\norm{\nablagz\vphi_2}_{L^\infty(0,\tT;\L^\infty(\gz))}\norm{\vphi_1-\vphi_2}_{L^2(0,\tT;H^2(\gz))})
    \\&\leq C(T)(\tT^{\frac14}\norm{\vphi_1}_{L^4(0,\tT;H^2(\gz))}\norm{\vphi_1-\vphi_2}_{\XT^2}+\tT^{\frac14}\norm{\vphi_1-\vphi_2}_{L^4(0,\tT;H^2(\gz))}\norm{\vphi_2}_{\XT^2})\\&
    \leq C(T,\widetilde{T}, R)\Vert(\v_1 - \v_2, \vphi_1 -\vphi_2)\Vert_{\XT}.
\end{align*}
We can now pass to estimate the second row of the vector $\mathcal F(\v_1,\vphi_1)-\mathcal F(\v_2,\vphi_2)$, which has to be estimated in $L^q(0,\tT;L^p(\gz))$. It is immediate to deduce, since $\H^1(\gz)\hookrightarrow\L^p(\gz)$, recalling \eqref{uniform} and the decomposition \eqref{decompositionB} (we denote by $DB$ the Jacobian matrix of $B$ with respect to $\vphi,\nablagz\vphi$, as we did for the operator $\A_\ast$),
\begin{align*}
&\norm{B(t;\vphi_1,\nablagz\vphi_1)(\v_1)-B(t;\vphi_2,\nablagz\vphi_2)(\v_2)}_{\Lqp}\\&\leq 
C(T)(\norm{DB}_{L^\infty((0,\tT)\times [-CR,CR]\times [-CR,CR]^3)}\norm{(\vphi_1-\vphi_2,\nablagz(\vphi_1-\vphi_2))}_{L^\infty(0,\tT;\L^\infty(\Gamma_0))})(1+\norm{\v_1}_{\Lqpb})
\\&+ C(T)\norm{\v_1-\v_2}_{\Lqpb}\norm{B(t;\vphi_2,\nablagz\vphi_2)}_{L^\infty(0,\tT;\L^\infty(\gz))}
    \\&\leq C(T,R)\tT^{\frac1q}(\norm{\vphi_1-\vphi_2}_{\XT^2}+\norm{\v_1-\v_2}_{L^\infty(0,\tT;\H^1(\gz))})\\&
  \leq C(T,\widetilde{T}, R)\Vert(\v_1 - \v_2, \vphi_1 -\vphi_2)\Vert_{\XT}. 
\end{align*}
Note that here we have also exploited the continuity of $\Psi$ up to its fourth derivative, as well as the fact that $\norm{\vphi_i}_{W^{1,\infty}(\Gamma_0)}\leq CR$, $i=1,2$. Recalling the decomposition \eqref{decompositionB}, we have also defined 
\begin{align*}
&\norm{B(t;\vphi_2,\nablagz\vphi_2)}_{L^\infty(0,\tT;\L^\infty(\Gamma_0))}\\&:=\sup_{(x,t)\in\Gamma_0 \times [0,T]}{\norma{B_{0}(t;\vphi_2(x,t),\nablagz\vphi_2(x,t))}}+\sup_{(x,t)\in\Gamma_0 \times [0,T]}\sup_{\u\in \R^3,\ \u\not=\mathbf 0}\frac{\norma{B_{1}(t;\vphi_2(x,t),\nablagz\vphi_2(x,t))(\u)}}{\norma{\u}},   
\end{align*} 
as well as 
\begin{align*}
&\norm{DB}_{L^\infty((0,\tT)\times [-CR,CR]\times [-CR,CR]^3)}\\&:=\sup_{(t,y,\mathbf z)\in(0,\tT)\times [-CR,CR]\times [-CR,CR]^3}{\norma{DB_{0}(t;y,\mathbf z)}}\\&+\sup_{(t,y,\mathbf z)\in(0,\tT)\times [-CR,CR]\times [-CR,CR]^3}\sup_{\u\in \R^3,\ \u\not=\mathbf 0}\frac{\norma{DB_{1}(t;y,\mathbf z)(\u)}}{\norma{\u}}.
\end{align*}

\noindent Now, concerning $f_0$
 and $\Delta^2_{\gz}\vphi$, it is immediate to observe that
 $$
f_0(0;(\nablagz(\nablagz(\nablagz(\nablagz \vphi)_i)_j)_k)_l)=\Delta_{\Gamma_0}^2\vphi,
 $$
 and
 $$
 \norma{\partial_t f_0(t;(\nablagz(\nablagz(\nablagz(\nablagz \vphi)_i)_j)_k)_l)}\leq C(T)\sum_{i,j,k,l}\vert (\nablagz(\nablagz(\nablagz(\nablagz (\vphi))_i)_j)_k)_l\vert,
 $$
 so that 
\begin{align*}
   & \Vert -f_0(t;(\nablagz(\nablagz(\nablagz(\nablagz (\vphi_1-\vphi_2))_i)_j)_k)_l)+\Delta_{\Gamma_0}^2(\vphi_1-\vphi_2)\Vert_{\Lqp}\\&\leq 
  \norm{t \int_0^1 \partial_t f_0(st;(\nablagz(\nablagz(\nablagz(\nablagz (\vphi_1-\vphi_2))_i)_j)_k)_l)ds}_{\Lqp}\\&
  \leq C(T)\tT\norm{\vphi_1-\vphi_2}_{L^q(0,\tT;W^{4,p}(\gz))}\leq C(T,\widetilde{T}, R)\Vert(\v_1 - \v_2, \vphi_1 -\vphi_2)\Vert_{\XT}.
  \end{align*}
 Furthermore, concerning $f_1$ and $f_2$, since $2<q\leq 4$ by assumption, and recalling \eqref{fundamental},
\begin{align*}
   &\Vert-f_1(t;((\nablagz(\nablagz(\nablagz \vphi_1)_i)_j)_k))+f_2(t; ((\nablagz(\nablagz \vphi_1)_i)_j)\\&+f_1(t;((\nablagz(\nablagz(\nablagz \vphi_2)_i)_j)_k))-f_2(t; ((\nablagz(\nablagz \vphi_2)_i)_j)\Vert_{\Lqp}\\&
   \leq C(T)\norm{\vphi_1-\vphi_2}_{L^q(0,\tT;W^{3,p}(\gz))}\leq C(T)\tT^{\frac{q-2}{4q}}\norm{\vphi_1-\vphi_2}_{\XT^2}\\&\leq C(T,\widetilde{T}, R)\Vert(\v_1 - \v_2, \vphi_1 -\vphi_2)\Vert_{\XT}.
\end{align*}
 In conclusion, for $f_3$ we have, similarly, recalling again \eqref{fundamental} with $q\in(2,4]$, $(\vphi_1-\vphi_2)(0)=0$, and the fact that $\Psi\in C^4(\R)$,
 \begin{align*}
        &\norm{\Psi''(\vphi_1)f_3(t;(\nablagz(\nablagz\vphi_1)_i)_j)-\Psi''(\vphi_1)f_3(t;(\nablagz(\nablagz\vphi_2)_i)_j)}_{\Lqp}\\&\leq
 C(T)\max_{r\in [-2CR,2CR]}\vert\Psi^{(iii)}(r)\vert \norm{\vphi_1-\vphi_2-(\vphi_1-\vphi_2)(0)}_{L^\infty(0,\tT;L^\infty(\gz))}\norm{\vphi_1}_{L^q(0,\tT;W^{2,p}(\gz))}\\&+
 C(T)\max_{r\in [-CR,CR]}\vert\Psi^{(ii)}(r)\vert\norm{\vphi_1-\vphi_2}_{L^q(0,\tT;W^{2,p}(\gz))}\\&
 \leq C(T,R)\tT^{\thetaq}\norm{\vphi_1-\vphi_2}_{C^{\thetaq}([0,\tT];W^{2,p}(\gz))}+C(T,R)\tT^\frac{q-2}{4q}\norm{\vphi_1-\vphi_2}_{\XT^2}\\&\leq 
 C(T,\widetilde{T}, R)\Vert(\v_1 - \v_2, \vphi_1 -\vphi_2)\Vert_{\XT}.
 \end{align*}
 Therefore, putting together all the estimates, we have shown that there is a constant $C(T,\widetilde{T}, R)>0$ such that
$$\Vert\mathcal{F}(\v_1, \vphi_1) -\mathcal{F}(\v_2, \vphi_2)\Vert_{\YT}
\leq C(T,\widetilde{T}, R)\Vert(\v_1 - \v_2, \vphi_1 -\vphi_2)\Vert_{\XT},$$
for all $(\v_i, \vphi_i ) \in \XT$ with $\Vert(\v_i, \vphi_i )\Vert_{\XT}
\leq R$, $R>0$, and $i = 1, 2$.
This concludes the proof of Proposition \ref{prop1}.
 \subsection{Proof of Proposition \ref{prop2}}
 \label{l2p}In order to prove this proposition, we can follow closely the arguments of \cite[Theorem 4]{AWe}, namely, we will make use of the auxiliary \cite[Theorem 5]{AWe} related to monotone operator theory, whose proof can be found in \cite[Chapter IV, Theorem 6.1]{Showalter}.
Let us first observe that the two rows of $\mathcal{L}$ are independent from each other, so that we can study them separately. Let us start with the velocity part. We need now to prove the following lemma, which is the surface counterpart of \cite[Lemma 3]{AWe}.
\begin{lemma}
\label{L1}
Under the same assumptions of Theorem \ref{thm1} concerning $\Gamma_0, \nu, \trho$, for every $\u_0 \in \H^1(\gz)\cap \L^2_\sigma(\gz)$, $\f \in\LQs$, $\vphi_0 \in W^{1,
\infty}(\gz)$, with $\vert \tvphi_0\vert\leq 1$ for any $x\in\gz$, and every $0 < \tilde{T} \leq T$ there exists a
unique solution $\u\in \Zut$ to the following problem
\begin{align}
    \begin{cases}
P^{\gz}\partial_t(\trho_0\u)-2P^{\gz}\P_0\mathrm{div}_{\gz}(\nu(\tvphi_0)\E_S^{\gz}(\u))+\omega\u=\f,&\quad\text{ on }\gz\times(0,\tT),\\
\mathrm{div}_{\gz}\u=0,&\quad\text{ on }\gz\times(0,\tT),\\
\u(0)=\u_0, &\quad\text{ on }\gz
    \end{cases}
    \label{p0}
\end{align}
where $\trho_0=\trho(\tvphi_0)$.
\end{lemma}
 \begin{proof}
     We aim at applying \cite[Theorem 5]{AWe}, thus we adopt the same notation, and follow the proof of \cite[Lemma 3]{AWe} for the main steps. In particular, we set $\mathcal{B}\u:=P^{\gz}(\trho_0\u)$, $E=E_b':=\L^2_\sigma(\gz)$, so that $E_b'$ is a Hilbert space endowed with the seminorm
     $$
     \vert \u\vert_b^2:=(\mathcal{B}\u)(\u)=\int_{\gz}P^{\gz}(\trho_0\u)\cdot \u\ds =\int_{\gz}\trho_0\vert \u\vert^2 \ds\cong \norm{\u}^2.
     $$
     Then clearly $E_b=\L^2_\sigma(\gz)$. Furthermore, we also introduce the operator $A:\mathfrak{D}(A)\to \L^2_\sigma(\gz)$, by 
     \begin{align*}
         (A\u)(v):=\begin{cases}
             -\int_{\gz}\left(2P^{\gz}\P_0\mathrm{div}_{\gz}(\nu(\tvphi_0)\E_S^{\gz}(\u))+\omega\u\right)\cdot \v\ds,\quad \text{if }\u\in \mathfrak{D}(A),\\
             \emptyset,\quad  \text{if }\not\in\u\in \mathfrak{D}(A),
         \end{cases}
     \end{align*}
     for any $\v\in \L^2_\sigma(\gz)$. Here we have defined $\mathfrak{D}(A):=\H^2(\gz)\cap \L^2_\sigma(\gz)$. To stick with the notation of \cite[Theorem 5]{AWe}, we then introduce the relation 
     $$
    \mathcal{A}:=\{(\u,\v):\ \v=A\u,\u\in\mathfrak{D}(A)\}\subset E\times E_b'.
     $$
    We now set $\psi: \L_\sigma^2(\gz)\to [0,+\infty]$ to satisfy
     \begin{align*}
         \psi(\u):=\begin{cases}
             \int_{\gz}\left(\nu(\tvphi_0)\E_S(\u):\E_S(\u)+\frac12\omega\vert \u\vert^2\right)\ds,\quad &\text{if }\u\in \H^1(\gz)\cap \L^2_\sigma(\gz)=:\mathfrak{D}(\psi),\\
             +\infty,\quad &\text{ else.}
         \end{cases} 
     \end{align*}
     Clearly $\psi(0)=0$ and $\u_0\in \mathfrak{D}(\psi)$ by assumption. We also notice that $\psi$ is trivially convex and lower-semicontinuos. Moreover, denoting by $\partial\psi$ the subdifferential of $\psi$, it also holds that $\mathfrak{D}(A)=\mathfrak{D}(\partial\psi)$, as well as $\partial\psi(\u)=\{A\u\}$ for any $\u\in \mathfrak{D}(A)$. The proof of these two facts can be easily adapted from the one of \cite[Lemma 4]{AWe}. In particular, this is possible since the regularity result of \cite[Lemma 5]{AWe} used in the proof is here substituted by its surface counterpart of Lemma \ref{regularity}. Indeed, the fact that $\mathfrak{D}(A)\subseteq \mathfrak{D}(\partial\psi)$ is trivial. About the opposite inclusion, let us assume $\u\in \mathfrak{D}(\partial\psi)$, and let $\w\in \partial\psi(\u)\subset\L^2_\sigma(\gz)$. This means that, for any $\v\in \L^2_\sigma(\gz)$ it holds 
     \begin{align}
         \psi(\u)\leq \psi(\v)+\int_{\gz}\w\cdot (\u-\v)\ds.
         \label{subdiff}
     \end{align}
Choosing now $\v:=\u+t\tilde{\w}$, $t>0$, for some $\tilde{\w}\in \mathfrak{D}(\psi)$, since $\v\in \mathfrak{D}(\psi)$, from \eqref{subdiff}, after some algebraic manipulations, we get
\begin{align*}
    &\psi(\u)-\psi(\v)\\&=-2t\int_{\gz}\left(\nu(\tvphi_0)\E_S(\u):\E_S(\tilde{\w})+\frac12\omega \u\cdot \tilde{\w}\right)\ds-t^2\int_{\gz}\left(\nu(\tvphi_0)\E_S(\tilde{\w}): \E_S(\tilde{\w})+\frac12\omega  \vert\tilde{\w}\vert^2\right)\ds\\&\leq -t\int_{\gz}\w\cdot\tilde{\w}\ds,
\end{align*}
    so that dividing by $-t$ and letting $t\to 0^+$ we deduce in the end, since $\tilde{\w}\in \mathfrak{D}(\psi)$ is arbitrary, that
\begin{align}
    \int_{\gz}\left(2\nu(\tvphi_0)\E_S(\u):\E_S(\tilde{\w})+\omega \u\cdot \tilde{\w}\right)\ds=\int_{\gz}\w\cdot \tilde{\w}\ds,\quad \forall \tilde{\w}\in \H^1(\gz)\cap \L^2_\sigma(\gz).
    \label{identity}
\end{align}
But by Lemma \ref{regularity} we know that there exists a unique solution $\u\in \mathfrak{D}(A)$ such that \eqref{identity} holds. Therefore, integrating by parts, we also get     
    \begin{align*}
\int_{\gz}\left(-2P^{\Gamma_0}\text{div}_{\gz}(\nu(\tvphi_0)\E_S(\u))+\omega \u\right)\cdot \tilde{\w}\ds=\int_{\gz}\w\cdot \tilde{\w}\ds,\quad \forall \tilde{\w}\in \H^1(\gz)\cap \L^2_\sigma(\gz),
\end{align*}
and by the density of $\H^1(\gz)\cap \L^2_\sigma(\gz)$ in $\L^2_\sigma(\gz)$ we immediately infer that $\w=-2P^{\Gamma_0}\text{div}_{\gz}(\nu(\tvphi_0)\E_S(\u))+\omega \u$, so that $\u\in \mathfrak{D}(A)$ and
$\w=A\u$, allowing to conclude that $\mathfrak{D}(\partial\psi)\subseteq\mathfrak{D}(A)$ and therefore $\partial\psi(\u)=\{A\u\}$ for any $\u\in \mathfrak{D}(\partial\psi)$. 

Now we have all the ingredients to apply \cite[Theorem 5]{AWe}. In particular, we obtain the existence of $\u\in W^{1,2}(0,\tT;\L^2_\sigma(\gz))$, wth $\u(t)\in \H^2(\gz)\cap \L^2_\sigma(\gz)$ for almost any $t\in(0,\tT)$ and any $0<\tT\leq T$ satisfying \eqref{p0}. To be precise, concerning initial datum, from the statement of \cite[Theorem 5]{AWe} we only obtain that $\mathcal{B}\u(0)=\mathcal{B}\u_0$, but this entails
$$
\int_{\gz} \trho_0\u(0)\cdot \v\ds=\int_{\gz} \trho_0\u_0\cdot \v\ds,\quad \forall \v\in \L^2_\sigma(\gz),
$$
from which, since $\trho_0\geq \rho_*>0$ (recall that $\tvphi_0\in[-1,1]$ by assumption), we immediately get $\u(0)=\u_0$. To conclude the proof of Lemma \ref{L1} (existence part), we need to show that $\u\in \Zut$, i.e., we need additionally $\u\in L^2(0,\tT;\H^2(\gz))$. This is immediate, since $\u$ solves 
$$
-2P^{\gz}\P_0\mathrm{div}_{\gz}(\nu(\tvphi_0)\E_S^{\gz}(\u))+\omega\u=\f-P^{\gz}\partial_t(\trho_0\u)=:\f_1\in L^2(0,\tT;\L^2_\sigma(\gz)),
$$
so that the result comes directly from an application of Lemma \ref{regularity}, namely from \eqref{estStokes}, recalling $\tvphi_0\in W^{1,\infty}(\gz)$. Concerning the uniqueness of solutions to \eqref{p0}, the result is trivial and can be obtained, for instance, by considering two different solutions $\u_1,\u_2$, take the difference on the equations they solve, test the resulting equation by $\u_1-\u_2$ and recall that $\nu\geq \nu_*$ by assumption, and $\sqrt{\omega\norm{\cdot}^2+2\nu_*\norm{\E_S(\cdot)}^2}$ is an equivalent norm on $\H^1(\gz)$. See also the proof of \cite[Lemma 3]{AWe} for similar arguments. This finishes the proof.
 \end{proof}
We can now pass to discuss the invertibility second row of the operator $\mathcal{L}$ defined in \eqref{L}. In this case, the result directly comes from the fact that the operator $\Delta_{\gz}^2+\omega$ has, by assumption, the $L^q$-$L^p$ maximal regularity property (see, for instance, \cite[Teorem 6.4.3 (ii)]{SP}): given an initial datum in a suitable trace space, i.e., $\tvphi_0\in (L^p(\gz),W^{4,p}(\gz))_{1-\frac1q,q }$, and a function $f\in L^q(0,\tT;L^p(\gz))$ (for $p>4$, $q\in(2,4]$ given in the definition of $\Zdt$, there exists a unique $\vphi\in \Zdt$ such that
\begin{align*}
    \begin{cases}
        \partial_t\vphi+\Delta_{\gz}^2\vphi=f,&\quad\text{ on }\gz\times(0,\tT),\\
        \vphi(0)=\tvphi_0,&\quad\text{ on }\gz.
    \end{cases}
\end{align*}
To sum up this result with the one of Lemma \ref{L1}, we have obtained, by the bounded inverse theorem, that, for any $0<\tT\leq T$ the operator $\mathcal{L}$ in \eqref{L} is invertible as a linear operator from $\YT$ to $\XT$, i.e., there exists $C(\tT)>0$ at this stage possibly depending on $\tT$, such that 
$$
\norm{\mathcal{L}^{-1}}_{\mathcal{L}(\YT,\XT)}\leq C(\tT),\quad \forall 0<\tT\leq T.
$$
To conclude the proof of Proposition \ref{prop2} it is then enough to show that the constant above does not actually change with $\tT$. This can be obtained by a simple extension argument (see for instance the proof of \cite[Lemma 7]{AWe}), leading to show that
$$
\norm{\mathcal{L}^{-1}}_{\mathcal{L}(\YT,\XT)} \leq \norm{\mathcal{L}^{-1}}_{\mathcal{L}(Y_T,X_T)}\leq C(T),\quad \forall 0<\tT\leq T,
$$
concluding the proof of the proposition.
\section{Proof of Theorem \ref{strong1a}}
\label{proof1}
In order to prove the theorem, we consider the equivalent problem \eqref{detract} and look for $(\V,\pi,\vphi)$ which satisfies the equations. 
In order to obtain back a solution to \eqref{mainp} it is enough to set $\v:=\V+\widehat{\u}$, since, as already observed, thanks to the regularity of the flow map, $\widehat{\u}\in C^4_{\W^{3,p}}$ for any $p\in[2,+\infty)$.
Now observe that, by assumption, the initial data are such that 
$$
\V_0=\v_0-\widehat{\u}(0)\in \H^1_\sigma(\gz),\quad  \vphi_0\in B^{4-\frac 4q}_{p,q}(\gz),$$
for some $p>4$ and $q\in(2,4]$.
Therefore, the initial data $\u_0=\V_0$ and $\tvphi_0=\vphi_0$, and the functions $\Psi,\nu,\rho$ are sufficiently regular to apply Theorem \ref{thm1}. This means that there exists $0<\tT\leq T$ and $(\u,\tvphi)\in X_{\tT}$ satisfying problem \eqref{systf}. Furthermore, if 
\begin{align}
   \exists \delta_0\in(0,1):\ \norm{\vphi_0}_{L^\infty(\gz)}\leq 1-2\delta_0,
   \label{delta}
\end{align}
then there exists $T_*\in(0,\tT]$ such that also \eqref{strspr} holds.
Now we can set $\V=\phi_t(\u)$ and $\vphi=\phi_t(\tvphi)$, and, thanks to the regularity of the flow map, we immediately deduce that $(\V,\vphi)$ satisfies problem \eqref{detract} and is  such that 
\begin{align*}
    \V\in  L^2_{\H^2\cap \L^2_\sigma(\tT)},\quad \dt \V\in L^2_{\L^2(\tT)},\quad \vphi\in L^q_{W^{4,p}(\tT)},\quad \dt\vphi\in L^q_{L^p(\tT)},
\end{align*}
and, if \eqref{delta} holds, then
\begin{align}
\sup_{t\in[0,T_*]}\norm{\vphi(t)}_{L^\infty\gt}\leq 1-\delta_0.
\label{separation}
\end{align}
Observe that the regularity on $\dt\V\in L^2_{\L^2(\tT)}$ can be obtained as follows: since $\phi_{-t}(\dts\V)=\dtn \u\in L^2(0,\tT;\L^2_\sigma(\gz))$, by the regularity of the flow map we deduce $\dts\V\in L^2_{\L^2_\sigma(\tT)}$, then, since it holds (see \eqref{relation_base1}) $\dt\V=\dts\V-\A(\dt\A^{-1})\V$, we deduce that $\dt\V\in L^2_{\L^2}$. Concerning $\dt\vphi$ we simply have $\dtn\tvphi=\tphimt(\dt\vphi)$ and thus the regularity can be simply deduced by the regularity of the flow map and of $\partial_t\tvphi$.

In conclusion, by repeating the same proof as in \cite[Section 7]{ES} (see also \cite[Section 4.1]{V2}), we infer that there exists $\pi\in L^2_{L^2(\tT)}$ such that, in the sense of distributions, 
\begin{align*}
    \nablag \pi=\f&:=-\rho \P\dt \v -((\rho\v+\J_\rho)\cdot\nabla_\Gamma)\v-\rho \Vn\mathbf{H}\v-\Vn\H\J_\rho+2\P\divg(\nu(\vphi)\E_S(\v))\\&\quad 
     -\P\divg(\nablag\vphi\otimes\nablag\vphi)+2\P\divg(\nu(\vphi)\Vn\mathbf{H})+\frac\rho2\nablag(\Vn)^2.
\end{align*}
It is then immediate to deduce that ${\f}\in{L^2_{\L^2(\tT)}}$, so that we infer $\pi\in L^2_{H^1(\tT)}$. This concludes the proof.

\appendix
\section{On the regularity for a Laplace equation on evolving surfaces}
In this section we aim to find higher-order time regularity results to the problem
\begin{align}
-\Delta_\Gamma g(t)=f(t),\quad \overline{f}(t)\equiv0,\quad\text{ on }\Gamma(t),\quad \forall t\in[0,T],
\label{laplac}
\end{align}
with $\overline{g}(t)\equiv 0$ for any $t\in[0,T]$.
Some preliminary results have already been obtained in \cite[Lemma A.1]{ES}. In particular, it is shown that, under suitable regularity assumptions on the flow map (which are satisfied also by our assumptions in Section \ref{regflowmap1}), if $f\in C^0_{W^{1,p}}\cap C^1_{L^{p}}$, for any $p\geq 2$, then also the unique solution $g$ is such that $g\in C^0_{W^{3,p}}\cap C^1_{W^{1,p}}$. Here we actually prove the following regularity result.
\begin{lemma}
Assume that the flow map $\Phi^n_t$ enjoys the regularity stated in Section \ref{regflowmap1}. Then, if $f\in C^1_{W^{2,p}}$ for any $p\geq2$, the unique solution $g$ to \eqref{laplac} is such that $g\in C^1_{W^{4,p}}$. If we further define $\tilde{\u}:=\nablag g$, then it holds $\tilde{\u}\in C^{0}_{\W^{3,p}}$ and $\P\dt\tilde{\u}\in C^{0}_{\W^{3,p}}$ for any $p\geq2$.
\label{appendixb}
\end{lemma}
\begin{proof}
    Existence and uniqueness of a solution $g$ to \eqref{laplac} can be easily obtained by an application of Lax-Milgram Lemma.
    Then, by elliptic regularity, it is immediate to deduce that 
$$
\norm{g}_{C^0_{W^{4,p}}}\leq C(T)\norm{f}_{C^0_{W^{2,p}}},\quad \forall p\geq 2.
$$
 Notice that the constant $C$  depends on $T$ and not on the specific $t\in[0,T]$. This can be shown  exploiting the elliptic regularity theory applied on each surface $\gam(t)$ as well as the regularity on the flow map, in a similar way as done in the proof of Lemma \ref{regularity}. Let us consider the regularity in time. From \cite[Lemma A.1]{ES} it holds $g\in C^1_{W^{1,p}}$. We can now write the problem in weak formulation as: for any $\eta\in H^1(\Gamma(t))$,
\begin{align*}
\ints{\Gamma}\nablag g\cdot \nablag \eta\ds=\ints{\Gamma}f\eta\ds, \quad \forall t\in[0,T].
\end{align*}
Let us consider now test functions $\eta\in C^1_{H^1}$.
 By taking the time derivative of this formulation (which is rigorous, thanks to the regularity of $g$), we obtain, recalling \cite[Proposition 2.8]{DE},
 \begin{align}
  &\nonumber\ints{\Gamma}\nablag \dt g\cdot \nablag \eta \ds+   \ints{\Gamma}\nablag g\cdot \nablag \dt\eta \ds+\ints{\Gamma}\nablag g\cdot \nablag \eta \divg\VVn\ds-\ints{\Gamma}\nablag  g\cdot {\E}_s(\VVn)\nablag \eta \ds\\&
=\ints{\Gamma}\dt f\eta\ds+\ints{\Gamma}\ f\dt\eta\ds+ \ints{\Gamma}f\eta\divg\VVn\ds.
 \end{align}
 Since then $\dt\eta\in C^0_{H^1}$ is still an admissible test function, we deduce that it holds in the end, after integrating by parts,
 \begin{align}
  &\nonumber\ints{\Gamma}\nablag \dt g\cdot \nablag \eta \ds\\&
=-\ints{\Gamma}\divg(\divg\VVn\nablag g )\eta\ds-\ints{\Gamma}\divg({\E}_s(\VVn)\nablag  g) \eta \ds+\ints{\Gamma}\dt f\eta\ds+ \ints{\Gamma}f\eta\divg\VVn\ds.
 \end{align}
 This entails that, formally, $\dt g$ satisfies the problem
 $$
 -\Delta_\Gamma \dt g=-\divg(\divg\VVn\nablag g )-\divg({\E}_s(\VVn)\nablag  g)+\dt f+f\divg\VVn,\quad \overline{\dt g}(t)=\ints{\Gamma}g\divg\VVn\ds,
 $$
 so that, again by elliptic regularity theory uniform over $[0,T]$, we can immediately deduce, since the right-hand side belongs to $C^0_{W^{2,p}}$, that $\dt g\in C^{0}_{W^{4,p}}$ for any $p\geq 2$. In conclusion, setting $\tilde{\u}:=\nablag g$, it is immediate to deduce $g\in C^0_{\W^{3,p}}$ for any $p\geq 2$. Moreover, we can observe by standard computations (see, e.g., \cite[Lemma 2.6]{Dzk}) that
 $$
 \P\dt\tilde{\u}=\nablag\dt g-\nablag^T\VVn\nablag g,
 $$
 entailing from the above regularity that $\P\dt\tilde{\u}\in C^{0}_{\W^{3,p}}$ for any $p\geq2$, thus concluding the proof.
\end{proof}
\section{Computations for the pullback equations}
\label{pullback1}
In this appendix section we show all the precise computations leading to system \eqref{syst1}. Note that in the following we will use the Einstein summation convention. Let us start from the equation for the velocity.
 First we have
 \begin{align}
 (\nablag \V)_{ij}=(\nabla_\Gamma \tphi_t(\A_{lr}))_j \tphi_t(\u_r)\P_{il}+\tphi_t(\A_{ir}(\nablagz\u_r)_p)\Dm_{pj}. \end{align}
 This entails
 \begin{align*}
     \divg(\nablag\V)_i&=[\nablag((\nabla_\Gamma \tphi_t(\A_{lr}))_j\tphit(\P_{il}))]_j\tphi_t(\u_r)+(\nabla_\Gamma \tphi_t(\A_{lr}))_j\tphi_t(\P_{il})\tphit(\nablagz\u_r)_k\Dm_{kj}\\&+[\nablag(\tphi_t(\A_{ir})\Dm_{pj})]_j\tphi_t(\nablagz\u_r)_p
     +\tphit(\A_{ir}(\nablagz(\nablagz\u_r)_p)_q)\Dm_{qj}\Dm_{pj}.
 \end{align*}
 Therefore, we end up with
 \begin{align*}
     \phmt(\P\divg(\nablag\V))_q&=(\P_0\tphimt(\A^{-1}))_{qi}\tphimt([\nablag((\nabla_\Gamma \tphi_t(\A_{lr}))_j\tphit(\P_{il}))]_j)\u_r\\&
      +(\P_0\tphimt(\A^{-1}))_{qi}(\nabla_\Gamma \tphi_t(\A_{lr}))_j\tphi_t(\P_{il})\tphit(\nablagz\u_r)_k\tphimt(\Dm_{kl})
     \\&+(\P_0\tphimt(\A^{-1}))_{qi}\tphimt([\nablag(\tphi_t(\A_{ir})\Dm_{pj})]_j)(\nablagz\u_r)_p
    \\&+(\P_0\tphimt(\A^{-1}))_{qi}\A_{lr}(\nablagz(\nablagz\u_r)_p)_k\tphimt(\Dm_{kj}\Dm_{pj}).
 \end{align*}
 In the end, we can write 
 \begin{align*}
\phmt(\P\divg(\nu(\vphi)\nablag\V))_q&=\phmt(\nu'(\vphi)\nablag\V\nablag\vphi)+\phmt(\nu(\vphi)\P\divg(\nablag\V))\\&
=\nu'(\tvphi)(\nablag \tphit(\A_{lr}))_j\u_r\P_{il}(\nablagz\tvphi)_k\tphimt(\Dm_{kj}\A^{-1}_{qi})\\&
\quad + \nu'(\tvphi)\A_{ir}(\nablagz \u_r)_p(\nablagz\tvphi)_k\tphimt(\Dm_{pj}\Dm_{kj}\A^{-1}_{qi})
\\& \quad+\nu(\tvphi)(\P_0\tphimt(\A^{-1}))_{qi}\tphimt([\nablag((\nabla_\Gamma \tphi_t(\A_{lr}))_j\tphit(\P_{il}))]_j)\u_r\\&
      \quad+\nu(\tvphi)(\P_0\tphimt(\A^{-1}))_{qi}(\nabla_\Gamma \tphi_t(\A_{lr}))_j\tphi_t(\P_{il})(\nablagz\u_r)_k\tphimt(\Dm_{kl})
     \\& \quad+\nu(\tvphi)(\P_0\tphimt(\A^{-1}))_{qi}\tphimt([\nablag(\tphi_t(\A_{ir}\Dm_{pj})]_j)(\nablagz\u_r)_p
    \\& \quad+\nu(\tvphi)(\P_0\tphimt(\A^{-1}))_{qi}\A_{lr}(\nablagz(\nablagz\u_r)_p)_k\tphimt(\Dm_{kj}\Dm_{pj}).
 \end{align*}
Analogously, one can find, recalling that (see, for instance, \cite[Lemma 2.1]{V1}) $\P\divg(\nablag^T\V)=\nablag(\divg \V)+K\V=K\V$, with $K$ as the Gaussian curvature, that 
 \begin{align*}
\phmt(\P\divg(\nablag^T\V))=
\tphimt(K)\u,
 \end{align*}
so that in the end
\begin{align*}
    &\phmt(\P\divg(\nu(\vphi)\nablag^T\V))_q=\phmt(\nu'(\vphi)\nablag^T\V\nablag\vphi)+\nu(\tvphi)\phmt(\P\divg(\nablag^T\V))\\&=
    \nu'(\tvphi)(\nablag \tphit(\A_{lr}))_i\u_r\P_{lj}(\nablagz\tvphi)_k\tphimt(\Dm_{ki}\A^{-1}_{qj})\\&
\quad + \nu'(\tvphi)\A_{jr}(\nablagz \u_p)_r(\nablagz\tvphi)_k\tphimt(\Dm_{pi}\Dm_{ki}\A^{-1}_{qj})\\&\quad+\nu(\tvphi)\tphimt(K)\u_q.
\end{align*}
Concerning now the transport terms, one has
\begin{align*}
    &\phmt(\rho(\V\cdot \nablag)\V)_q=\phmt((\nablag\V)\rho\V)_q\\&=\tilde{\rho}(\P_0)_{qr}\tphimt(\A^{-1}_{rl}(\nablag \tphit(\A_{lk}))_j)\u_k\A_{ji}\u_i+\tilde{\rho}(\tphimt(\Dm)\A\u\cdot\nablagz)\u.
\end{align*}
 Similarly, we can obtain
 \begin{align*}
    &\phmt(((\J_\rho+\rho\widehat{\u})\cdot \nablag)\V)=\phmt((\nablag\V)(\J_\rho+\rho\widehat{\u}))\\&=(\P_0)_{qr}\tphimt(\A^{-1}_{rl}(\nablag \tphit(\A_{lk}))_j)\u_k\tphimt((\J_\rho+\rho\widehat{\u}))_j+(\tphimt(\Dm(\J_\rho+\rho\widehat{\u}))\cdot\nablagz)\u.
\end{align*}
We need now to explicitate $\J_\rho$. To this aim, we have 
$$
\tilde{\mu}=\Psi'(\tvphi)-\tphimt(\Delta_\Gamma\vphi).
$$
Now, it holds
\begin{align}
\tphimt(\Delta_\Gamma\vphi)=(\nablagz(\nablagz\tvphi)_m)_k\tphimt(\Dm\D^{-T})_{mk}+(\nablagz\tvphi)_m\tphimt(\nablag \Dm_{ms})_s.
    \label{deltaphi}
\end{align}
This means that
\begin{align*}
    \nonumber\tphimt(\nablag\mu)_j&=\Psi''(\tvphi)(\nablagz\tvphi)_l\tphimt(\Dm)_{lj}\\&\quad-(\nablagz(\nablagz(\nablagz \tvphi)_m)_k)_l\tphimt(\Dm_{lj}\Dm_{mq}\Dm_{kq})\\&\quad
    -(\nablagz(\nablagz\tvphi)_m))_k\tphimt((\nablag(\Dm\D^{-T})_{mk})_j)\\&\quad -(\nablagz(\nablagz\tvphi)_m)_l\tphimt(\Dm_{lj}(\nablag\Dm_{ms})_s)\\&\quad 
    -(\nablagz\tvphi)_m\tphimt((\nablag(\nablag \Dm_{ms})_s)_j),
\end{align*}
which gives the explicit representation of $\tilde{\J}_\rho=-\frac{\widetilde{\rho}_1-\widetilde{\rho}_2}2\tphimt(\nablag\mu)$.
Proceeding in the computations,
In conclusion, the only nontrivial term left in the pullback is the Korteweg force. In particular, we obtain  
\begin{align*}
\phmt(\P\divg(\nablag\vphi\otimes \nablag \vphi))_{j}&=(\P_0)_{jm}\tphimt(\A^{-1}_{mi})(\nablagz\tvphi)_p(\nablagz(\nablagz\tvphi)_r)_k\tphimt(\Dm_{ri}\Dm_{kl}\Dm_{pl})\\&\quad+
(\P_0)_{jm}\tphimt(\A^{-1}_{mi})(\nablagz\tvphi)_p(\nablagz\tvphi)_r\tphimt(\Dm_{ri}(\nablag\Dm_{pl})_l)\\&\quad+
(\P_0)_{jm}\tphimt(\A^{-1}_{mi})(\nablagz\tvphi)_r(\nablagz(\nablagz \tvphi)_p)_k\tphimt(\Dm_{ri}\Dm_{pl}\Dm_{kl})\\&\quad
+(\P_0)_{jm}\tphimt(\A^{-1}_{mi})(\nablagz\tvphi)_r(\nablagz\tvphi)_p\tphimt(\Dm_{pl}(\nablag(\Dm_{ri})_l)).
\end{align*}
We can now pass to consider the pullback Cahn--Hilliard equation. Namely, we have 
\begin{align*}
\tphimt(\Delta^2_\Gamma\vphi)&=(\nablagz(\nablagz(\nablagz(\nablagz \tvphi)_m)_k)_r)_p\tphimt(\Dm_{pi}\Dm_{ri}\Dm_{mq}\Dm_{kq})\\&\quad+
(\nablagz(\nablagz(\nablagz\tvphi)_m)_k)_r\tphimt((\nablag(\Dm_{mq}\Dm_{kq}\Dm_{rl}))_l)\\&\quad+
(\nablagz(\nablagz(\nablagz\tvphi)_m)_r)_s)\tphimt(\Dm_{si}\Dm_{ri}(\nablag(\Dm_{mj}))_j)\\&\quad+
(\nablagz(\nablagz\tvphi)_m)_r\tphimt((\nablag(\Dm_{rl}(\nablag\Dm_{mj})_j))_l)\\&\quad
+(\nablagz(\nablagz(\nablagz(\tvphi)_m)_k)_s\tphimt(\Dm_{si}(\nablag(\Dm_{mq}\Dm_{kq}))_i)\\&
\quad+
(\nablagz(\nablagz\tvphi)_m)_k\tphimt((\nablag(((\nablag(\Dm_{mq}\Dm_{kq}))_l))_l)\\&
\quad +(\nablagz (\nablagz \tvphi)_m)_k\tphimt(\Dm_{kl}(\nablag(\nablag\Dm_{ms})_s)_l)\\&
 \quad + (\nablagz \tvphi)_m\tphimt((\nablag (\nablag(\nablag\Dm_{ms})_l)_l).
\end{align*}
Moreover, we have
\begin{align*}
     \tphimt(\Delta_\Gamma \Psi'(\vphi))&=\Psi'''(\tvphi)\nablagz^T\tvphi\tphimt(\Dm(\Dm)^{T})\nablagz\tvphi\\&\quad+ \Psi''(\tvphi)(\nablagz(\nablagz\tvphi)_i)_m\tphimt(\D_{ml}\D_{il})\\&
     \quad+\Psi''(\tvphi)(\nablagz\tvphi)_i\tphimt((\nablag \Dm_{ij})_j).
\end{align*}
Furthermore, the advective term becomes
\begin{align*}
    \tphimt(\V\cdot \nablag \vphi)=\A\u\cdot (\nablagz \tvphi\tphimt(\Dm)).
\end{align*}
These computations, recalling \eqref{dtv2}-\eqref{dtv1} and the definition of $\nablagphi$ in \eqref{press}, lead to the following system \eqref{syst1bis}:
\begin{align}
\label{syst1bis}
    \begin{cases}
\widetilde{\rho}\partial_t\u+\widetilde{\rho}\b_0(t)\u+\widetilde{\rho}\b_1(t)\u\cdot\u+\widetilde{\rho}(\b_2(t)\u\cdot\nablagz)\u\\
+\b_3(t)\u\cdot \widetilde{\J}_\rho+\widetilde{\rho}\b_4(t)\u+\widetilde{\rho}(\b_5(t)\cdot\nablagz)\u+(\b_6(t)\widetilde{\J}_\rho\cdot \nablagz)\u
 \\-\nu(\tvphi)\P_0\a_0(t;\nablagz((\nablagz\u_i)_j)_k)-\nu'(\tvphi)\a_{01}(t)\u\cdot\nablagz\tvphi-\nu'(\tvphi)\a_{02}(t;(\nablagz\u_i)_j,\nablagz\tvphi)
\\-\nu(\tvphi)\a_{03}(t;(\nablagz\u_i)_j)-\nu(\tvphi)\a_1(t)\u-\nu(\tvphi)a_2(t)\u+\widetilde{\rho}\b_7(t)\u+\b_7(t)\widetilde{\J}_\rho
+\b_8(t;\nablagz\tvphi,(\nablagz(\nablagz\tvphi)_i)_j)\\+\b_{9}(t)\nablagz\tvphi\cdot \nablagz\tvphi
+\widetilde{\rho}\b_{10}(t)+\b_{11}(t)\nablagz\tvphi+\nu(\tvphi)\b_{12}(t)+\widetilde{\rho}\b_{13}(t)+\nabla_{\Gamma(t)}^\phi\widetilde{\pi}=\mathbf{0},\\\, \\
\mathrm{div}_{\Gamma_0}\u=0,\\ \, \\
\partial_t\tvphi+\c_0(t)\u\cdot \nablagz\tvphi+\c_1(t)\cdot \nablagz\tvphi\\+f_0(t;(\nablagz(\nablagz(\nablagz(\nablagz \tvphi)_i)_j)_k)_l)\\+f_1(t;((\nablagz(\nablagz(\nablagz \tvphi)_i)_j)_k))+f_2(t; ((\nablagz(\nablagz \tvphi)_i)_j)+\c_2(t)\cdot\nablagz \tvphi\\
-\Psi''(\tvphi)f_3(t;(\nablagz(\nablagz\tvphi)_i)_j)-\Psi'''(\tvphi)\c_3(t)\nablagz\tvphi\cdot \nablagz\tvphi
-\Psi''(\tvphi)\c_4(t) \nablagz\tvphi=0,\\ \\
\widetilde{\mu}=\Psi'(\tvphi)-c_5(t;((\nablagz(\nablagz \tvphi)_i)_j)-\c_6(t)\nablagz\tvphi,\\ \\
\widetilde{\J}_\rho=\Psi''(\tvphi)\di_0(t)\nablagz\tvphi+\di_1(t;((\nablagz(\nablagz(\nablagz \tvphi)_i)_j)_k))+\di_2(t; ((\nablagz(\nablagz \tvphi)_i)_j))+\di_3(t)\nablagz\tvphi,
    \end{cases}
\end{align}
where all the functions $\mathbf a_i,\mathbf b_i,\mathbf c_i,\mathbf d_i,f_i$ are (bi)linear in their arguments (apart from time $t$). Here the boldface stands for a $3\times 3\times 3$ tensor (or possibly a $3\times 3$ matrix or even a vector in $\R^3$. With the notation $\b\v_1\cdot \v_2$, if $\b$ is a three-dimensional tensor, by its application to vectors $\v_1,\v_2$ we mean, in coordinates,$(\b\v_1\cdot \v_2)^k=\b_{ij}^k\v_{1,i}\v_{2,j}$, for $k=1,2,3$.  Moreover, when we indicate, for instance, $d_3(t; ((\nablagz(\nablagz \tvphi)_i)_j)$, we mean that $d_3$ depends linearly on \textit{all} $(\nablagz(\nablagz \tvphi)_i)_j$, for $i,j=1,2,3$.

Namely, we can characterize precisely all the functions introduced so far. We have:
\begin{align*}
&(\a_0)_q:=(\tphimt(\A^{-1}))_{qi}\A_{lr}(\nablagz(\nablagz\u_r)_p)_k\tphimt(\Dm_{kj}\Dm_{pj}),\\& (\a_{01})_q:=(\P_0\tphimt(\A^{-1}))_{qi}\tphimt([\nablag(\tphi_t(\A_{ir}\Dm_{pj})]_j)(\nablagz\u_r)_p\\&\quad\quad\quad\quad+(\P_0\tphimt(\A^{-1}))_{qi}(\nabla_\Gamma \tphi_t(\A_{lr}))_j\tphi_t(\P_{il})(\nablagz\u_r)_k\tphimt(\Dm_{kl}),\\&
(\a_{02})_{qr}:=(\P_0\tphimt(\A^{-1}))_{qi}\tphimt([\nablag((\nabla_\Gamma \tphi_t(\A_{lr}))_j\tphit(\P_{il}))]_j),\\&
(\a_{03})_{kr}^q:=(\nablag \tphit(\A_{lr}))_j\P_{il}\tphimt(\Dm_{kj}\A^{-1}_{qi})),\\&
(\a_1)_q:=\A_{ir}((\nablagz \u_r)_p+(\nablagz \u_p)_r)(\nablagz\tvphi)_k\tphimt(\Dm_{pj}\Dm_{kj}\A^{-1}_{qi}),\\&
a_2=\tphimt(K),
    \\&\b_0:= \P_0\tphi_{-t}(\dt(\A^{-1}))\A,\quad 
    (\b_1)_{ki}^q:=(\P_0)_{qr}\tphimt(\A^{-1}_{rl}(\nablag \tphit(\A_{lk}))_j)\A_{ji},\\&
\b_2:=\tphimt(\Dm)\A,\quad (\b_3)_{jk}^q= (\P_0)_{qr}\tphimt(\A^{-1}_{rl}(\nablag \tphit(\A_{lk}))_j),\\&
  (\b_4)_{qk}= (\P_0)_{qr}\tphimt(\A^{-1}_{rl}(\nablag \tphit(\A_{lk}))_j)\tphimt(\widehat{\u}),\quad
  \b_5:=\tphimt(\Dm\widehat{\u}),\quad \b_6:=\tphimt(\Dm),\\&
  \b_7:=\tphimt(v_\n\H),\quad (\b_8)^j:=(\P_0)_{jm}\tphimt(\A^{-1}_{mi})(\nablagz\tvphi)_p(\nablagz(\nablagz\tvphi)_r)_k\tphimt(\Dm_{ri}\Dm_{kl}\Dm_{pl}),\\&(\b_9)_{pr}^j:=(\P_0)_{jm}\tphimt(\A^{-1}_{mi})\tphimt(\Dm_{ri}(\nablag\Dm_{pl})_l),\quad \b_{10}:= \tphimt( \A^{-1}(\P\dt \widehat{\u} +(\widehat{\u}\cdot\nabla_\Gamma)\widehat{\u}+ \Vn\mathbf{H}\widehat{\u})),\\&
  \b_{11}:=-2\tphimt(\A^{-1}(v_\n\H+\E_s(\widehat{\u}))(\Dm)^{T}),\\&
\b_{12}:=-2\tphimt(\A^{-1}\P\divg(v_\n\H+\E_s(\widehat{\u}))),\quad \b_{13}=\frac12\tphimt(\A^{-1}\nablag v_\n^2),\\&
  f_0:=(\nablagz(\nablagz(\nablagz(\nablagz \tvphi)_m)_k)_r)_p\tphimt(\Dm_{pi}\Dm_{ri}\Dm_{mq}\Dm_{kq}),\\&
  f_1:=(\nablagz(\nablagz(\nablagz\tvphi)_m)_k)_r\tphimt((\nablag(\Dm_{mq}\Dm_{kq}\Dm_{rl}))_l)\\&\quad\quad+
(\nablagz(\nablagz(\nablagz\tvphi)_m)_r)_s)\tphimt(\Dm_{si}\Dm_{ri}(\nablag(\Dm_{mj}))_j)\\&\quad\quad+(\nablagz(\nablagz(\nablagz(\tvphi)_m)_k)_s\tphimt(\Dm_{si}(\nablag(\Dm_{mq}\Dm_{kq}))_i),
\\& f_2:=(\nablagz(\nablagz\tvphi)_m)_k\tphimt((\nablag(((\nablag(\Dm_{mq}\Dm_{kq}))_l))_l)+(\nablagz (\nablagz \tvphi)_m)_k\tphimt(\Dm_{kl}(\nablag(\nablag\Dm_{ms})_s)_l),\\& f_3:=(\nablagz(\nablagz\tvphi)_i)_m\tphimt(\D_{ml}\D_{il}),\\&
 \c_0:=\tphimt(\D^{-T})\A,\quad \c_1:=\tphimt(\D^{-T}(\widehat{\u}+\V_\n)),\quad
 (\c_2)_m:= (\nablagz \tvphi)_m\tphimt((\nablag (\nablag(\nablag\Dm_{ms})_l)_l),\\&
\c_3:=\tphimt(\Dm\D^{-T}),\quad (\c_4)^i:=\tphimt((\nablag \Dm_{ij})_j),\\&
c_5:=(\nablagz(\nablagz\tvphi)_m)_k\tphimt(\Dm\D^{-T})_{mk},\quad (\c_6)_m:=\tphimt(\nablag \Dm_{ms})_s,\\&
\di_0:=-\frac{\widetilde{\rho}_1-\widetilde{\rho}_2}{2}\tphimt(\D^{-T}),\quad (\di_1)_j:=\frac{\widetilde{\rho}_1-\widetilde{\rho}_2}{2}((\nablagz(\nablagz(\nablagz \tvphi)_m)_k)_l\tphimt(\Dm_{lj}\Dm_{mq}\Dm_{kq})),\\&
(\di_2)_j:=\frac{\widetilde{\rho}_1-\widetilde{\rho}_2}{2}((\nablagz(\nablagz\tvphi)_m))_k\tphimt((\nablag(\Dm\D^{-T})_{mk})_j)+(\nablagz(\nablagz\tvphi)_m)_l\tphimt(\Dm_{lj}(\nablag\Dm_{ms})_s)),\\&
(\di_3)_{jm}:=\frac{\widetilde{\rho}_1-\widetilde{\rho}_2}{2}\tphimt((\nablag(\nablag \Dm_{ms})_s)_j).
\end{align*}

Now we can summarize some of the terms, for the sake of clarity: first, we set $\A_\ast:[0,T]\times \R\times \R^3\times \R^3\times \R^9\to \R^3$ as

\begin{align*}
\A_\ast(t;\tvphi,\nablagz\tvphi)(\u,\nablagz\u)&:= \widetilde{\rho}\b_0(t)\u+\widetilde{\rho}\b_1(t)\u\cdot\u+\widetilde{\rho}(\b_2(t)\u\cdot\nablagz)\u
\\&+\b_3(t)\u\cdot(\Psi''(\tvphi)\di_0(t)\nablagz\tvphi+\di_3(t)\nablagz\tvphi)
+\widetilde{\rho}\b_4(t)\u+\widetilde{\rho}(\b_5(t)\cdot\nablagz)\u\\&+(\b_6(t)(\Psi''(\tvphi)\di_0(t)\nablagz\tvphi+\di_3(t)\nablagz\tvphi)\cdot \nablagz)\u\\&+\b_7(t)(\Psi''(\tvphi)\di_0(t)\nablagz\tvphi+\di_3(t)\nablagz\tvphi)\\&+\b_{9}(t)\nablagz\tvphi\cdot \nablagz\tvphi+
\widetilde{\rho}\b_{10}(t)+\b_{11}(t)\nablagz\tvphi+\nu(\tvphi)\b_{12}(t)+\widetilde{\rho}\b_{13}(t)\\&
-\nu(\tvphi)\a_1(t;(\nablagz\u_i)_j)-\nu(\tvphi)\a_2(t)\u.
\end{align*}
Then we also introduce $B:[0,T]\times \R\times \R^3\times \R^3\to \R$ as
\begin{align*}
     B(t;\tvphi,\nablagz\tvphi)(\u)&:=\c_0(t)\u\cdot \nablagz\tvphi+\c_1(t)\cdot \nablagz\tvphi+\c_2(t)\cdot\nablagz\tvphi\\&-\Psi'''(\tvphi)\c_3(t)\nablagz\tvphi\cdot \nablagz\tvphi-\Psi''(\tvphi)\c_4(t) \nablagz\tvphi,
\end{align*}
and notice that $\A_\ast$ can be split into the sum of five parts: 
\begin{align}
\nonumber&\A_\ast(t;\tvphi,\nablagz\tvphi)(\u,\nablagz\u)\\&\nonumber=\A_0(t;\tvphi,\nablagz\tvphi)+\A_1(t;\tvphi,\nablagz\tvphi)(\u)+\A_2(t;\tvphi,\nablagz\tvphi)(\nablagz\u)\\&+\A_3(t;\tvphi,\nablagz\tvphi)(\u,\u)+\A_4(t;\tvphi,\nablagz\tvphi)(\u,\nablagz\u),\label{decomposition}
\end{align}
which are linear (or bilinear) in the arguments $\u,\nablagz\u$ for $\A_i$, $i=0,1,2,4$, and quadratic in $\u$ for $\A_3$. Analogously, for $B$ we can write
\begin{align}
B(t;\tvphi,\nablagz\tvphi)(\u)=B_0(t;\tvphi,\nablagz\tvphi)+B_1(t;\tvphi,\nablagz\tvphi)(\u),   \label{decompositionB}
\end{align}
where $B_1$ is linear in $\u$. Substituting this decomposition in \eqref{syst1bis}, together with the definitions of $\widetilde{\mu}$ and $\widetilde{\J}_{\rho}$, we end up with system \eqref{syst1}.
\medskip

\textbf{Acknowledgments.} \rev{The authors thank the anonymous referee for the appropriate comments and useful remarks.} Part of this work was done while AP was visiting HA and HG at the Department of Mathematics of the University of Regensburg, whose hospitality is kindly acknowledged. The authors gratefully acknowledge the support by the Graduiertenkolleg 2339 IntComSin of the Deutsche
Forschungsgemeinschaft (DFG, German Research Foundation) – Project-ID 321821685. 
AP is a member of Gruppo Nazionale per l’Analisi Matematica, la Probabilità e le loro Applicazioni (GNAMPA) of
Istituto Nazionale per l’Alta Matematica (INdAM).  This research was funded in part by the Austrian Science Fund (FWF) \href{https://doi.org/10.55776/ESP552}{10.55776/ESP552}.
For open access purposes, the authors have applied a CC BY public copyright license to
any author accepted manuscript version arising from this submission.\\

\textbf{Conflict of interests.} On behalf of all authors, the corresponding author states that there is no conflict of interest.\\

\textbf{Data availability.} Data sharing not applicable to this article as no datasets were generated or analysed during
the current study.

\bibliography{Bibliography}
\bibliographystyle{abbrv}

\end{document}